\RequirePackage[l2tabu, orthodox]{nag}

\documentclass[a4paper,noamsfonts]{amsart}

\usepackage[T1]{fontenc}
\usepackage[tt/osf=false]{cfr-lm}
\usepackage[british]{babel}
\usepackage{hieroglf}

\usepackage{amsmath}
\usepackage{amsthm}
\usepackage{amssymb}
\usepackage[foot]{amsaddr}
\usepackage{mathtools}

\usepackage{tikz}
\usetikzlibrary{cd}
\usetikzlibrary{decorations.markings}

\usepackage{multirow}
\usepackage[dvipsnames]{xcolor}
\usepackage{graphicx}
\usepackage{booktabs}
\usepackage{pinlabel}

\usepackage{subcaption}

\DeclareCaptionLabelFormat{xx}{#2}
\usepackage{enumitem}
\setlist[enumerate]{format=\lstyle\textup}
\usepackage[backend=biber, style=numeric, maxnames=8, isbn=false, sortcites=true, backref]{biblatex}
\AtEveryBibitem{\clearlist{language}}
\usepackage{hyperref}
\hypersetup{
    colorlinks,
    linkcolor={red!50!black},
    citecolor={blue!50!black},
    urlcolor={blue!80!black}
}

\usepackage{zref-clever}
\usepackage{microtype}

\zcsetup{noabbrev, cap}
\zcsetup{countertype={thm=theorem}}
\newtheorem{thm}{Theorem}
\numberwithin{thm}{section}
\newtheorem*{thm*}{Theorem}
\newtheorem*{mainthm*}{Main Theorem}
\AddToHook{env/lem/begin}{\zcsetup{countertype={thm=lemma}}}
\newtheorem{lem}[thm]{Lemma}
\AddToHook{env/prp/begin}{\zcsetup{countertype={thm=proposition}}}
\newtheorem{prp}[thm]{Proposition}
\AddToHook{env/cor/begin}{\zcsetup{countertype={thm=corollary}}}
\newtheorem{cor}[thm]{Corollary}
\zcRefTypeSetup{conjecture}{
Name-sg = Conjecture ,
name-sg = conjecture ,
Name-pl = Conjectures ,
name-pl = conjectures ,
}
\AddToHook{env/conj/begin}{\zcsetup{countertype={thm=conjecture}}}

\zcRefTypeSetup{question}{
Name-sg = Question ,
name-sg = question ,
Name-pl = Questions ,
name-pl = questions ,
}
\AddToHook{env/question/begin}{\zcsetup{countertype={thm=question}}}
\newtheorem{question}[thm]{Question}
\newtheorem{subclaim}{Subclaim}
\zcRefTypeSetup{subclaim}{
Name-sg = Subclaim ,
name-sg = subclaim ,
Name-pl = Subclaims ,
name-pl = subclaims ,
}

\theoremstyle{definition}
\AddToHook{env/defn/begin}{\zcsetup{countertype={thm=definition}}}
\newtheorem{defn}[thm]{Definition}
\AddToHook{env/ex/begin}{\zcsetup{countertype={thm=example}}}
\newtheorem{ex}[thm]{Example}
\zcRefTypeSetup{cons}{
Name-sg = Construction ,
name-sg = construction ,
Name-pl = Constructions ,
name-pl = constructions ,
}
\AddToHook{env/cons/begin}{\zcsetup{countertype={thm=cons}}}
\newtheorem{cons}[thm]{Construction}
\newtheorem{prob}[thm]{Problem}
\zcRefTypeSetup{prob}{
Name-sg = Problem ,
name-sg = problem ,
Name-pl = Problems ,
name-pl = problems ,
}
\AddToHook{env/prob/begin}{\zcsetup{countertype={thm=prob}}}

\theoremstyle{remark}
\AddToHook{env/rem/begin}{\zcsetup{countertype={thm=remark}}}
\newtheorem{rem}[thm]{Remark}
\AddToHook{env/notation/begin}{\zcsetup{countertype={thm=notation}}}

\DeclarePairedDelimiter{\abs}{\lvert}{\rvert}

\newcommand{\df}{\textit}
\newcommand{\union}{\cup}
\newcommand{\inter}{\cap}

\newcommand{\mc}{\mathcal}

\newcommand{\mf}{\mathfrak}

\newcommand{\IN}{\mathbb{N}}
\newcommand{\IZ}{\mathbb{Z}}
\newcommand{\IR}{\mathbb{R}}

\newcommand{\IC}{\mathbb{C}}

\newcommand{\IS}{\mathbb{S}}
\newcommand{\IQ}{\mathbb{Q}}
\newcommand{\IH}{\mathbb{H}}
\newcommand{\PSL}{\mathsf{PSL}}
\newcommand{\GL}{\mathsf{GL}}
\newcommand{\SL}{\mathsf{SL}}
\newcommand{\Farey}{\mathrm{SW}}
\newcommand{\Jorgensen}{\mathrm{TW}}
\newcommand{\Goer}{\mathfrak{G}}
\newcommand{\KK}{\mathcal{K}}
\newcommand{\TP}{\mathcal{T\!P}}
\newcommand{\TA}{\mathcal{T\!A}}
\newcommand{\Tw}{\mathfrak{T}}

\DeclareMathOperator{\ceil}{ceil}

\DeclareMathOperator{\Hol}{Hol}
\DeclareMathOperator{\Hom}{Hom}
\DeclareMathOperator{\Teich}{Teich}

\DeclareMathOperator{\Aut}{Aut}

\DeclareMathOperator{\tr}{tr}

\DeclareMathOperator{\Mod}{Mod}
\DeclareMathOperator{\Out}{Out}
\DeclareMathOperator{\QH}{QH}
\DeclareMathOperator{\Inn}{Inn}
\DeclareMathOperator{\PML}{\mathcal{PML}}
\DeclareMathOperator{\hconv}{conv_{\IH^3}}

\begin{document}

\title{Kleinian groups with tunnel number one}
\author[A. Elzenaar]{Alex Elzenaar}
\address{School of Mathematics, Monash University, Melbourne}
\email{alexander.elzenaar@monash.edu}
\thanks{The author was supported by an Australian Government Research Training Program Scholarship during the period that this work was undertaken. Part of this research was performed
in Jan--Mar 2026 while the author was visiting SLMath which is supported by NSF Grant No. DMS-2424139. This document was prepared without the use of any generative AI}


\begin{abstract}
  We give a structural theory for the holonomy groups of hyperbolic links with tunnel number $1$ in arbitrary $3$-manifolds. This is a first step towards
  understanding the canonical triangulations of such links, as it generalises the combinatorial structure theory of two-bridge link groups that is an input
  to the Sakuma--Weeks triangulation.
\end{abstract}

\maketitle

{\small\tableofcontents}

\section{Introduction}
\begin{figure}
  \labellist
  \small\hair 2pt
  \pinlabel {$\tau$} [br] at 254 432
  \endlabellist
  \centering
  \includegraphics[width=\textwidth]{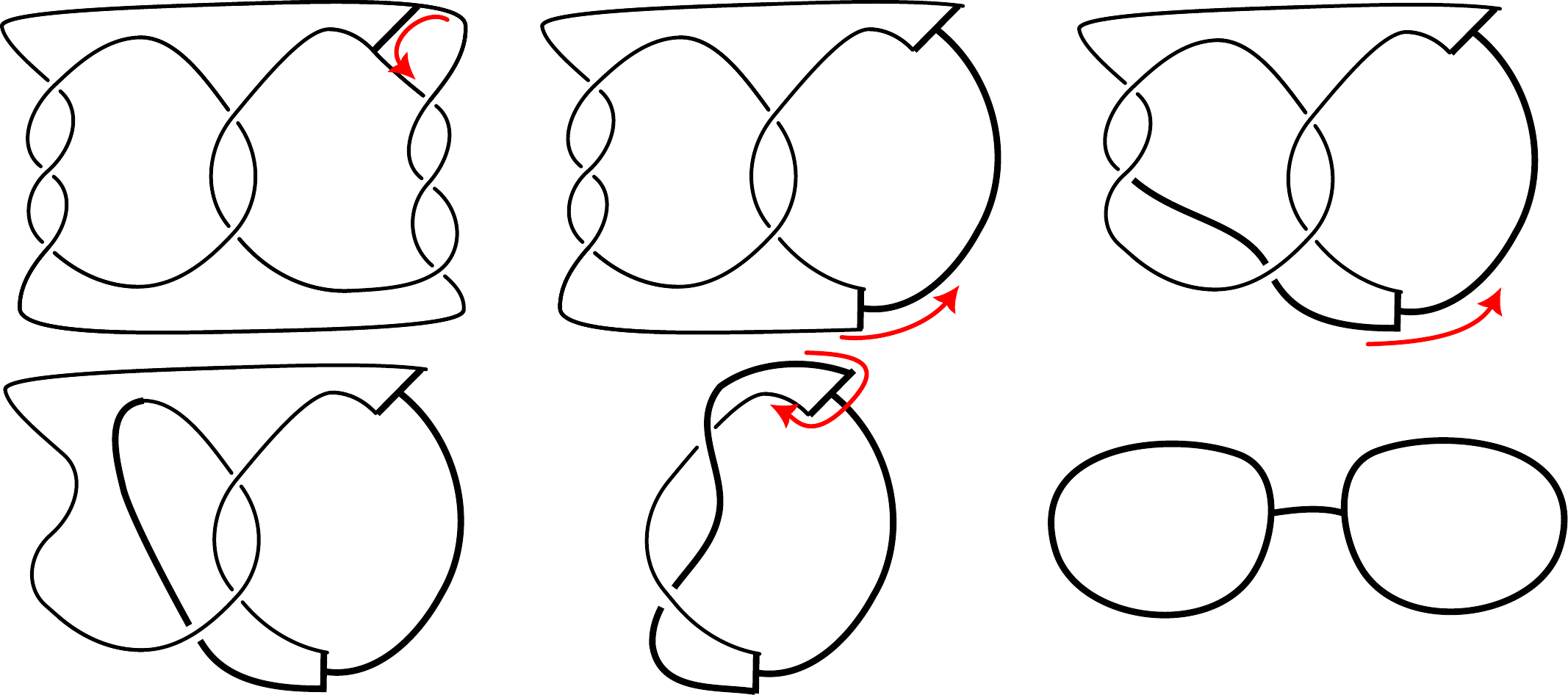}
  \caption{The \href{https://katlas.org/wiki/8_5}{$ 8_5 $ knot} $\mf{k}$ admits an unknotting tunnel $\tau$; the complement of $ \mf{k} \cup \tau $ can be untwisted via an ambient
  isotopy of $ \IS^3 $. The red arrows are supposed to depict `handle slides', where one of the three `arms' adjacent to a vertex is slid along one of the adjacent arms and around
  a handle in order to untwist part of the diagram.\label{fig:eight_five_unknotting}}
\end{figure}

A tunnel number $1$ link is a link $\mf{k} $ which admits an \df{unknotting tunnel}: an arc $ \tau $ disjoint from $ \mf{k} $ apart from its endpoints such that the complement of $ \mf{k} \union \tau $ can
be deformed by isotopy into an unknotted genus $2$ handlebody~\cite{sakuma98}. In \zcref{fig:eight_five_unknotting} we show an example of such an isotopy.
Several natural infinite families of hyperbolic $3$-manifolds are tunnel number $1$ link complements in different manifolds, for instance all two-bridge links are tunnel number $1$ links
in the $3$-sphere (\zcref{ex:2bridge}) and the simplest punctured torus bundles over the circle are tunnel number $1$ links in lens spaces (\zcref{ex:punc_torus_bund}).

Tunnel number $1$ knots are known to have many interesting geometric and topological properties, for instance tunnel number $1$ knots are always prime~\cite{norwood82} (this
result is false for links but the non-prime tunnel number $1$ links are heavily restricted~\cite{morimoto94,gordon95}), essentially never fibre over the circle~\cite{dunfield06},
and satisfy the Berge conjecture in $ \IS^3 $~\cite{li25berge}. There is also extensive modern work on the interaction between an unknotting tunnel (which is a topological object) and the hyperbolic
metric on a hyperbolic tunnel number $1$ link in $ \IS^3 $~\cite{lackenby14,cooper10,burton14,dang19}, often motivated by a conjecture of Adams~\cite{adams95} that
every unknotting tunnel in a hyperbolic manifold is isotopic to a geodesic arc.

Tunnel number $1$ links are intimately related to genus $2$ Heegaard splittings of $3$-manifolds, and in this paper we make use of this to give a structure theory for their
holonomy groups. There are two main pieces of work which we draw upon for motivation here.

Firstly, a classic paper by Sakuma and Weeks~\cite{sakuma95} gives a combinatorial
structure theory for every $2$-bridge link: a `filtration' by embedded four-punctured spheres that separate the manifold into two `caps' and finitely many copies of $ S_{0,4} \times [0,1] $
such that the entire complexity of the link is encoded in the gluing data between these pieces. This gluing data is essentially given by a path $P$ in the Farey complex from $ \infty $
to the rational number associated to the $2$-bridge link by Schubert and Conway: this point of view is taken by Purcell~\cite[Chapter~10]{purcell} who encodes the entire structure
explicitly in terms of a word made up of $R$'s and $L$'s that encodes $P$, and by Chesebro~\cite{chesebro25} who does the same thing but encodes everything in terms of `funnels' (which
are the same as the combinatorial objects called `pivot sequences' arising in the theory of once-punctured torus groups due to Minsky~\cite{minsky99} and the theory of four-punctured sphere groups due to
Ohshika and Miyachi~\cite{ohshika10}, a further motivation for our work here).

A second classic paper, this time by Cho and McCullough~\cite{cho09}, studies all tunnel number $1$ links in the $3$-sphere, and parameterises them by a tree arising from the curve complex
of the genus $2$ surface. This parameterisation reduces to the parameterisation of $2$-bridge links by the Farey tree. One of the most celebrated parts of this paper is the classification
of all tunnel number $1$ links modulo the action of the Goeritz group (roughly, modulo trivial symmetries). However, it also gives a systematic process for `building' an arbitrary tunnel
number $1$ link in $ \IS^3 $ from combinatorial data only.

In this paper, we give a combinatorial structure parameterising all tunnel number $1$ links in \emph{arbitrary} closed $3$-manifolds,
generalising the work of Cho and McCullough. This parameterisation comes from a filtration generalising that of Sakuma and Weeks for two-bridge links,
thereby giving a topological structure theory for tunnel number $1$ links in arbitrary $3$-manifolds, exhibiting them as a union of topologically simple
`slices' together with gluing data that can be interpreted topologically. For general $3$-manifolds the structure of the Goeritz group is very badly understood,
and so we are unable to generalise the full power of Cho and McCullough's work to the general setting: our parameterisation `sees' the same link multiple times. However,
this is not an obstruction to the main question of interest to us here:

\vspace{.5ex}
\begin{center}
\begin{minipage}{0.9\textwidth}\itshape\centering
Given a tunnel number $1$ link, what is its essential combinatorial structure? That is, can one find a parameterisation of all tunnel number $1$ links by finite combinatorial data
such that the link can be reconstructed explicitly, in a topologically meaningful way, from the given data?
\end{minipage}
\end{center}
\vspace{.5ex}

We note that `a link diagram', while finite and combinatorial, is not a great answer to this question since an arbitrary link diagram will not give a tunnel number $1$ link---link diagrams
do not form a natural parameter space.

To state our answer to the question more formally, if $ \mf{k} \subset A $ is a hyperbolic tunnel number $1$ link in some closed $3$-manifold $A$ then we give a filtration of groups
\begin{equation}\label{eq:filtration}
  \Gamma_0 < \Gamma_1 < \cdots < \Gamma_n = \pi_1(A \setminus \mf{k})
\end{equation}
which reflects a filtration of the link into a stack of submanifolds separated by genus two surfaces. The generators of the groups $ \Gamma_i $ can be computed
via a word-combinatorial recursion based on the topological structure. Unfortunately, the effectiveness of this for computation is reduced by a dependence on the combinatorics
of the complex of curves of the genus two surface, and this complex is difficult to compute with. Even in the case of links in the sphere studied by Cho and McCullough~\cite{cho09}
practical calculations were limited to torus links and two-bridge links for effectively the same reason (see \S\S14--15 of Cho and McCullough's paper and Ishihara~\cite{ishihara11}).
Similar computational issues have been discussed by Johannson~\cite{johannson95}.

Though it is not essential for large swathes of the theory, we will normally restrict ourselves to \emph{hyperbolic} tunnel number $1$ manifolds, because part of our motivation is to understand the Kleinian groups associated
to tunnel number $1$ links. The following version of our main result is phrased with this restriction.
\begin{mainthm*}[{\zcref{thm:tno_filtration}}]
  Let $ M $ be any hyperbolic tunnel number $1$ manifold. There exists an explicit filtration of the holonomy group of $ M $
  \begin{displaymath}
     \Gamma_0 < \cdots < \Gamma_n = \Hol(M) < \PSL(2,\IC)
  \end{displaymath}
  where each $ \Gamma_i $ is a $ S_{2,0}$-surface representation and where $ \Gamma_0 $ is a quotient of $ (\IZ \oplus \IZ)*(\IZ \oplus \IZ) $ by killing zero or more
  generators. Each $ \Gamma_i $ is the homotopy group of an embedded genus $2$ pleated surface, and
  these surfaces are nested (possibly not properly so). The filtration is `essentially unique' in the sense that all possible filtrations arising from our machinery lie within a
  uniformly bounded distance in $ \Teich(S_{2,0}) $ of each other.
\end{mainthm*}

\begin{figure}\vspace{1.5ex}
  \labellist
  \small\hair 2pt
  \pinlabel {\color{red}$\mu$} [t] at 32 211
  \pinlabel {\color{red}$\lambda$} [b] at 33 83
  \pinlabel {\color{red}$X$} [l] at 569 199
  \pinlabel {\color{red}$Y$} [l] at 569 104
  \pinlabel {generating words get shorter} [b] at 277 261
  \pinlabel {Heegaard diagrams get more complex} [t] at 277 4
  \pinlabel {\parbox{8em}{\centering rational\\braid}} [t] at 68 50
  \pinlabel {\parbox{8em}{\centering Dehn\\filling}} [t] at 147 50
  \pinlabel {$\cdots$} [t] at 228 50
  \pinlabel {$\cdots$} [t] at 318 50
  \pinlabel {$\cdots$} [t] at 398 50
  \pinlabel {\parbox{8em}{\centering rational\\braid}} [t] at 479 50
  \pinlabel {$\Gamma_n = \langle X,Y \rangle $} [l] at 503 262
  \pinlabel {$\Gamma_0 = \IZ \oplus \IZ = \langle \mu, \lambda \rangle $} [r] at 79 262
  \endlabellist
  \centering
  \includegraphics[width=\textwidth]{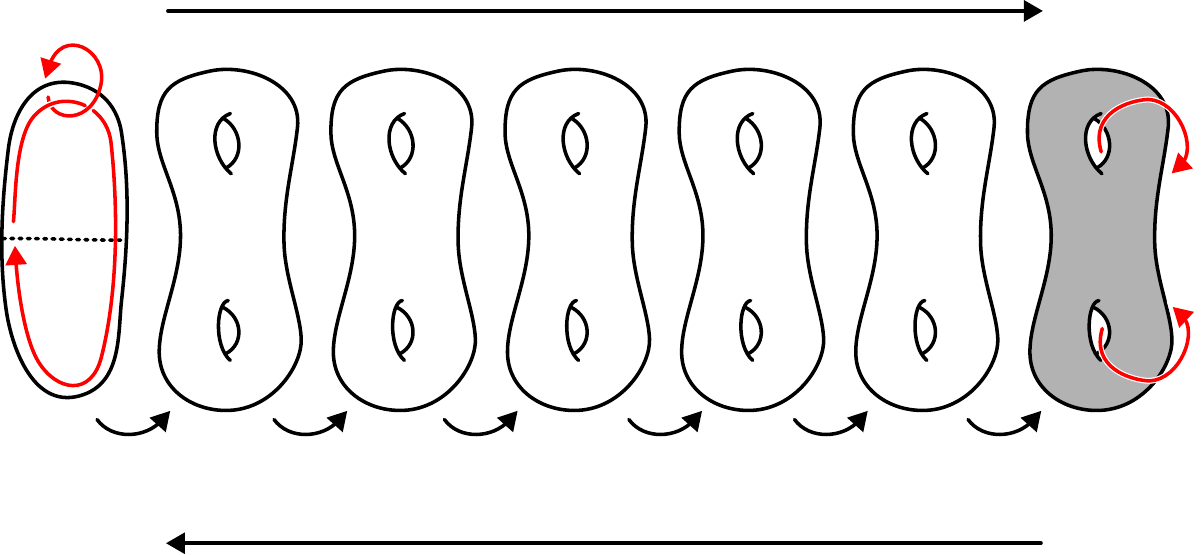}\vspace{1.5ex}
  \caption{The internal structure of a tunnel number $1$ knot.\label{fig:stack}}
\end{figure}
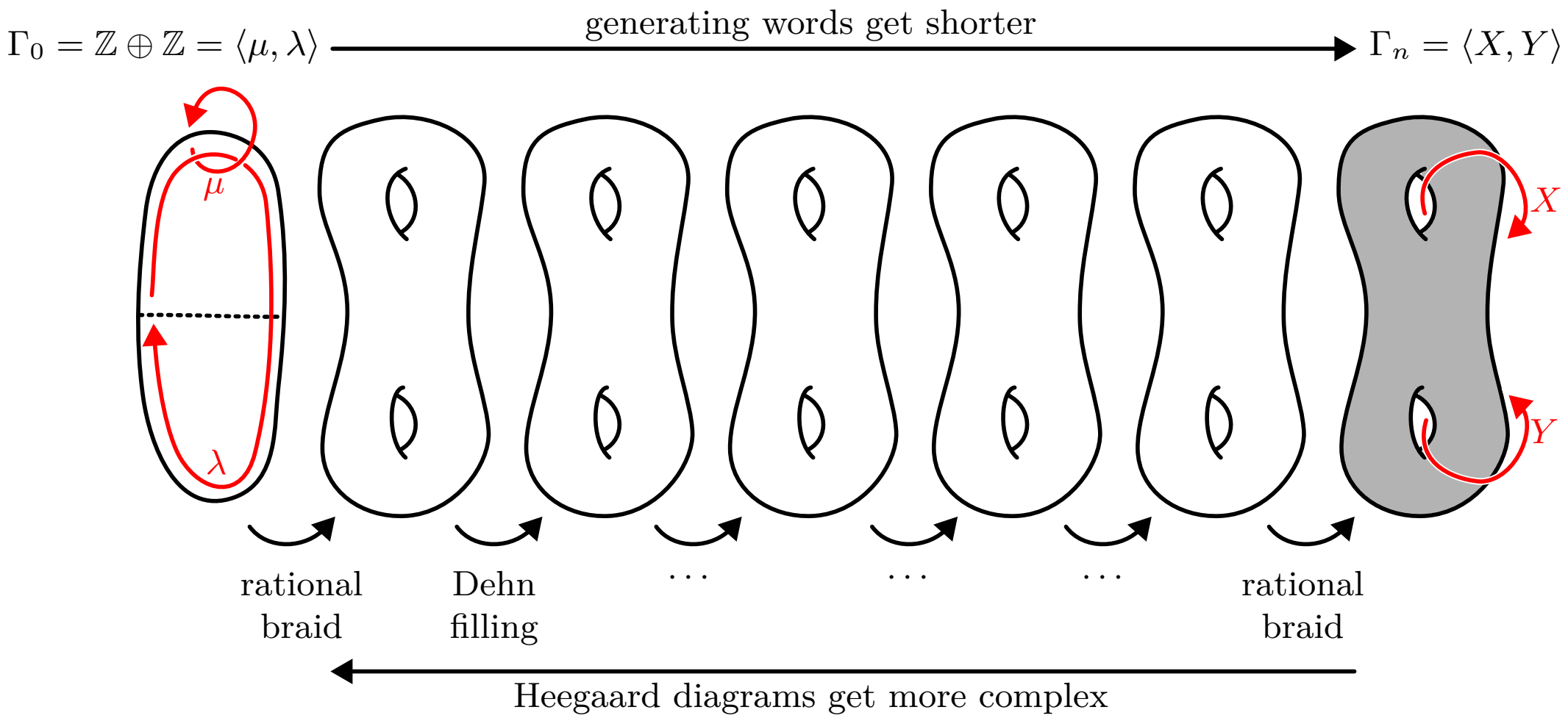

As we already alluded to, the group filtration and surfaces in the main theorem give rise to a topological structure theory for tunnel number $1$ link groups, which we explain in
\zcref{sec:tno_in_sphere} and \zcref{sec:tno_general}. In rough language, every tunnel number $1$ link is built from an alternating sequence of (i) Dehn
fillings along tori intersecting the link in two points and which change the topology of the ambient $3$-manifold, and (ii) certain `rational braiding' moves,
akin to the moves used to build $2$-bridge links, which take place inside a fixed ambient $3$-manifold. We show a schematic version of this in \zcref{fig:stack};
the idea is that the link complement is obtained by gluing thickened genus two surfaces in a stack according to combinatorially simple surface homeomorphisms,
and then to glue a solid genus $2$ handlebody onto one end and a trivial knot or separated two-component link drilled from a genus $2$ handlebody onto the other end.
It may be helpful to refer back to this picture while reading the paper.

\subsection{Motivations and directions for future work}
We have two main motivations for this study. First is the following problem~\cite{elzenaar25h,elzenaar26pr}:
\begin{prob}\label{prob:general}
  Let $ \Sigma $ be a genus $2$ surface bounding a genus $2$ handlebody $ \mc{H} $. Let $ \{\alpha,\beta,\gamma\} $ be a $2$-simplex in the curve complex of $ \Sigma $ so that
  none of the three curves are homotopically trivial in $ \mc{H} $. Classify the discrete, non-elementary, non-Fuchsian representations $ \rho : \pi_1(\mc{H}) \to \PSL(2,\IC) $ such that
  \begin{equation}\label{eq:paralocus}
    \tr^2 \rho(\alpha) = \tr^2 \rho(\beta) = \tr^2 \rho(\gamma) = 4.
  \end{equation}
\end{prob}
The natural conjecture is that every $ \rho $ satisfying \eqref{eq:paralocus} together with the additional assumptions
\begin{enumerate}[label=(\alph*)]
  \item $ \rho(\gamma) = 1 $, and
  \item both $ \rho(\alpha) $ and $ \rho(\beta) $ are parabolic
\end{enumerate}
is obtained by gluing a $2$-handle to $ \mc{H} $ along $ \gamma $, or in other words is the holonomy representation of a tunnel number $1$ link
complement. More generally, it is natural to conjecture whether every $ \rho $ satisfying \eqref{eq:paralocus} and the single additional condition
\begin{enumerate}[resume,label=(\alph*)]
  \item$ \{\alpha,\beta,\gamma\} \inter \ker \rho \neq \emptyset$
\end{enumerate}
is a Dehn filling of a tunnel number $1$ link complement. One possible approach to the proof of these conjectures would be to build an algebraic filtration for $ \rho(\pi_1(\mc{H})) $
matching that given by \eqref{eq:filtration} for some tunnel number $1$ link and then to show that the corresponding geometric pieces of each manifold are the same. The
main technical problem here that would need to be addressed by future work is the construction of appropriate disjoint pleated surfaces in the unknown
manifold $ \IH^3/\rho(\pi_1(\mc{H})) $.

Our second motivation is an ongoing project to construct fundamental domains for holonomy groups of link complements from combinatorial data. There are several classes of links where
this is possible, for instance two-bridge links~\cite{sakuma95,chesebro25} and fully augmented links~\cite[Chapter~7]{purcell}. These classes all have in common that they admit various
embedded geodesic surfaces along which they can be cut to produce simpler manifolds, and the general fundamental domain arises from combining fundamental domains for the simpler
pieces. Our results here provide a topological decomposition for tunnel number $1$ manifolds that can be used to translate between the algebra of their holonomy groups
and the combinatorics of their link diagrams. Again, a major problem for future work is to understand whether the groups $ \Gamma_i $ preserve nested pleated surfaces with rational pleating
lamination. In this direction we note that Sakai and Sakuma~\cite{sakai24} have shown that subgroups of geometrically finite Kleinian groups generated by two parabolic elements are
geometrically finite, which can be viewed as the simplest version of the setup in the current paper (c.f.\ \zcref{sec:two_bridge}).

We also want to mention a further application; this is not so practical and its main interest is that it places our work in some historical context within the
study of knot group representations. In the early 1970s, Riley set up a number of computer experiments to study the position of representations of two-bridge links in the deformation
space of subgroups of $ \PSL(2,\IC) $ generated by a pair of parabolic elements. 
Since every such subgroup is conjugate to one of the form
\begin{displaymath}
  \Big\langle X = \begin{bmatrix} 1 & 1 \\ 0 & 1 \end{bmatrix}\!,\; Y_z = \begin{bmatrix} 1 & 0 \\ z & 1 \end{bmatrix} \Big\rangle
\end{displaymath}
this is a one parameter space indexed by $ z \in \IC $. These experiments played a major role
in Thurston's initial conjectures on existence of hyperbolic structures on knot complements~\cite{riley13}. Essentially Riley's idea was to construct a system of
words $W_{p/q} $ in $ X $ and $ Y_z $ indexed by slopes of two-bridge links $ p/q \in \IQ $, so that the value of $ z $ corresponding to the two-bridge link $ p/q $ is obtained
as a root of the polynomial which is in the top-left entry of $ W_{p/q} $~\cite[Theorem~2]{riley72}. A modern point of view on this has been given by Chesebro~\cite{chesebro25}.
Our results in the current paper are of a similar, but far less practical, nature. Namely, one can choose a $\IC^3$-parameter space of the character variety of representations of the free
group into $ \PSL(2,\IC) $ (for instance, that found in Elzenaar~\cite[Definition~3.2]{elzenaar25h}) and then find the triples corresponding to the holonomy representation of
any tunnel number $1$ link by solving algebraic equations in the entries of the words which we will describe. The practical disadvantages of this are numerous: (i) the word recursion
is difficult to compute since the complex of curves of the genus $2$ surface is much more complicated than that of the four punctured sphere; (ii) even if one can actually compute the relevant polynomials, identifying the correct (`geometric')
root is difficult; and (iii) for all tunnel number $1$ links that are combinatorially simple enough that it is plausible to use computational techniques, one can already use the
computer programme \texttt{SnapPy}~\cite{SnapPy} to find holonomy representations very quickly.

\subsection{Leitfaden}
We begin with \zcref{sec:curve_coding}, where we recall the idea of `curve coding' which connects symbolic dynamics to geometric group theory and topology. We deal
with this fairly precisely as some technical details in the latter parts of the paper depend on these ideas. The casual reader may wish to skip this section.

The first section of the paper where we actually deal with tunnel number $1$ objects is \zcref{sec:heegaard_general}. Here, we define carefully what a tunnel number $1$ link or
manifold is and relate it to Heegaard splittings of genus $2$ (\zcref{sec:tno_manifolds}), we recall the notion of the curve complex of a surface together with some useful results
(\zcref{sec:curve_complex}), and we discuss paths in the curve complex that connect a given tunnel number $1$ link to a fixed basepoint and give an outline of the main construction (\zcref{sec:heegaard}).
This outline, which introduces the notion of a \df{ns-geodesic}, is very important for the development of the theory.

In \zcref{sec:goeritz}, we relate the curve complex to the Teichm\"uller space of the genus $2$ surface and we discuss some connections to representation space. We use this
to say a few things about the fibres of the map that sends Heegaard diagrams to the associated tunnel number $1$ link. The remainder of the paper does not depend on this material,
and the reader primarily interested in knot theory may want to skip this section on a first reading.

In \zcref{sec:fan_edge}, we consider paths in the curve complex consisting of $2$-simplices so that adjacent simplices in the path intersect in an edge. We study this from the
point of view of the Farey graph, and we give a recursion on words representing vertices of the simplices in such a path. This is used to prove our first filtration result
(\zcref{thm:filtration_in_component}). This in turn sets up \zcref{sec:tno_in_sphere}, where we study tunnel number $1$ links in the $3$-sphere via our machinery. This section
can be viewed as a reinterpretation of some results of Cho and McCullough~\cite{cho09}.

In \zcref{sec:fan_vertex}, we allow more general paths in the curve complex where adjacent simplices are now allowed to intersect in a single vertex. We again give a recursion
on words representing the vertices of the path and a resulting filtration \zcref{thm:filtration_global}; this is based on the notion of the `tripants graph' studied by Betts,
Larsen, Utley, and Vanis~\cite{betts21}, which allows us to give a recursion based on the standard theory of dynamics on the once-punctured torus. We use this second filtration
theorem in \zcref{sec:tno_general} to study tunnel number $1$ links in arbitrary ambient closed $3$-manifolds, and in particular we show that every such link is built from the
unknot via a sequence of `rational braiding moves' (essentially braid moves) and suitable Dehn fillings that arise from the walk in the tripants graph along sets of curves on a twice-punctured
torus. In this final section we also say a few words about general tunnel number $1$ manifolds (\zcref{sec:tno_manifold_general}).

In \zcref{app:pleated_surfaces} we recall a folklore result in hyperbolic geometry about geometrically infinite surface subgroups of
geometrically finite $3$-manifold groups.

\subsection{Notational conventions}
We use blackletter type (e.g.\ $ \mf{k}, \mf{l}, \mf{b} $) to denote links in a closed $3$-manifold, i.e.\ smoothly embedded submanifolds isomorphic
to disjoint unions of copies of $ \IS^1$. Generally we will use the letter $A$ to denote the ambient $3$-manifold of a link. The topological surface
of genus $g$ with $n$ punctures is denoted $ S_{g,n} $, and its curve complex (\zcref{defn:curve_complex}) by $ \KK(S_{g,n}) $. The subset of $m$-simplices,
for $ m \in \IN $, is denoted by $ \KK^{(m)}(S_{g,n}) $. Generally we will reserve the letter $ \Delta $ for top-level simplices of the curve complex,
and we will use Greek letters (e.g.\ $ \alpha, \beta, \gamma $) for vertices of the curve complex (freely identified with isotopy classes of simple
closed curves on the relevant surface). We will define a distinguished top-level simplex $ \Delta_0 $ of $ \KK(S_{2,0}) $ (\zcref{defn:deltanought}).
In several places we will construct from the curve complex some other complexes which we will call `dual complexes';
the reader should take care with the definitions as these dual complexes are usually \emph{not} the standard abstract dual complex. The mapping class group $ \Mod(S_{g,n}) $
acts on $ \KK(S_{g,n}) $ by isometry.

We use $ \mc{F} $ for the Farey tree/graph (which is a uniform trivalent tree), and $ \mc{F}^* $ for the Farey complex (the abstract simplicial complex associated to the Farey
triangulation of $ \IH^2 $). Some further notation which is used less frequently includes $ \Goer(A,\Sigma) $ for the Goeritz group of a Heegaard splitting
of a manifold $A$ by a compact surface $ \Sigma $ (\zcref{defn:goeritz}), $ \Tw $ for the twist group, $\mc{D}(A,\Sigma)$ for the nonseparating disc complex of a distinguished-by-context handlebody
of the Heegaard splitting (\zcref{defn:disc_complex}), $\tilde{\mc{T}}(A,\Sigma)$ for a `dual tree' to $\mc{D}(A,\Sigma)$ (the lift of the Cho--McCullough tree, \zcref{defn:lifted_cm_tree}),
and $ \TP(T) $ for the tripants graph of a twice-holed torus $T$ (\zcref{defn:tripants_graph}).

If $M$ is a hyperbolic manifold, then we denote the holonomy group of $ M $ by $ \Hol(M) $; so $ \Hol(M) $ is a discrete subgroup of $ \PSL(2,\IC) $ isomorphic to $\pi_1(M) $
and $ M $ is isometric to $ \IH^3/\Hol(M) $. Further, if $ S \subset M $ is a submanifold (e.g.\ an embedded surface) then we write $ \Hol(S) $ for the image of the canonical
map $ \pi_1(S) \hookrightarrow \pi_1(M) \stackrel{\sim}{\smash{\longrightarrow}\rule{0pt}{0.4ex}} \Hol(M) $.

\subsection{Acknowledgments}
I thank Jessica Purcell for suggesting the study of tunnel number $1$ links from this point of view and for discussing the mathematics. I thank Connie On Yu Hui, Gaven Martin,
and Jeroen Schillewaert for additional discussions about this work at various points. I thank Annie Holden for correspondence about handlebody groups.

\section{Cutting sequences and curve coding}\label{sec:curve_coding}
In this section we will briefly recall the concept of curve coding for surface groups; the idea originated with Koebe and Morse, see Series~\cite{series86}, and can be viewed as
a generalisation of the notion of a Sturmian sequence~\cite[Chapter 6]{fogg02} on a more complicated (non-Euclidean) grid. Let $S$ be a closed topological surface of genus $g$.
Then $S$ is homeomorphic to the connected sum $ T_1 \# \cdots \# T_g $ where each $ T_i $ is a torus. It is well-known that a generating set for $ \pi_1(S) $ is given by
\begin{displaymath}
  \mc{A} = \{ A_1, B_1, \ldots, A_g, B_g \}
\end{displaymath}
where each $A_i $ is basepoint-free homotopic on $S$ to an oriented meridian $ \overline{A_i}$ of $ T_i $ and each $ B_i $ is basepoint-free homotopic to an oriented longitude $\overline{B_i}$ of $T_i $. The words `meridian'
and `longitude' really have no intrinsic meaning, and all we require is that $ \overline{A_i}$ intersects $\overline{B_i}$ exactly once (with positive crossing orientation),
and that $\{\overline{A_i},\overline{B_i}\}$ is a basis for the $\IZ$-module $ H_1(T_i) $.

\begin{lem}
  Let $ \gamma $ be a closed loop on $S$ which intersects each homology basis $\{\overline{A_i},\overline{B_i}\}$ transversely. Fix a basepoint $ \bullet $ on $ \gamma $, and an
  orientation of $ \gamma $. This defines a word $[\bar{\gamma}]$ in the alphabet
  \begin{displaymath}
    \overline{\mc{A}} = \{ \overline{A_1}^{\pm 1}, \overline{B_1}^{\pm 1}, \ldots, \overline{A_g}^{\pm 1}, \overline{B_g}^{\pm 1} \};
  \end{displaymath}
  namely, walk along $\gamma$ from the basepoint $ \bullet $ in the direction given by the orientation, and write down the symbol of each meridian or longitude of some $ T_i $ as
  it is crossed with positive crossing orientation, or the inverted symbol when a meridian or longitude is crossed with negative orientation. Define the word $ [\gamma] $ to be
  the image of $ [\overline{\gamma}] $ under the map on words induced by the map of alphabets from $\overline{\mc{A}}$ to $\mc{A}$ given by
  \begin{displaymath}
    \overline{A_i}^{\pm 1} \mapsto B_i^{\pm 1} \quad \text{and}\quad \overline{B_i}^{\pm 1} \mapsto A_i^{\pm 1}.
  \end{displaymath}
  Then $ [\gamma] $, interpreted as an element of $ \pi_1(S)$, is the group element representing the oriented curve $ \gamma $. As $\gamma$ varies over all representatives
  of its isotopy class, the word $ [\gamma] $ varies by insertion or deletion of subwords of the form $ XX^{-1} $; as $\bullet$ moves on $\gamma$, the word $[\gamma]$ varies
  by cyclic permutations. If the orientation of $\gamma$ is reversed, then $[\gamma] $ changes by reversal of its letters. \qed
\end{lem}

Suppose now that $ G $ is a genus $g$ compression body group (for example, a Kleinian group which lies in the closure of the Schottky space $ \mc{S}_g$). Then there is a surjective
map $ \pi_1(S) \to G $, and so we can do curve coding in $G$.

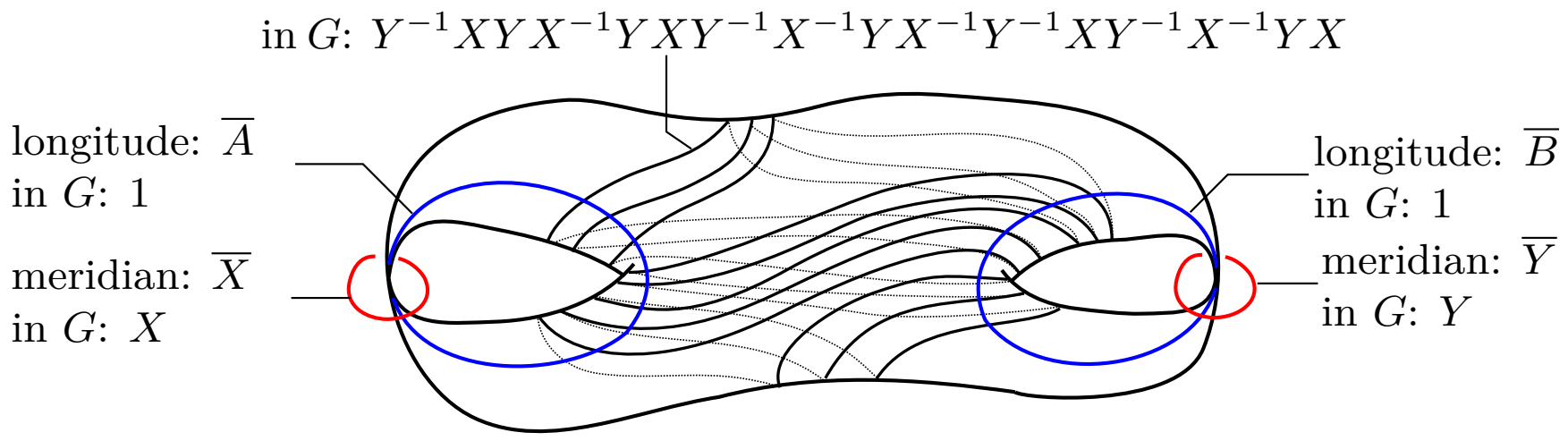
\begin{figure}\vspace{1em}
  \labellist
  \small\hair 2pt
  \pinlabel {\parbox{2cm}{meridian: $ \overline{X} $\\in $G$: $X$}} [r] at 0 53
  \pinlabel {\parbox{2cm}{longitude: $ \overline{A} $\\in $G$: $1$}} [r] at 0 105
  \pinlabel {\parbox{2cm}{meridian: $ \overline{Y} $\\in $G$: $Y$}} [l] at 394 59
  \pinlabel {\parbox{2cm}{longitude: $ \overline{B} $\\in $G$: $1$}} [l] at 391 101
  \pinlabel {\parbox{6cm}{in $G$: $Y^{-1}XYX^{-1}YXY^{-1}X^{-1}YX^{-1}Y^{-1}XY^{-1}X^{-1}YX$}} [b] at 145 145
  \endlabellist
  \centering
  \includegraphics[width=.6\textwidth]{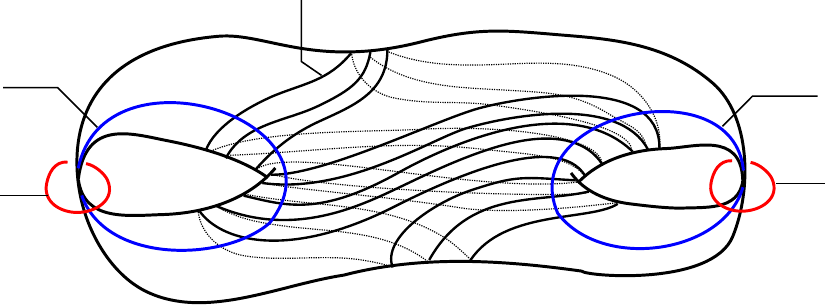}
  \caption{Curve coding in a group on the boundary of genus two Schottky space.\label{fig:curve_cutting_whitehead}}
\end{figure}

\begin{ex}
  Let $G = \langle X,Y \rangle $ be Kleinian where $X $ and $ Y $ are parabolics, and suppose $ \Omega(G)/G = S_{0,4} $. Then $ \IH^3/G$ is homeomorphic to a genus $2$ handlebody
  with hyperbolic structure that has two rank $1$ cusps represented by the elements $X$ and $Y$ in $G$~\cite{keen94}. We may choose longitudes for these cusps
  which are trivial in the handlebody $ \IH^3/G $ but not in the surface; say they are represented by the letters $ A,B \in \pi_1(S_{2,0}) $. In particular we can take our homology
  bases for the two handles to be $ \{ \overline{X},\overline{A} \} $ and $ \{ \overline{Y}, \overline{B}\} $. Now an application of curve coding tells us that:
  \begin{itemize}
    \item The word in $ G $ corresponding to the meridian loop $ \overline{X} $ is $X$, since the loop $\overline{X}$ meets $ \overline{A} $ exactly once and then we replace $ \overline{A} \mapsto X $
    when we change alphabets from $ \{\overline{X},\overline{A},\overline{Y},\overline{B}\} $ to $\{A,X,B,Y\}$; then the group element $X \in \pi_1(S_{2,0})$ is preserved by the map $ \pi_1(S_{2,0}) \to G $.
    \item The word in $ G $ corresponding to the longitude loop $ \overline{A} $ is $1$, since the loop $\overline{A}$ meets $ \overline{X} $ exactly once and then we replace $ \overline{X} \mapsto A $
    when we change alphabets from $ \{\overline{X},\overline{A},\overline{Y},\overline{B}\} $ to $\{A,X,B,Y\}$; then the group element $A \in \pi_1(S_{2,0})$ is killed by the map $ \pi_1(S_{2,0}) \to G $.
    \item The word in $G$ corresponding to the $3/8$ loop is
          \begin{displaymath}
            Y^{-1}XYX^{-1}YXY^{-1}X^{-1}YX^{-1}Y^{-1}XY^{-1}X^{-1}YX
          \end{displaymath}
          since when walking along the curve we write down the cutting sequence
          \begin{displaymath}
            \overline{B}^{-1}\overline{A}\overline{B}\overline{A}^{-1}\overline{B}\overline{A}\overline{B}^{-1}\overline{A}^{-1}\overline{B}\overline{A}^{-1}\overline{B}^{-1}\overline{A}\overline{B}^{-1}\overline{A}^{-1}\overline{B}\overline{A}
          \end{displaymath}
          and then we change alphabets from $ \{\overline{X},\overline{A},\overline{Y},\overline{B}\} $ to $\{A,X,B,Y\}$.
  \end{itemize}
  See \zcref{fig:curve_cutting_whitehead}.
\end{ex}

In reality, to avoid the annoyance of changing alphabets, we normally label the surface with the letters $ A_1, B_1, \ldots, A_n, B_n $ so that every time the geodesic $\gamma$
crosses $A_i$ it picks up a letter $ A_{i}^{\pm1} $ (depending on orientation) and every time it crosses $B_i$ it picks up a letter $ B_{i}^{\pm1} $. Formally, we are labelling the
curve $ \overline{A_i} $ with $ B_i $ and we are labelling $ \overline{B_i} $ with $ A_i $. This has the advantage of being much clearer to do calculations, but the disadvantage of
the confusion between two different loops with the same label: $ A_i $ represents both an element of the group and the projection of its axis onto the surface, \emph{and} the loop
on the surface such that if a curve crosses it then the word of the curve picks up an $ A_i^{\pm 1} $. We will give one example with this convention now.

\begin{figure}\vspace{1.5ex}
  \labellist
  \small\hair 2pt
  \pinlabel {$MP^{-1}MQ^{-1}$} [r] at 124 115
  \pinlabel {$M$} [l] at 324 74
  \pinlabel {$1$} [l] at 317 121
  \pinlabel {$P$} [b] at 19 103
  \pinlabel {$Q$} [r] at 24 39
  \endlabellist
  \centering
  \includegraphics[width=.6\textwidth]{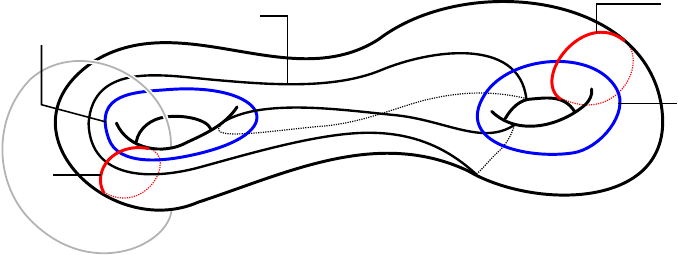}
  \caption{Curve coding in a $(1;2)$-compression body group. The coloured loops labelled with single letters are the loops defining cutting sequences. The black curve has word $MP^{-1}MQ^{-1}$
  since it crosses the curves labelled $M$, $P$, $M$, then $Q$ in that order, and each time it crosses $M$ it does so with the same orientation. Note that the loop labelled $M$ is actually
  represented by the trivial word in the fundamental group due to our convention!\label{fig:two_cpt_lnk}}
\end{figure}
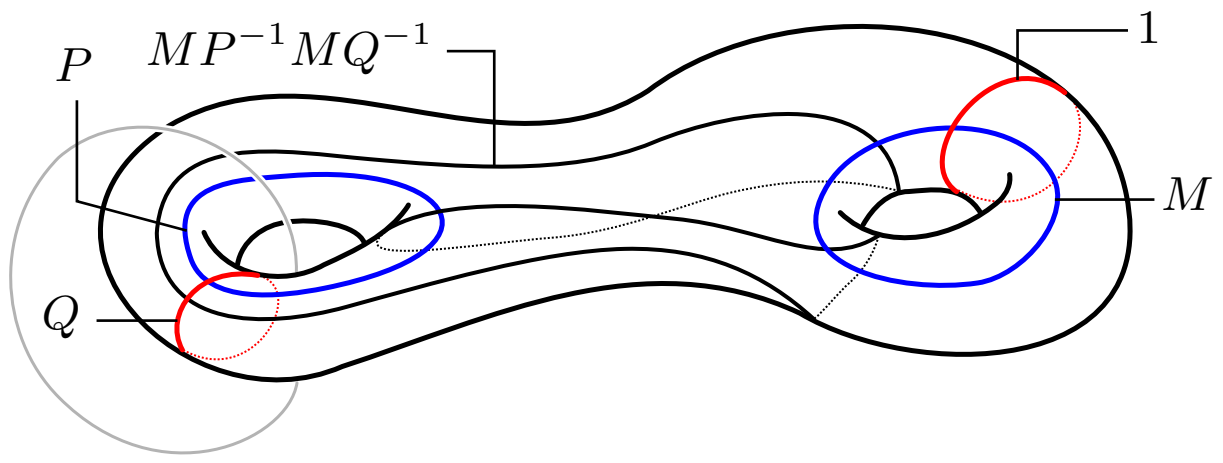

\begin{ex}
  A $(1;2)$-compression body is obtained by gluing a $1$-handle onto a thickened torus, and has fundamental group $ G = \langle P, Q, M : [P,Q]=1\rangle$. The genus $2$ surface
  surjects onto $G$ by killing one generator, and the subgroup $ \langle P,Q\rangle $ represents a genus $1$ boundary component of the compression body. In \zcref{fig:two_cpt_lnk}
  we give an example of curve coding.
\end{ex}

\section{Heegaard splittings and tunnel number} \label{sec:heegaard_general}

\subsection{Tunnel number one manifolds}\label{sec:tno_manifolds}
If $ M $ is a manifold then we write $ N(A) $ for a regular neighbourhood of a subset $ A \subset M $.
\begin{defn}
  A link $ \mf{k} \subset \IS^3 $ is said to have \df{tunnel number $1$} if there exists some embedded 1-cell $ \tau \subset \IS^3 \setminus \mf{k} $
  such that $ \IS^3 \setminus N(\mf{k} \union \tau) $ is a genus $2$ handlebody. The cell $ \tau $
  is called an \df{unknotting tunnel} for $\mf{k}$.
\end{defn}

\begin{figure}
  \labellist
  \small\hair 2pt
  \pinlabel {$\tau$} [r] at 0 63
  \pinlabel {$\tau'$} [l] at 328 57
  \endlabellist
  \centering
  \includegraphics[width=.5\textwidth]{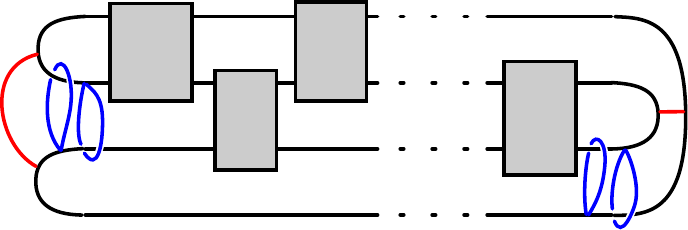}
  \caption{A two-bridge link with its unknotting tunnels; the shaded boxes represent twist regions. The upper and lower tunnels $\tau$ and $\tau'$ are unknotting tunnels for links of
  both $1$ and $2$ components; the other four are tunnels only when the link is a knot.\label{fig:2bridge_braided}}
\end{figure}

\begin{ex}\label{ex:2bridge}
  A \df{two-bridge link} is a link which admits a diagram in $ \IR^2 $ with exactly two local maxima in the $y$-ordinate~\cite[Chapter~12]{burde3e},
  Every two-bridge link is a tunnel number $1$ link. Up to isotopy, every two-bridge knot has at most six unknotting tunnels~\cite{kobayashi99,goda00} (and
  the isotopy classes within these six are classified by by Morimoto--Sakuma~\cite[Theorem~5.2]{morimoto91}),
  and every nontrivial two-bridge link of two components has exactly two unknotting tunnels~\cite{adams96}. These tunnels are shown in \zcref{fig:2bridge_braided}.
\end{ex}

\begin{ex}[Non-hyperbolic examples]
  A \df{torus knot} is a knot which embeds on the boundary of the genus $1$ Heegaard splitting of  $ \IS^3 $.
  Every torus knot is tunnel number $1$, and in general has three tunnels up to isotopy~\cite{boileau88}.
  In addition, the tunnel number $1$ satellite links are classified by Morimoto and Sakuma~\cite{morimoto91}.
\end{ex}

\begin{figure}
  \labellist
  \small\hair 2pt
  \pinlabel {\parbox{6ex}{\centering $a$\\half\\twists}} at 36 166
  \pinlabel {\parbox{6ex}{\centering $b$\\half\\twists}} at 178 166
  \pinlabel {\parbox{6ex}{\centering $c$\\half\\twists}} at 326 166
  \pinlabel {$\tau$} [br] at 282 270
  \endlabellist
  \centering
  \includegraphics[width=.4\textwidth]{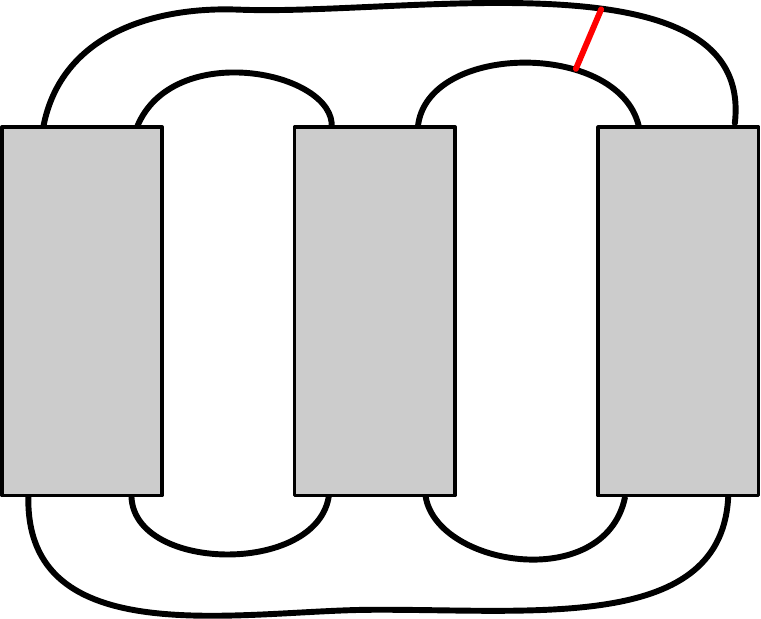}
  \caption{The pretzel link $ \mf{P}(a,b,c) $ with unknotting tunnel $\tau$.\label{fig:pretzel}}
\end{figure}
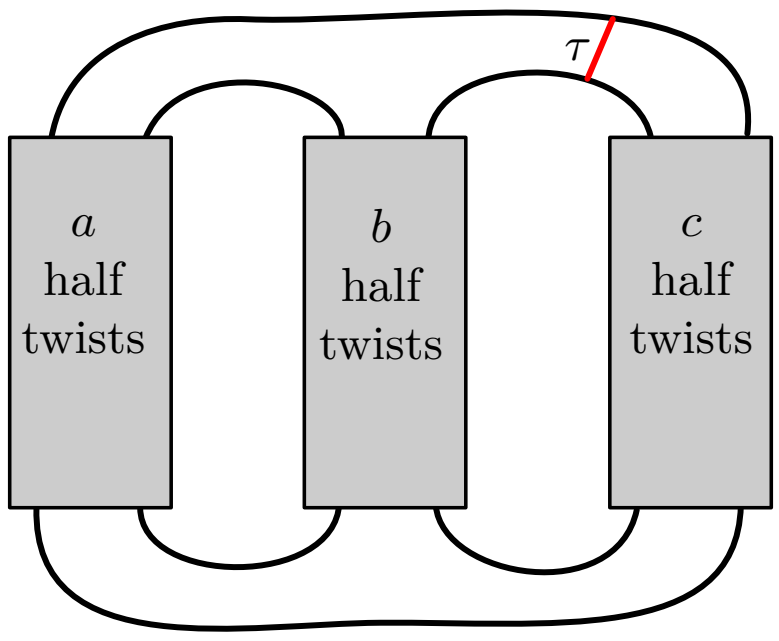

\begin{ex}
  If $ a,b,c \in \IZ $ then the \df{pretzel link} $ \mf{P}(a,b,c) $ is defined as shown in \zcref{fig:pretzel}. A sufficient condition for $ \mf{P}(a,b,c) $ to have tunnel number $1$ is
  that $ b = 2 $ and $ a \equiv c \equiv 1 \pmod{2} $~\cite[Theorem~2.2, condition (2)]{morimoto96}. We showed an unknotting tunnel for $ \mf{P}(3,2,3) $ in \zcref{fig:eight_five_unknotting}.
\end{ex}

\begin{rem}
  Further information on the classical (i.e.\ knot-theoretic) theory of unknotting tunnels, and further examples, may be found in the survey article of Sakuma~\cite{sakuma98}.
  A full list of tunnel number $1$ prime knots with at most ten crossings was compiled by Morimoto, Sakuma, and Yokota~\cite[Theorem 2.6]{morimoto96}, and of alternating knots
  of $11$ or $12$ crossings by Castellano-Mac{\'{\i}}as and Owad~\cite{castellano22}.
\end{rem}

If $ \mf{k} $ is a tunnel number $1$ link with unknotting tunnel $ \tau $, then both $ N(\mf{k} \union \tau) $ and $ \IS^3 \setminus N(\mf{k} \union \tau) $
are genus $2$ handlebodies. In particular, the pair $ (\mf{k},\tau) $ determines a Heegaard splitting of $ \IS^3 $. This point of view is classical (e.g.\ it is
used by the work of Morimoto and Sakuma~\cite{morimoto91}) and was used by Cho and McCullough~\cite{cho09} to impose a natural combinatorial tree structure
on the set of all isotopy classes of pairs $ (\mf{k},\tau) $, where $ \tau $ is an unknotting tunnel for $ \mf{k} $.

\begin{defn}
  If $A$ is a closed $3$-manifold, then a \df{genus $2$ Heegaard diagram} for $A$ is a set $ D= \{\sigma_1, \rho, \sigma_2\}$
  of three disjoint isotopy classes  of curves on the boundary of a genus $2$ handlebody $ \mc{H} $ such that $A$ is homeomorphic to the manifold
  obtained by gluing $2$-handles to $ \mc{H} $ along the three elements of $D$ and then gluing $3$-balls onto the resulting spherical boundary components.
  If such a diagram is given then the handlebody $ \mc{H} $ embeds into $ A $, and the three $2$-handles bounded by the elements of $D$ embed into $ A \setminus \mc{H} $.
  We will refer to $ \mc{H} $ as the \df{exterior handlebody} of the Heegaard splitting of $A$ determined by $D$.
\end{defn}

\begin{figure}
  \centering
  \begin{subfigure}{.48\textwidth}
    \centering\includegraphics[height=4cm]{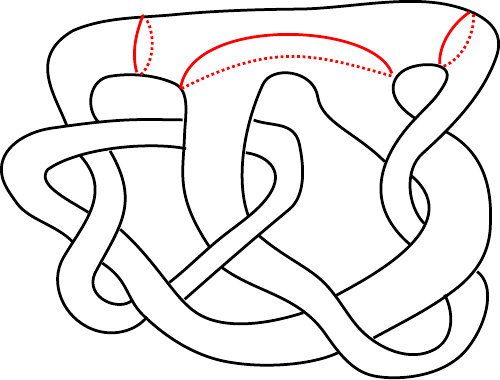}
    \caption{}
  \end{subfigure}
  \begin{subfigure}{.48\textwidth}
    \centering\includegraphics[height=4cm]{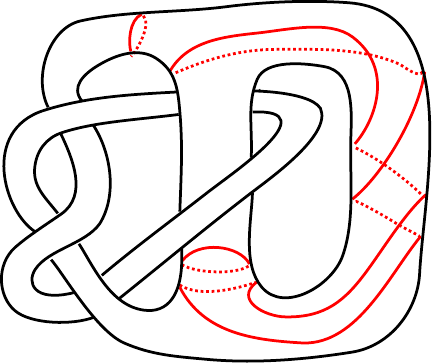}
    \caption{}
  \end{subfigure}\\
  \begin{subfigure}{.48\textwidth}
    \labellist
    \tiny\hair 2pt\color{red}
    \pinlabel {$1$} [r] at 73 158
    \pinlabel {$1$} [tr] at 92 140
    \pinlabel {$7$} [t] at 162 164
    \pinlabel {$6$} [r] at 217 127
    \pinlabel {$2$} [r] at 204 91
    \pinlabel {$2$} [tl] at 199 45
    \pinlabel {$8$} [l] at 112 125
    \pinlabel {$1$} [r] at 95 99
    \pinlabel {$1$} [r] at 118 44
    \pinlabel {$10$} [r] at 115 80
    \pinlabel {$6$} [r] at 110 12
    \pinlabel {$12$} [tl] at 133 21
    \pinlabel {$1$} [t] at 151 30
    \endlabellist
    \centering\includegraphics[height=4cm]{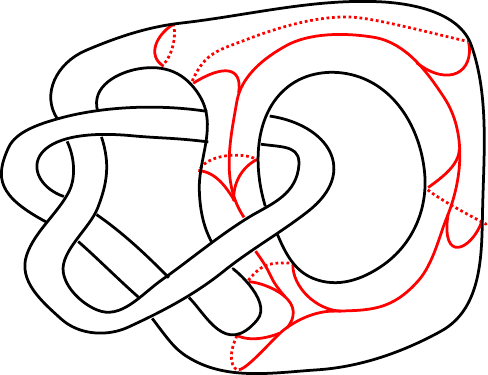}
    \caption{}
  \end{subfigure}%
  \begin{subfigure}{.48\textwidth}
    \labellist
    \tiny\hair 2pt\color{red}
    \pinlabel {$1$} [r] at 20 127
    \pinlabel {$1$} [r] at 59 135
    \pinlabel {$18$} [b] at 138 138
    \pinlabel {$2$} [r] at 180 81
    \pinlabel {$14$} [l] at 190 73
    \pinlabel {$2$} [t] at 178 52
    \pinlabel {$12$} [l] at 188 47
    \pinlabel {$1$} [b] at 156 42
    \pinlabel {$1$} [l] at 167 37
    \pinlabel {$5$} [r] at 158 33
    \pinlabel {$6$} [t] at 154 18
    \pinlabel {$12$} [r] at 105 18
    \pinlabel {$18$} [r] at 140 34
    \pinlabel {$24$} [t] at 127 53
    \pinlabel {$2$} [r] at 76 90
    \pinlabel {$20$} [l] at 85 106
    \endlabellist
    \centering\includegraphics[height=4cm]{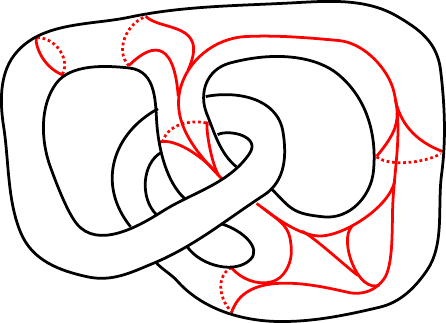}
    \caption{}
  \end{subfigure}\\
  \begin{subfigure}{\textwidth}
    \labellist
    \tiny\hair 2pt\color{red}
    \pinlabel {$1$} [r] at 21 128
    \pinlabel {$1$} [l] at 59 147
    \pinlabel {$18$} [b] at 118 145
    \pinlabel {$2$} [r] at 149 85
    \pinlabel {$14$} [l] at 166 90
    \pinlabel {$2$} [b] at 157 54
    \pinlabel {$12$} [l] at 173 37
    \pinlabel {$5$} [r] at 127 30
    \pinlabel {$1$} [b] at 133 40
    \pinlabel {$1$} [l] at 140 32
    \pinlabel {$18$} [b] at 104 60
    \pinlabel {$6$} [r] at 109 23
    \pinlabel {$12$} [r] at 52 31
    \pinlabel {$24$} [r] at 90 76
    \pinlabel {$24$} [bl] at 128 61
    \pinlabel {$20$} [l] at 86 115
    \pinlabel {$2$} [r] at 61 103
    \endlabellist
    \centering\includegraphics[height=4cm]{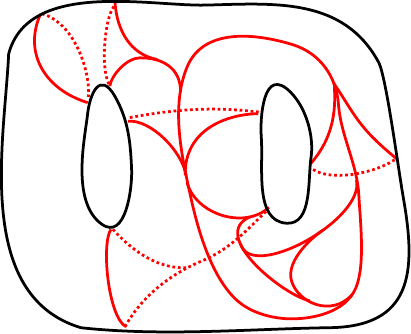}
    \caption{}
  \end{subfigure}%
  \caption{The Heegaard diagram for the $8_5$ knot.\label{fig:8_5_twisting}}
\end{figure}
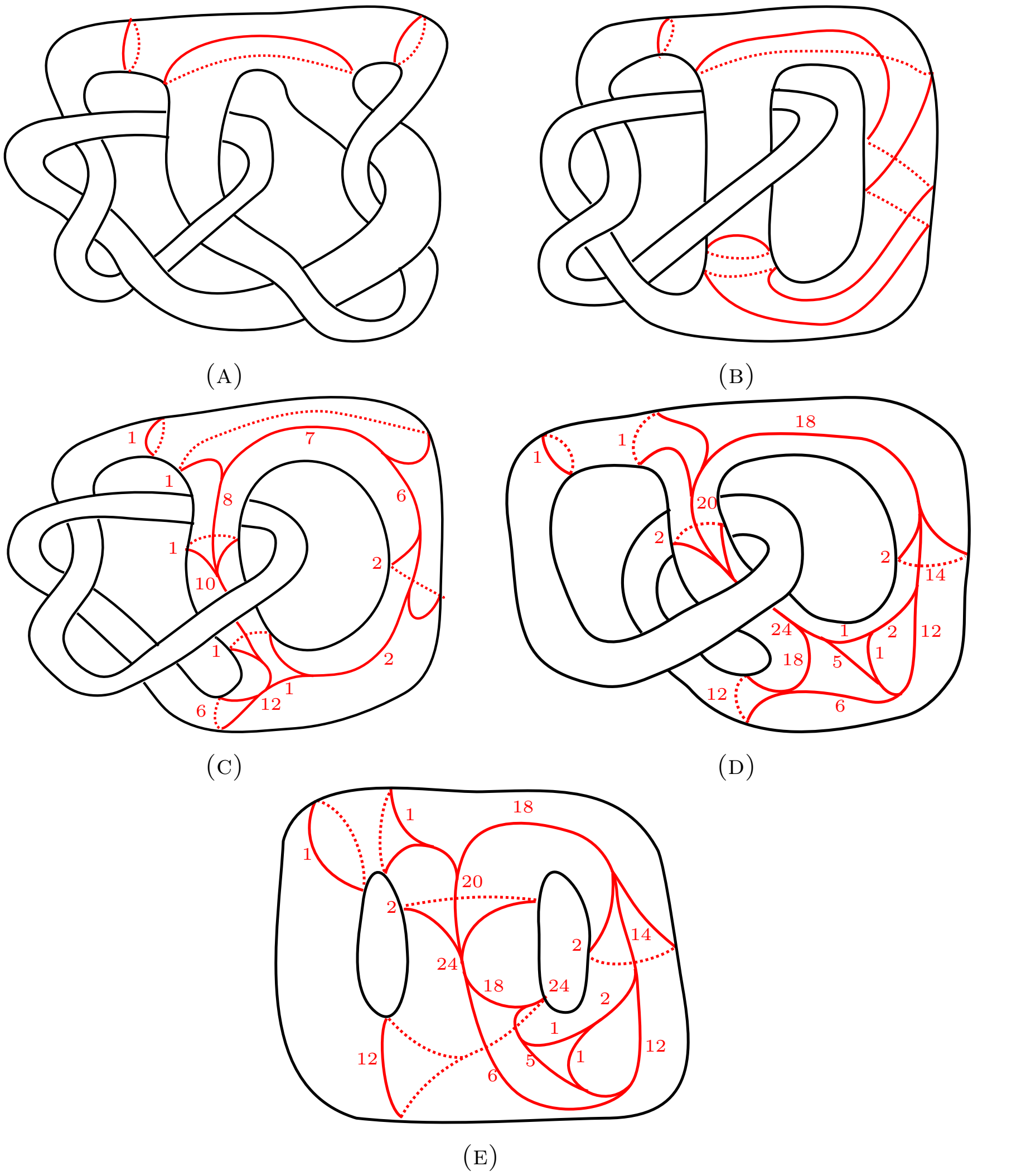

\begin{ex}
  In \zcref{fig:8_5_twisting}, we show the computation of the Heegaard diagram associated to the unknotting tunnel of the \href{https://katlas.org/wiki/8_5}{$ 8_5 $ knot} which we saw
  in \zcref{fig:eight_five_unknotting}. Because the curves are so complicated, we have used Thurston's train track notation~\cite[\S 8.9]{thurstonN} for the last three
  images. By work of Morimoto, Sakuma, and Yokota~\cite[Theorem~2.6]{morimoto96} this knot is one of the two non-two-bridge tunnel number $1$ knots of minimal
  crossing number (the other is the \href{https://katlas.org/wiki/8_10}{$ 8_{10} $ knot}); thus, this is the simplest example where
  the Cho--McCullough machinery is actually needed to study the tunnel. The complexity of even this simplest
  example explains why we do not give explicit computations of the various filtrations in this paper. This also justifies our claim in \cite[Example~2.4]{elzenaar25h}
  that for the $8_5$ knot the words involved are `very long'.
\end{ex}

We can now generalise the notion of tunnel number $1$ to links in arbitrary closed $3$-manifolds.
\begin{defn}\label{lem:tunnel_number_1_is_heegaard}
  Let $ A $ be a closed $3$-manifold. A link $ \mf{k} \subset A $ is called \df{tunnel number $1$} if either of the following equivalent conditions holds:
  \begin{itemize}
    \item There exists a $1$-cell $ \tau \subset A $ such that $A \setminus N(\mf{k} \union \tau) $ is a genus $2$ handlebody.
    \item There exists a genus $2$ Heegaard diagram for $A$ consisting of three curves $ D= \{\sigma_1, \rho, \sigma_2\} $
          such that $ A \setminus N(\mf{k}) $ is isotopic (in $A$) to $ \mc{H} \union \Delta_\rho $ where $ \Delta_\rho $
          is the $2$-handle with boundary $ \rho $ and $ \mc{H} $ is the exterior handlebody of the Heegaard splitting.
  \end{itemize}
\end{defn}

If $ A $ is a closed $3$-manifold which contains a tunnel number $1$ link, then it has Heegaard genus $2$. By a theorem of Viro~\cite{viro72}, Birman
and Hilden~\cite[Theorem~5]{birman75}, and Montesinos~\cite[Corollary 2]{montesinos75}, this implies that $A$ is a two-fold cover of $ \IS^3 $ branched over a $3$-bridge link.

\begin{figure}\vspace{1em}
  \labellist
  \small\hair 2pt
  \pinlabel {$2$-handle boundary} [b] at 229 154
  \endlabellist
  \centering
  \includegraphics[width=.6\textwidth]{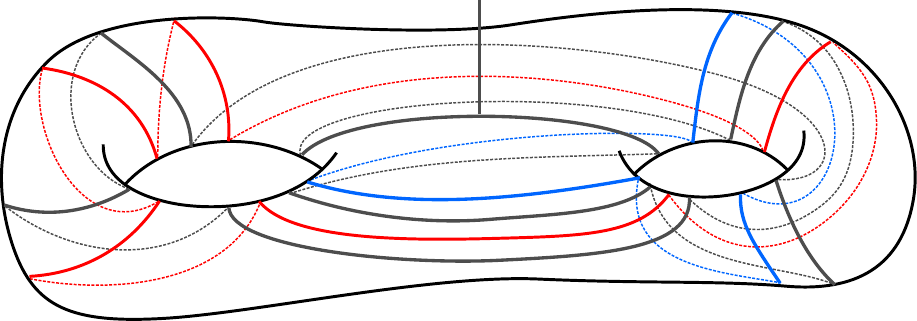}
  \caption{The manifold \texttt{m003(0,0)} is a tunnel number $1$ link in a lens space, and arises from gluing a $2$-handle to the unbounded handlebody along the indicated curve.\label{fig:m003}}
\end{figure}

\begin{ex}\label{ex:m003}
  Let $M$ be the manifold \texttt{m003(0,0)} in the \textsc{SnapPy} census of orientable cusped $3$-manifolds~\cite{SnapPy}. Then $ M $
  is a link complement in the lens space $ L(3,10) $~\cite[Example~3.4]{kegel24} and is tunnel number $1$ arising from the Heegaard diagram shown in \zcref{fig:m003}~\cite[Theorem~3.8]{elzenaar26pr}.
\end{ex}

\begin{ex}\label{ex:punc_torus_bund}
  If $ S $ is a once-punctured torus and $ f \in \SL(2,\IZ) $ is a mapping class of $S $ then the mapping torus
  \begin{displaymath}
    M_f = S \times [0,1] / ((x, 0) \sim (f(x), 1))
  \end{displaymath}
  is tunnel number $1$ in a lens space if $f$ is either elliptic ($\tr^2 f < 4 $) or if $f$ or its inverse is conjugate to
  \begin{displaymath}
    \pm\begin{bmatrix}
      1&1\\0&1
    \end{bmatrix}
    \quad\text{or}\quad
    \pm\begin{bmatrix}
      1&1\\0&1
    \end{bmatrix}^n
    \begin{bmatrix}
      1&0\\1&1
    \end{bmatrix}
  \end{displaymath}
  for some $n\in\IN$. This result is attributed to Klaus Johannson by Sakuma~\cite[\S 1.7]{sakuma98}. \zcref[S]{ex:m003} is a special case, where
  \begin{displaymath}
    f = -\begin{bmatrix}
      1&1\\0&1
    \end{bmatrix}
    \begin{bmatrix}
      1&0\\1&1
    \end{bmatrix}.
  \end{displaymath}
\end{ex}

In general:
\begin{defn}\label{defn:tno_manifold}
  A \df{tunnel number $1$ manifold} is one of the following equivalent things:
  \begin{itemize}
    \item A $3$-manifold obtained from a genus $2$ handlebody $ \mc{H}^+ $ by gluing $2$-handles along $1$, $2$, or $3$ disjoint isotopy classes of curves on the boundary and then gluing $3$-balls onto
          any $2$-sphere boundary components.
    \item The result of gluing a genus $2$ compression body onto a genus $2$ handlebody along the genus $2$ surface boundary components of each.
    \item A Dehn filling of a tunnel number $1$ link in an arbitrary closed $3$-manifold.
  \end{itemize}
\end{defn}
These manifolds are the same as those which Johannson called \df{$N$\!-relator $3$-manifolds} of genus $2$~\cite{johannson95} (for $N \in \{1,2,3\} $).

\subsection{Curve complexes}\label{sec:curve_complex}
We recall the following definition which is due originally to Harvey~\cite{harvey81} and which has become central in modern Teichm\"uller theory.
\begin{defn}\label{defn:curve_complex}
  Let $ \Sigma $ be a Riemann surface (connected but possibly with punctures), and let $ \mu $ be the minimal possible intersection number of two distinct isotopy classes of essential closed curves
  on $ \Sigma $. If $ \chi(\Sigma) < 0 $ and $ \Sigma \not\in \{ S_{0,3}, S_{0,4}, S_{1,1} \} $ then $ \mu = 0 $. If $ \Sigma = S_{0,4} $ then $ \mu = 2 $, if $ \Sigma = S_{1,1} $ or $ \Sigma = S_{1,0} $
  then $ \mu = 1 $.

  The \df{curve complex} $ \KK(\Sigma) $ is the simplicial complex supported on the set of isotopy classes of essential closed curves on $ \Sigma $ such that a set of curves forms
  a simplex if the elements of the set have pairwise intersection number $ \mu $. If $ \Sigma = S_{0,n} $ for $ n \in {0,1,2,3} $, then $ \KK(\Sigma) $ is the complex $ \{\emptyset\} $
  (of dimension $-1$).
\end{defn}

We will primarily be interested in $ \KK(S_{2,0}) $ and in related curve complexes which can be embedded as subcomplexes in natural ways. The most important related complex is:
\begin{defn}\label{defn:farey_graph}
  The \df{Farey triangulation} $\mc{F}^*$ is the triangulation of $ \IH^2 $ with vertex set $ \IQ \union \{\infty\} $ where $ p/q $ and $ r/s $ are joined by an edge if and only if $ \abs{ps - qr} = 1 $.
  The \df{Farey graph} $\mc{F} $ is the dual graph to the $2$-simplices of $ \mc{F}^* $, and is an infinite trivalent tree.
\end{defn}

\begin{figure}
  \begin{subfigure}{.5\textwidth}\centering
    \includegraphics[height=20mm]{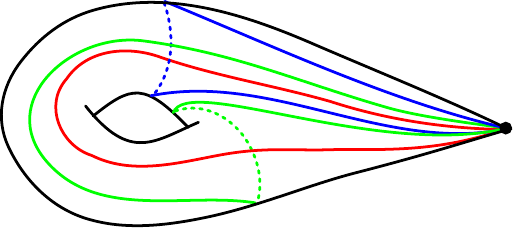}
    \caption{Ideal triangulation of $S_{1,1}$.}
  \end{subfigure}%
  \begin{subfigure}{.5\textwidth}\centering
    \includegraphics[height=20mm]{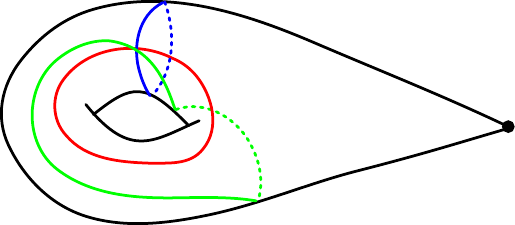}
    \caption{$2$-simplex of $\KK(S_{1,1})$}
  \end{subfigure}\\\vspace{1ex}
 \begin{subfigure}{.5\textwidth}\centering
    \includegraphics[height=30mm]{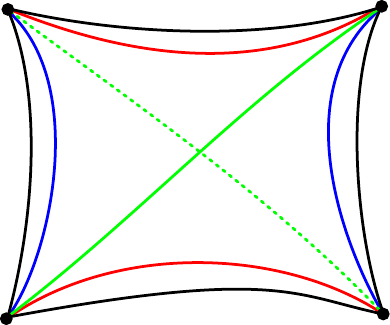}
    \caption{Ideal triangulation of $S_{0,4}$.}
  \end{subfigure}%
  \begin{subfigure}{.5\textwidth}\centering
    \includegraphics[height=30mm]{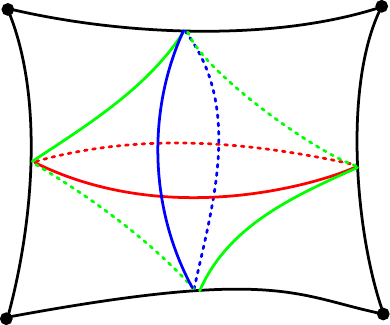}
    \caption{$2$-simplex of $\KK(S_{0,4})$}
  \end{subfigure}
  \caption{Ideal triangulations are obtained from simplices of the curve complex on $ S_{0,4} $ and $ S_{1,1}$ by `pushing' the curves parallel to themselves until they hit the punctures.\label{fig:farey_top_downstairs}}
\end{figure}

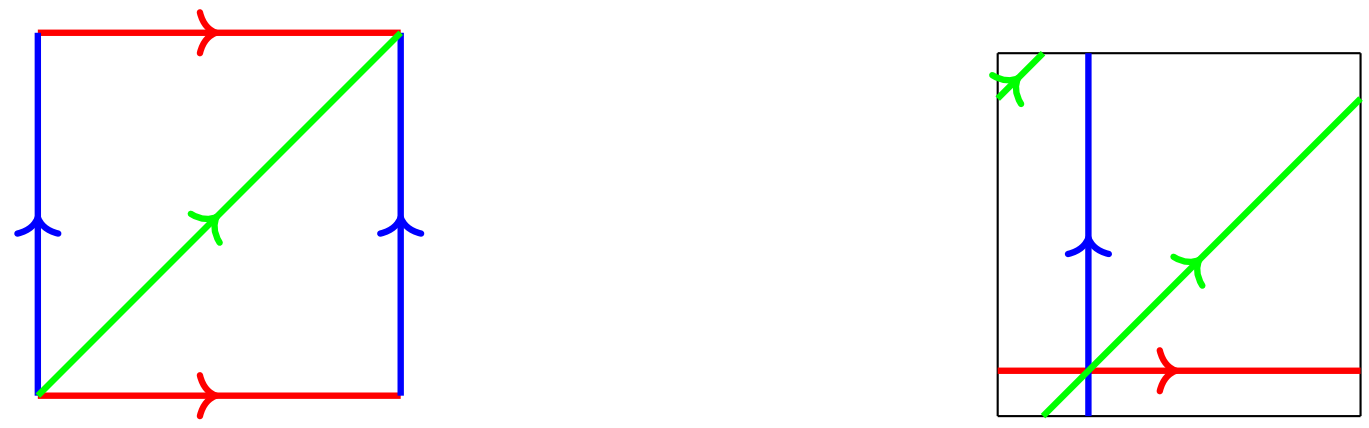
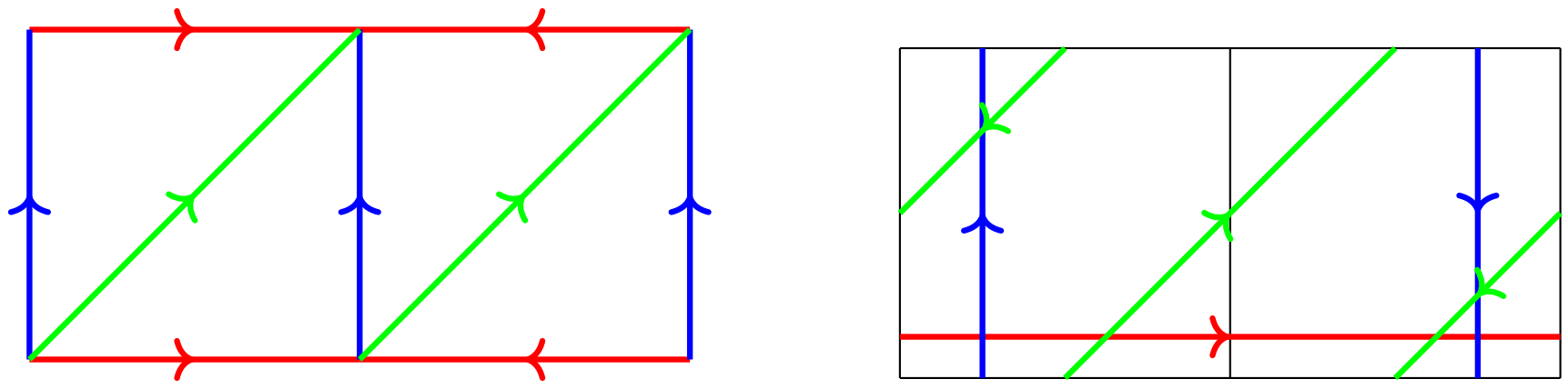
\begin{figure}
  \begin{subfigure}{.5\textwidth}\centering
    \begin{tikzpicture}[scale=.6]
      \begin{scope}[very thick,decoration={
          markings,
          mark=at position 0.5 with {\arrow{>}}}
          ]
          \draw[postaction={decorate}, red] (0,0)--(4,0);
          \draw[postaction={decorate}, blue] (4,0)--(4,4);
          \draw[postaction={decorate}, red] (0,4)--(4,4);
          \draw[postaction={decorate}, blue] (0,0)--(0,4);
          \draw[postaction={decorate}, green] (0,0)--(4,4);
      \end{scope}
    \end{tikzpicture}
    \caption{Ideal triangulation of $S_{1,1}$.}
  \end{subfigure}%
  \begin{subfigure}{.5\textwidth}\centering
    \begin{tikzpicture}[scale=.6]
      \draw (0,0)--(4,0);
      \draw (4,0)--(4,4);
      \draw (0,4)--(4,4);
      \draw (0,0)--(0,4);
      \begin{scope}[very thick,decoration={
          markings,
          mark=at position 0.5 with {\arrow{>}}}
          ]
          \draw[postaction={decorate}, red] (0,.5)--(4,.5);
          \draw[postaction={decorate}, blue] (1,0)--(1,4);
          \draw[postaction={decorate}, green] (.5,0)--(4,3.5);
          \draw[postaction={decorate}, green] (0,3.5)--(.5,4);
      \end{scope}
    \end{tikzpicture}
    \caption{$2$-simplex of $\KK(S_{1,1})$}
  \end{subfigure}\\%
 \begin{subfigure}{.5\textwidth}\centering
    \begin{tikzpicture}[scale=.6]
      \begin{scope}[very thick,decoration={
          markings,
          mark=at position 0.5 with {\arrow{>}}}
          ]
          \draw[postaction={decorate}, red] (0,0)--(4,0);
          \draw[postaction={decorate}, red] (8,0)--(4,0);
          \draw[postaction={decorate}, red] (0,4)--(4,4);
          \draw[postaction={decorate}, red] (8,4)--(4,4);
          \draw[postaction={decorate}, blue] (0,0)--(0,4);
          \draw[postaction={decorate}, blue] (4,0)--(4,4);
          \draw[postaction={decorate}, blue] (8,0)--(8,4);
          \draw[postaction={decorate}, green] (0,0)--(4,4);
          \draw[postaction={decorate}, green] (4,0)--(8,4);
      \end{scope}
    \end{tikzpicture}
    \caption{Ideal triangulation of $S_{0,4}$.}
  \end{subfigure}%
  \begin{subfigure}{.5\textwidth}\centering
    \begin{tikzpicture}[scale=.6]
      \draw (0,0)--(8,0);
      \draw (0,4)--(8,4);
      \draw (0,0)--(0,4);
      \draw (4,0)--(4,4);
      \draw (8,0)--(8,4);
      \begin{scope}[very thick,decoration={
          markings,
          mark=at position 0.5 with {\arrow{>}}}
          ]
          \draw[postaction={decorate}, red] (0,.5)--(8,.5);
          \draw[postaction={decorate}, blue] (1,0)--(1,4);
          \draw[postaction={decorate}, blue] (7,4)--(7,0);
          \draw[postaction={decorate}, green] (2,0)--(6,4);
          \draw[postaction={decorate}, green] (2,4)--(0,2);
          \draw[postaction={decorate}, green] (8,2)--(6,0);
      \end{scope}
    \end{tikzpicture}%
    \caption{$2$-simplex of $\KK(S_{0,4})$}
  \end{subfigure}
  \caption{Slope curves on the universal covers of $ S_{0,4} $ and $ S_{1,1}$ give either ideal triangulations or simplices of the curve complex, depending on whether they intersect the lifts of the punctures.\label{fig:farey_top_upstairs}}
\end{figure}

\begin{prp}\label{prp:farey}
  The Farey graph has two important topological interpretations.
  \begin{itemize}
    \item Let $ S =S_{1,1} $. Then $ \mc{F} $ parameterises ideal triangulations of $ S $ (triples of ideal arcs on $S$ which are mutually disjoint up to isotopy and which
          intersect only at the puncture), and $ \mc{F}^* $ is identified with $ \KK(S) $.
    \item Let $ S = S_{0,4} $. Then $ \mc{F} $ parameterises ideal triangulations of $ S $ (sextuples of ideal arcs on $S$ which are mutually disjoint up to isotopy and which
          intersect only at the punctures), and $ \mc{F}^* $ is identified with $ \KK(S) $. \qed
  \end{itemize}
\end{prp}
For some additional detail see \S 1.2 of Akiyoshi, Sakuma, Wada, and Yamashita~\cite{akiyoshi}, but we include some suggestive pictures in \zcref{fig:farey_top_downstairs} (on the surfaces)
and in \zcref{fig:farey_top_upstairs} (on fundamental domains for the surfaces) from which the proof of \zcref{prp:farey} can actually be deduced using some classic number
theory like \zcref{lem:farey_closeness} below.

These combinatorial data on $ \KK(S) $ can be `rigidified'
to give combinatorial hyperbolic decompositions of both punctured torus bundles and two-bridge link complements; this was originally observed by J\o{}rgensen~\cite{jorgensen03}
with additional important contributions by Minsky~\cite{minsky99}, Lackenby~\cite{lackenby03}, and Gu\'eritaud and Futer~\cite{gueritaud06}. For a detailed explanation of many
aspects of this work and a thoughtful survey of other related work see Akiyoshi, Sakuma, Wada, and Yamashita~\cite{akiyoshi}.

\begin{figure}
  \begin{subfigure}[t]{.48\textwidth}%
    \centering
    \includegraphics[width=\textwidth]{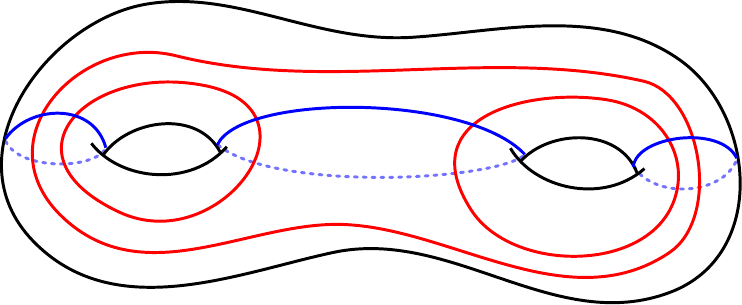}
    \caption{Two dual $2$-simplices.\label{fig:duals}}
  \end{subfigure}\hfill%
  \begin{subfigure}[t]{.48\textwidth}%
    \centering
    \includegraphics[width=\textwidth]{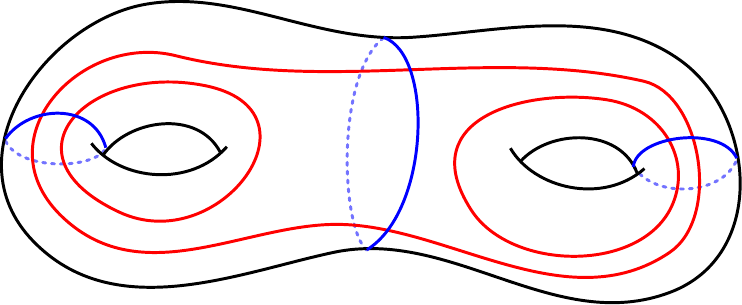}
    \caption{A $2$-simplex containing a separating curve, and a pseudodual $2$-simplex to it.\label{fig:duals2}}
  \end{subfigure}
  \caption{Duality and pseudoduality in $ \KK(S_{2,0})$.}
\end{figure}

We will also need this related definition:
\begin{defn}\label{defn:dual_simplex}
  Two simplices $ \Delta, \Delta' \in \KK^{(2)}(\Sigma)$ are \df{dual} if, whenever $ \alpha \in \Delta $ and $ \alpha' \in \Delta' $ are isotopy classes
  of curves, $ i(\alpha,\alpha') \in \{0,1\} $. Thus two $2$-simplices $ \{\alpha,\beta,\gamma\} $ and $ \{\alpha',\beta',\gamma'\} $ in $ \KK^{(2)}(S_{2,0}) $ are dual if, up to relabelling, they have
  the following intersection number matrix:
  \begin{displaymath}\setlength\arraycolsep{2pt}
    \begin{matrix}
                    & \alpha & \beta & \gamma\\\addlinespace
          \alpha'   & 0 & 1 & 1\\
          \beta'    & 1 & 0 & 1\\
          \gamma'   & 1 & 1 & 0
    \end{matrix}
  \end{displaymath}
  For example see \zcref{fig:duals}.
\end{defn}

Observe that two $2$-simplices can be dual only if the vertices of both consist of non-separating curves. Indeed, if a curve $\nu$ is separating then every other simple closed curve has even intersection number
with $\nu$. In our constructions later we will see non-separating curves only at the final stage of the process, and for this we will need an appropriate analogue of dual simplex.

\begin{defn}\label{defn:pseudodual_simplex}
  Let $ \Delta = \{ \alpha, \beta, \gamma \} \in \KK^{(2)}(\Sigma) $ be such that $\gamma$ is separating (hence both $ \alpha$ and $\beta$ are nonseparating). A \df{pseudodual} simplex to $\Delta$
  is a simplex $ \Delta' = \{ \alpha', \beta', \gamma' \} \in \KK^{(2)}(\Sigma) $ of non-separating curves which give the following intersection number matrix
  \begin{displaymath}\setlength\arraycolsep{2pt}
    \begin{matrix}
                    & \alpha & \beta & \gamma\\\addlinespace
          \alpha'   & 0 & 1 & 0\\
          \beta'    & 1 & 0 & 0\\
          \gamma'   & 0 & 0 & 2
    \end{matrix}
  \end{displaymath}
  For example see \zcref{fig:duals2}.
\end{defn}
Note that duality is symmetric, but pseudoduality is not.

\subsection{Heegaard diagrams and disc complexes}\label{sec:heegaard}
Maximal simplices of $\KK(S_{2,0})$ have $3$ vertices and so can be viewed as Heegaard diagrams as soon as a choice of identification of $ S_{2,0} $
with the boundary of a genus two handlebody $ \mc{H}^+ $ is fixed.

\begin{figure}
  \centering
  \includegraphics[width=.5\textwidth]{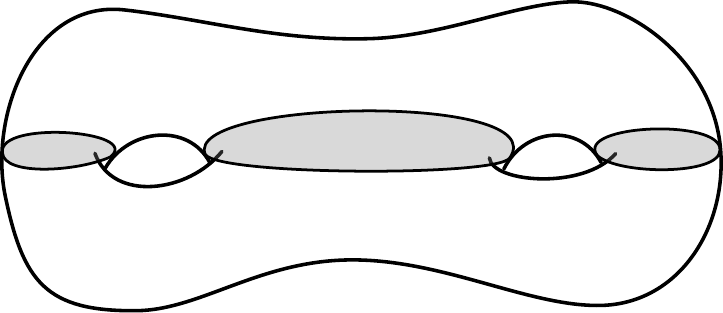}
  \caption{The simplex $ \Delta_0 \in \KK^{(2)}(\Sigma) $.\label{fig:deltanought}}
\end{figure}

With this in mind, we will fix a genus two handlebody $ \mc{H}^+ $ and its boundary $ \Sigma = \partial \mc{H}^+ $. It is well-known, but technically difficult to prove, that
the curve complex $ \KK(\Sigma) $ is a connected topological space (for one proof see Masur and Minsky~\cite[Lemma~2.1]{masur99}). We will choose a `canonical' basepoint
for the curve complex:
\begin{defn}\label{defn:deltanought}
  Let $\Delta_0 $ be the $2$-simplex shown in \zcref{fig:deltanought} which represents the unknot.
\end{defn}
The connectedness implies that, for any simplex $ \Delta \in \KK^{(2)}(\Sigma) $ (which we interpret as representing an arbitrary tunnel number $1$ knot in an arbitrary
ambient manifold) there exists a path of simplices $ \Delta_0, \Delta_1, \ldots, \Delta_n = \Delta $ such that for all $i$, $ \Delta_{i} $ and $ \Delta_{i+1} $ intersect
either in a single vertex or in an edge. If $ \Delta_n $ has vertices consisting of only nonseparating curves then we may assume that this path is geodesic and consists
of simplices whose vertices are nonseparating curves (see \cite[Proposition~2.4]{hovland24} and \cite[Corollary~3.3]{hamenstadt14}). If $ \Delta_n $ has a vertex which is
separating then we may assume $ \Delta_{n-1} $ does not contain this vertex and hence that the subgeodesic $ [\Delta_0, \Delta_{n-1}] $ consists only of simplices whose
vertices are nonseparating curves. We shall refer to a choice of a path as a \df{ns-geodesic} between $ \Delta_0 $ and $ \Delta_n $. (It is not clear to the author
whether one can always choose an ns-geodesic to be geodesic when $ \Delta_n $ has a vertex which is separating.)

If two adjacent simplices $ \Delta_i $ and $ \Delta_{i+1} $ intersect in a single vertex then we will call this vertex a \df{local cut vertex} for the ns-geodesic.
It turns out that the connected components of the topological space $ (\bigcup_{i=0}^n \Delta_i) \setminus C $, where $C$ is the set of all local cut vertices, correspond to sets of simplices
parameterising tunnel number $1$ links in different ambient manifolds. Our goal is to study the path $ (\Delta_i)_{i=0}^n $ by splitting it up into the components arising from
deleting local cut vertices in this way. In \zcref{sec:fan_edge} we will study the individual connected components (i.e.\ maximal subpaths of $(\Delta_i)_{i=0}^n $ such that consecutive
simplices intersect in an edge), and in \zcref{sec:fan_vertex} we will study the joins between these components (i.e.\ the places in the path where two consecutive simplices are joined
at a vertex alone).

The fact that the path is geodesic for $ [\Delta_0,\Delta_{n-1}]$ is important, as it implies that the words in the group corresponding to vertices of $ \Delta_k $ are longer than the words
corresponding to vertices of $ \Delta_{\ell} $ for $ \ell > k $. Indeed, by a result of Lickorish and Hempel~\cite{hempel01,lickorish62} (see also Watanabe~\cite{watanabe17}),
intersection number of curves coarsely increases with distance in the curve complex. We obtain nested \emph{topological} surfaces without this, since we construct them explicitly
via topological moves---the point is that, if we take a bad path, some of these moves might `cancel out' up to isotopy (like following a left-hand twist with a right-hand twist).
In any case, to obtain nested groups we need the word lengths in the generators $X$ and $Y$ of $ \pi_1(\mc{H}^+)$ to strictly increase as
we move from $ \Gamma_k $ to $ \Gamma_{k-1} $, and these word lengths are essentially sums of intersection numbers with $ \Delta_0 $ via curve coding. As we walk along a fan of
simplices intersecting in a single vertex, the distance in the curve complex does not
increase, but word length still increases as we are writing down Farey words and torus words (which we will define later in the paper) which do not contan any relation as the `shortest relation' in the group
corresponds to the compresson disc which has a longer word than any words seen in any $\Delta_k $ for $ k > n $). Actually, our strategy is to give recurrence relations that
show that each group $ \Gamma_k $ is generated by words in $ \Gamma_{k+1} $: this shows that the groups are nested without complicated analysis of the combinatorial group theory.

\section{The Goeritz group and the representation space}\label{sec:goeritz}
As described at the end of the previous section, our main construction takes a tunnel number $1$ manifold $M$ arising from a Heegaard diagram $ \Delta $
and produces a filtration depending on a choice of ns-geodesic $ [\Delta_0,\Delta] \subset \KK(S_{2,0}) $. In this section, which is not necessary for the
logical development of the remainder of the paper, we say some words about the number of distinct choices for a diagram $\Delta $ given the manifold $M$,
as well as the interaction between $\pi_1(M)$ and various subgroups of $ \Mod(S_{2,0}) $.

\subsection{The Goeritz group}
\begin{figure}
  \centering
  \labellist
  \small\hair 2pt
  \pinlabel {$q$} [b] at 410 263
  \pinlabel {$r$} [b] at 410 91
  \pinlabel {$s$} [r] at 613 325
  \pinlabel {$t$} [t] at 633 118
  \endlabellist
  \includegraphics[width=.8\textwidth]{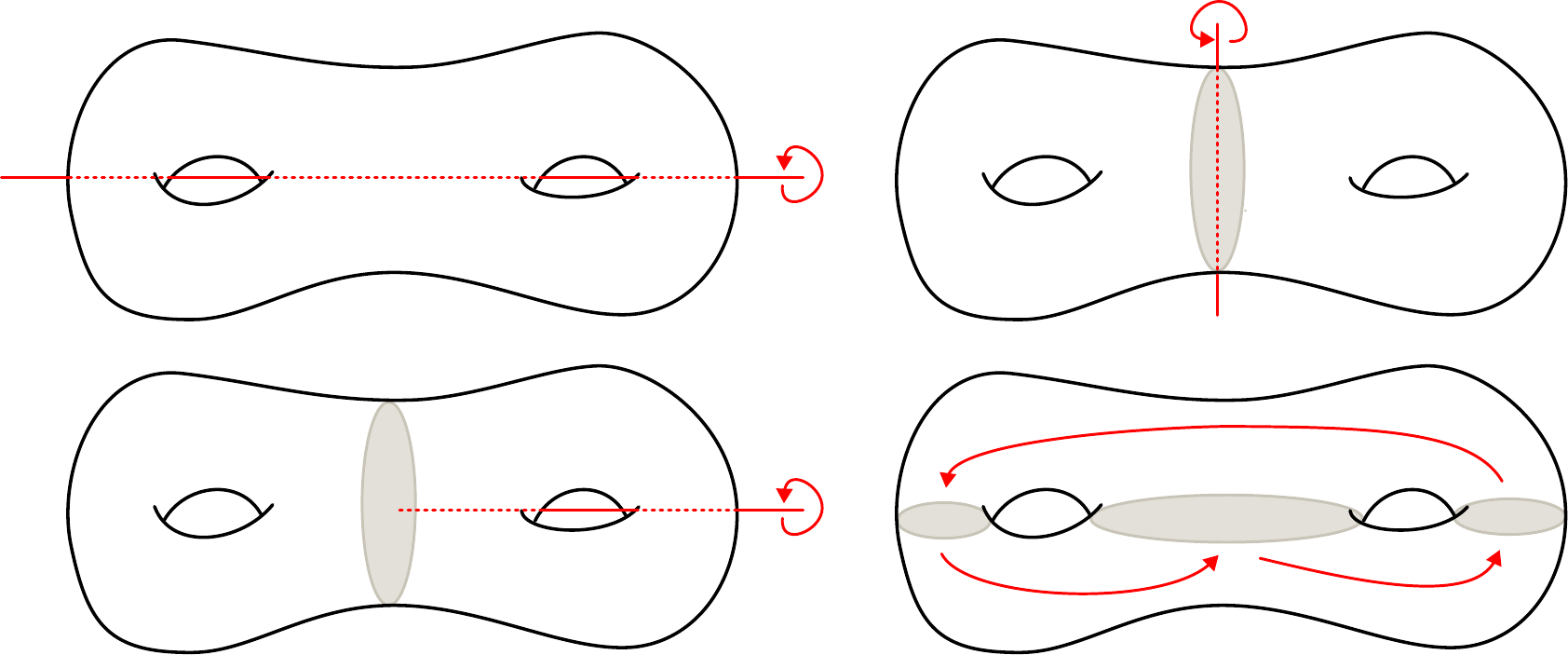}
  \caption{The four Akbas--Scharlemann generators of $ \Goer(\IS^3,\Sigma) $. The elements $ q $, $ r $, and $s$ are order $2$, and the element $t$ is order $3$.\label{fig:goeritz_gens}}
\end{figure}
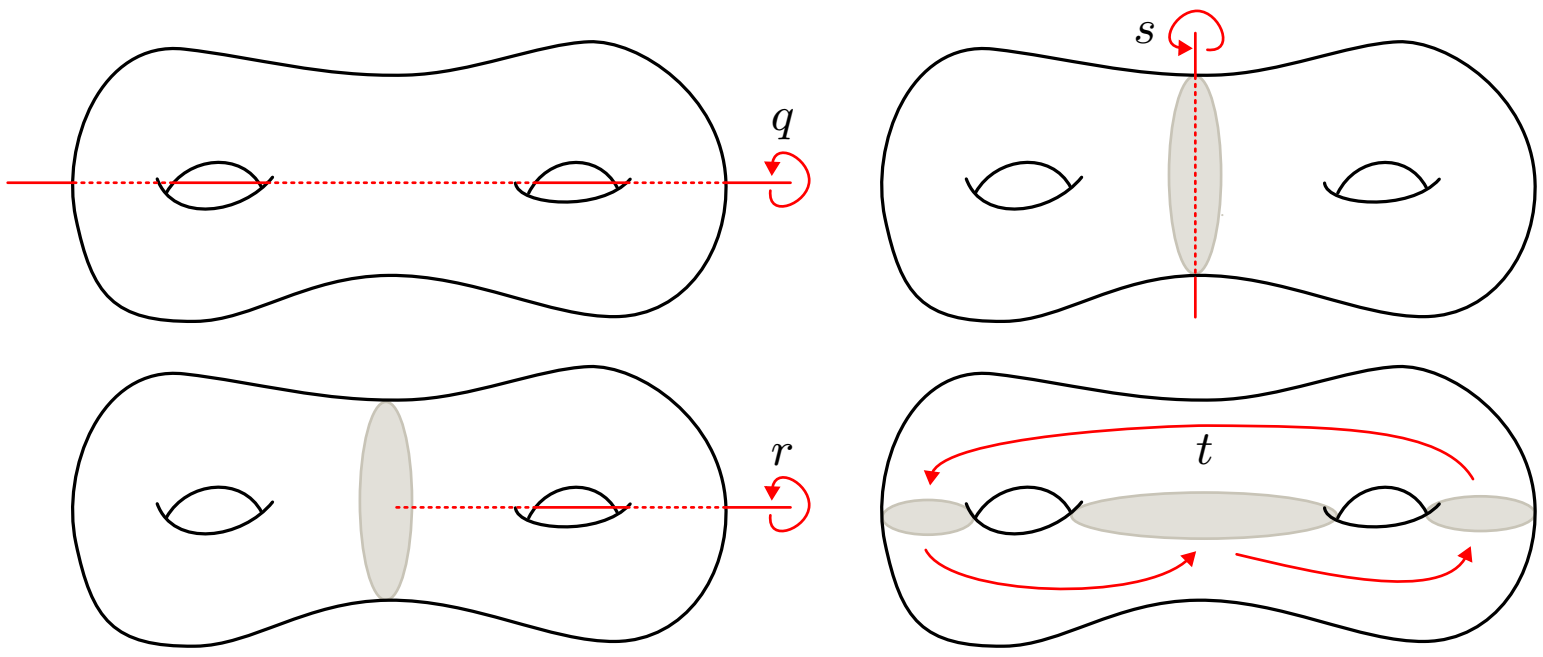

There is a natural automorphism group of the Heegaard splitting of a closed $3$-manifold $A$ by the genus $2$ surface $\Sigma$:
\begin{defn}\label{defn:goeritz}
  The \df{Goeritz group} of $ (A,\Sigma) $, which we will denote $ \Goer = \Goer(A,\Sigma) $, is the group of isotopy classes of
  homeomorphisms of $ A $ which preserve the two handlebodies $ \mc{H}^+ $ and $ \mc{H}^- $.
\end{defn}
The Goeritz group of $ \IS^3 $ (with its unique up to isotopy genus $2$ Heegaard splitting), first introduced by Goeritz~\cite{goeritz33}, is well-understood and explicit presentations
have been given by Akbas and Scharlemann~\cite{akbas08,scharlemann04,cho08}:
\begin{lem}\label{lem:goeritz_gens}
  One set of generators of $ \Goer(\IS^3, \Sigma) $ (where $\Sigma$ is the unique-up-to-isotopy genus $2$ Heegaard surface in $ \IS^3$) is $ \{q,r,s,t\} $ as shown
  in \zcref{fig:goeritz_gens}. \qed
\end{lem}

In addition, Goeritz groups of reducible splittings were proved to be finitely presented by Cho and Koda~\cite{cho19}. It is an open question whether the Goeritz group of a general Heegaard splitting is
even finitely generated.

Clearly elements of $ \Goer $ restrict to elements of $ \Mod(\mc{H}^+) $, the mapping class group of $ \mc{H}^+ $, and then to elements of $ \Mod(\Sigma) $. The induced
map $ \Mod(\mc{H}^+) \to \Mod(\Sigma) $ is an inclusion~\cite[Lemma~3.1]{hensel}, there is a similar inclusion $ \Mod(\mc{H}^-) \to \Mod(\Sigma) $, and so $ \Goer = \Mod(\mc{H}^+) \inter \Mod(\mc{H}^-) $
also embeds naturally into $ \Mod(\Sigma) $. There is a natural map $ \Mod(\mc{H}^+) \to \Out(\pi_1(\mc{H}^+)) $.

A \df{meridian} is a simple closed curve on $ \Sigma $ which bounds a disc in $ \mc{H}^+ $---note, this is the `opposite side' of $ \Sigma $ to the disc complex $ \mc{D} $.
The image of $\Mod(\mc{H}^+) \to \Mod(\Sigma) $ is characterised by the following property:
\begin{lem}[{\cite[Corollary~5.11]{hensel}}]\label{lem:image_of_modH}
  An element of $ \Mod(\Sigma) $ is the restriction of an element of $ \Mod(\mc{H}^+) $ if
  and only if it preserves the set of meridians (and hence permutes embedded discs in $ \mc{H}^+ $). \qed
\end{lem}

It was shown by Luft~\cite{luft78} that the kernel of $ \Mod(\mc{H}^+) \to \Out(\pi_1(\mc{H}^+)) $ is the \df{twist group} $\Tw$ generated by Dehn twists around
meridians. We remark that the twist group $ \Tw $ is not finitely generated~\cite[Theorem~1.1]{mccullough85}. In summary, we have the following groups:
\begin{equation}\label{eq:goeritz}\begin{tikzcd}[row sep=scriptsize]
  1\arrow[d]\\
  \Tw \arrow[d] & \Goer \arrow[dl,hook']\arrow[dr, hook] \\
  \Mod(\mc{H}^+) \arrow[d] \arrow[r, hook] & \Mod(\Sigma) & \Mod(\mc{H}^-) \arrow[l, hook']\\
  \Out(\pi_1(\mc{H}^+))\arrow[d]\\1
\end{tikzcd}\end{equation}

Finally, there is a natural map $ \mathrm{Glue} $ which sends a `flag' consisting of a $2$-simplex $ \Delta \in \KK^{(2)}(\Sigma) $ and a choice of vertex $ \alpha \in \Delta $ to
a pair $( \mf{k},\tau )$ consisting of (an isotopy class of) a tunnel number $1$ link (in some ambient $3$-manifold $A$) together with a tunnel $ \tau $, where $ A \setminus \mf{k} $ is obtained
by gluing a $2$-handle to $ \mc{H}^+ $ along $ \alpha $ and where $ \tau $ intersects that $2$-handle transversely once.

\subsection{The presentation space}
We recall some terminology following Johannson~\cite[Chapter~V]{johannson95}. If $G$ is a group, then for every pair $ F, H $ of test groups we can
consider two interesting spaces: the \df{presentation space}
\begin{displaymath}
  P(F,G) = (F \twoheadrightarrow G)/\Aut(F),
\end{displaymath}
and the \df{representation space}
\begin{displaymath}
  R(G,H) = \Inn(H) \backslash (G \hookrightarrow H)
\end{displaymath}
(here, $ \Inn(H) $ denotes $H$ acting on itself by conjugacy). In the current setting, we take $ G = \pi_1(M) $ for $M$ a hyperbolic tunnel number $1$ link complement, $F = \pi_1(\mc{H}^+) $,
and $ H = \PSL(2,\IC) $.

We will deal with the representation space in the next subsection.
The presentation space $ P(\pi_1(\mc{H}^+, \pi_1(M)) $ can be viewed as indexing the different Heegaard diagrams which give the same manifold group
up to group isomorphism. For prime knots, this determinacy is enough to determine the isotopy type~\cite{gordon89}---our primary motivations are from the point of view of the groups
rather than the manifolds anyway, so even in the case of links this space is still giving us all the desired information. Johansson has proved the following rigidity theorem:
\begin{thm}[{Johannson~\cite[Theorem~32.17]{johannson95}}]\label{thm:johannson}
  Let $M$ be a Haken $3$-manifold without nontrivial essential Stallings fibration. Then $M$ has finitely many Heegaard surfaces of Euler characteristic $ n $ for all $ n \leq -2 $.\
  In particular, $ P(F_n, M) $ is finite. \qed
\end{thm}
Any hyperbolic link complement is Haken (see e.g.\ the \textsc{Note} in \S 4.10 of Thurston~\cite{thurstonN}). Further, tunnel number $1$ link complements generically do not fibre over the circle~\cite{dunfield06}.
Thus, \zcref{thm:johannson} gives us a generic finiteness result for the number of different $2$-simplices $ \Delta \in \KK^{(2)}(S_{2,0}) $ which produce a fixed hyperbolic tunnel number $1$ link
manifold.
\begin{cor}\label{cor:johannson_cusps}
  For a generic hyperbolic tunnel number $1$ link complement $ A \setminus \mf{k} $, where `generic' takes one of the meanings in \cite{dunfield06}, there are only finitely many Heegaard diagrams $ \Delta $
  in the fibre of $ \mathrm{Glue} $ above $ (\mf{k},\tau) $, for each tunnel $ \tau $. \qed
\end{cor}
In fact, Johannson shows that this finite set can be computed~\cite[Appendix]{johannson95}.

Note that $ P(F_n, M) $ is a quotient by a group that is defined only in terms of abstract algebra, $ \Aut(F\{X,Y\}) $. We can replace this with
a group that arises naturally from the topology of $ \mc{H}^+ $. It is well-known that $ \Aut(F\{X,Y\}) $ is generated by the automorphisms $ S,R,U $ defined on generators by
\begin{displaymath}
  S(X,Y) = (Y,X),\; R(X,Y) = (X^{-1},Y),\; \text{and}\; U(X,Y) = (XY, Y)
\end{displaymath}
and that $ \Out(\pi_1(\mc{H}^+)) $ can be identified with $ \GL(2,\IZ) $~\cite{nielsen17} (c.f.~\cite[Corollary~N5 in~\S 3.6]{magnusCGTe1}) where the projection map is given by
\begin{displaymath}
  S \mapsto \begin{bmatrix} 0 & 1 \\ 1 & 0 \end{bmatrix},\; R \mapsto \begin{bmatrix} -1 & 0 \\ 0 & 1 \end{bmatrix},\; \text{and}\; U \mapsto \begin{bmatrix} 1 & 0 \\ 1 & 1 \end{bmatrix}.
\end{displaymath}
Observe the similarity between these generators and the generators of $ \Goer(\IS^3,\Sigma) $ given in \zcref{lem:goeritz_gens} above. A short computation shows that, if $ \Phi : \Goer \to \GL(2,\IZ) $
is the map induced by the diagram in~\eqref{eq:goeritz}, then $ \Phi(r) = R $, $ \Phi(s) = S $, and $ \Phi(t) = R^{-1} S^{-1} U $; i.e.\ the images of $ r,s,t $ generate $ \Out(\pi_1(\mc{H}^+)) $
and hence the image of $ \Goer$ generates $ \Out(\pi_1(\mc{H}^+)) $.

This also shows that the induced action of $ \Goer(\IS^3,\Sigma) $ on $ \Hom(\pi_1(\mc{H^+}), \PSL(2,\IC)) $, defined by restriction of homeomorphisms to $ \mc{H}^+$, generates the entire
set of automorphisms, and looking at the full action of $ \Mod(\mc{H}^+) $ gives no additional action.

\subsection{The representation space}
In the representation space $ R(\pi_1(M), \PSL(2,\IC)) $ (i.e.\ the $\PSL(2,\IC)$-\df{character variety} of $ \pi_1(M)$, viewed as a set with all algebraic structure forgotten), we
are mainly interested in the subspace of discrete representations. We know (by Mostow-Prasad rigidity~\cite[Theorem~3.32]{matsuzaki_taniguchi}) that inside this group there lies a unique `excellent'
representation: one which is discrete, faithful, and finite covolume. In general there may be other discrete representations, for example orbifold structures which lie on the same
algebraic component of the representation space as the excellent representation~\cite{hilden00}. One explicit example is the three-fold quotient of the Borromean rings, whose
holonomy group lies in $ R(\pi_1(\IS^3 \setminus \mf{l}), \PSL(2,\IC)) $ where $ \mf{l} $ is the link \href{https://katlas.org/wiki/L6a2}{$\mathrm{L6a2}$}~\cite[Theorem~3.6]{elzenaar26pr}.
In this subsection we will obtain some coarse information about the character variety of $ F\{X,Y\} $ around a fixed excellent representation of a tunnel number $1$ link group.

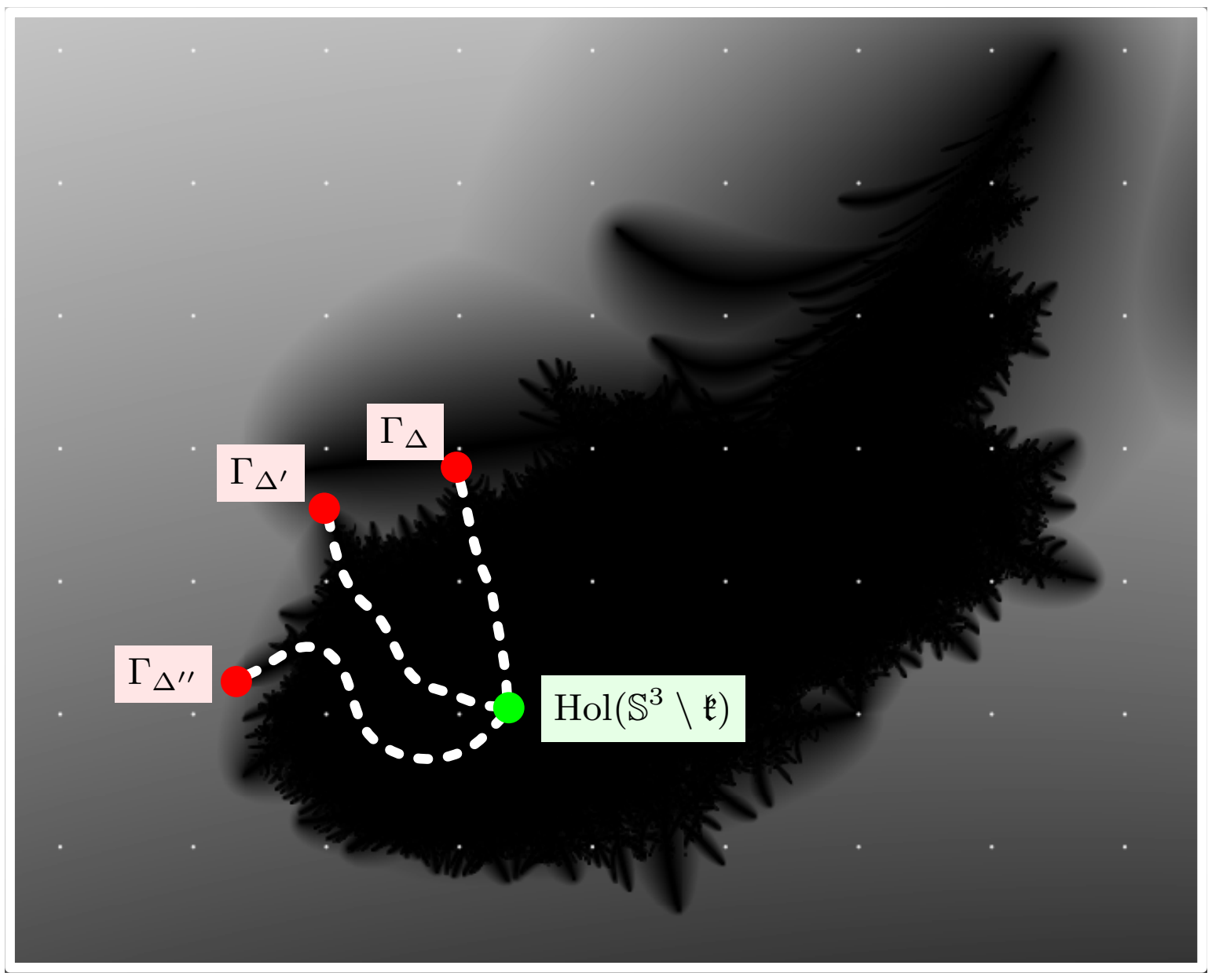
\begin{figure}
  \centering
  \labellist
  \small\hair 2pt
  \pinlabel {\colorbox{green!10}{$\Hol(\IS^3 \setminus \mf{k})$}} [l] at 76 38
  \pinlabel {\colorbox{red!10}{$\Gamma_{\Delta}$}} [br] at 64 73
  \pinlabel {\colorbox{red!10}{$\Gamma_{\Delta'}$}} [br] at 44 67
  \pinlabel {\colorbox{red!10}{$\Gamma_{\Delta''}$}} [r] at 31 43
  \endlabellist
  \includegraphics[width=.8\textwidth]{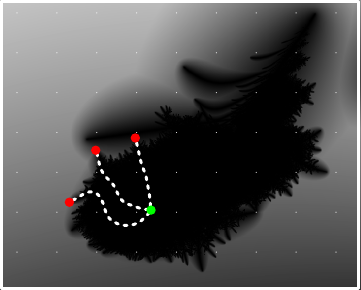}
  \caption{A schematic of the cusps corresponding to a fixed hyperbolic tunnel number $1$ link.\label{fig:defspace}}
\end{figure}

It is a result of the work of Ahlfors, Bers, Kra, and Maskit that if $ G $ is a genus $2$ Schottky group with $ \IH^3/G = \mc{H}^+ $, then the quasiconformal deformation space
$ \QH(G) $---which can be interpreted as the \df{Schottky space} $ \mc{S}_2 $ of marked, complete, convex cocompact hyperbolic metrics on $ \mc{H} $---is biholomorphically
identified with the quotient $ \Teich(\Sigma)/\Tw $~\cite[Example~5.28]{matsuzaki_taniguchi}.

Now by work of Keen, Maskit, and Series~\cite{keen91}, maximal cusps on the boundary of genus $2$ Schottky space $ \mc{S}_2 $ are entirely
determined as marked Kleinian groups up to $\PSL(2,\IC)$-conjugacy by the set of three loops which are pinched to parabolics; since
$ \mc{S}_2 $ is a quotient by $ \Tw $, the maximal cusps are in bijection with the orbits of $2$-simplices of $ \KK(\Sigma) $ under
the twist group $ \Tw$. Since all the groups are marked with generators $ X $ and $ Y $, the action of $ \Out(\mc{H}^+) $ (which we have seen is
generated entirely by the images of elements of the Goeritz group) then further permutes the different maximal cusps in the character variety which correspond to different choices
of generators for the same group. (For instance, this action gives rise to the $4$-fold symmetry of the Riley slice~\cite[\S 2.2]{ems22b}.)

In other words:
\begin{prp}
  Consider the action of $ \Mod(\Sigma) $ on flags consisting of a $2$-simplex $ \Delta \in \KK^{(2)}{\Sigma} $ and a vertex $ \alpha \in \Delta $. Two flags have the
  same image knot-tunnel pair under $ \mathrm{Glue} $ if and only if they lie in the same orbit of the subgroup of $ \Mod(\Sigma) $ generated by $ \Tw $
  and $ \Out(\pi_1(\mc{H}^+)) = \Goer(\IS^3,\Sigma) \subset \Mod(\mc{H}^+) $.
  \qed
\end{prp}
For a practical example of this in $ (1;2)$-compression body space, see~\cite[Example 2.3]{elzenaar25h}.

Now the corollary to Johannson's rigidity theorem, \zcref{cor:johannson_cusps}, tells us that $ \mathrm{Glue} $ is generically finite-to-one onto a given hyperbolic link $ \mf{k} $.
This means that $ \mf{k} $ is associated to only finitely many $2$-simplices $ \Delta $, and hence only finitely many cusp groups $ \Gamma_\Delta $. For each cusp group $ \Gamma_\Delta $, there is a path
in the character variety corresponding to a branch of the inverse image of $ [-2,2] $ under $ \tr \gamma $ (where $\gamma$ is the curve representing the loop along
which the $2$-handle is being glued) which joins the cusp to the corresponding link group, see \zcref{fig:defspace}. It has been conjectured~\cite[Preface, \S 3]{akiyoshi} that every
point on this curve is the holonomy group of a hyperbolic cone manifold that realises `drilling the unknotting tunnel'; this conjecture has been proved for the
special case of two-bridge links~\cite[Remark~3.10]{elzenaar25ej}.

\section{The fan of simplices intersecting in an edge}\label{sec:fan_edge}
Let $ (A,\Sigma) $ be a genus two Heegaard splitting of a closed $3$-manifold $A$. Let $ \mc{H}^+(A,\Sigma) $ and $ \mc{H}^-(A,\Sigma) $ be the two handlebodies of the splitting.
For the moment these two handlebodies are indistinguishable, but when we apply our discussion to the case of splittings arising from tunnel number $1$ links, the link will lie in
the negative handlebody.

\subsection{Dual graphs and the Farey graph}
We follow some of the machinery and terminology of Cho and McCullough~\cite{cho09}.
\begin{defn}\label{defn:disc_complex}
  The \df{disc complex} of $ (A,\Sigma) $, denoted $ \mc{D}(A,\Sigma) $, is the simplicial complex supported on the set of isotopy classes of embedded \emph{nonseparating} discs in $ \mc{H}^-(A,\Sigma) $,
  such that a set forms a simplex if the isotopy classes have representatives which are disjoint. By taking boundaries of the discs, $ \mc{D}(A,\Sigma) $ is a subcomplex of $ \KK(\Sigma) $.
\end{defn}
We now define a dual graph to the $2$-skeleton of $ \mc{D}(A,\Sigma) $.
\begin{defn}\label{defn:lifted_cm_tree}
  Let $ \mc{D}'(A,\Sigma) $ be the first barycentric subdivision of $ \mc{D}(A,\Sigma) $. This has three kinds of vertices: (i) vertices corresponding to $2$-simplices of $\mc{D}(A,\Sigma)$,
  i.e.\ to triples of discs; (ii) vertices corresponding to $1$-simplices of $\mc{D}(A,\Sigma)$, i.e.\ to pairs of discs; and (iii) the vertices of $ \mc{D}(A,\Sigma) $. Let $ \tilde{\mc{T}}(A,\Sigma) $
  be the subcomplex of $\mc{D}'(A,\Sigma)$ obtained by deleting the open stars of the vertices of $\mc{D}(A,\Sigma)$.
\end{defn}

\begin{lem}
  The complex $\tilde{\mc{T}}(A,\Sigma)$ is a tree.
\end{lem}
\begin{proof}
  Observe that $\tilde{\mc{T}}(A,\Sigma)$ is independent of the manifold $A$, depending only on $ \mc{H}^-(A,\Sigma) $. (That is, it is a complex that is assigned to the genus $2$ handlebody in the abstract, not to the
  splitting.) Thus the claimed result is an immediate consequence of Theorem~4.1 of \cite{cho09}.
\end{proof}
Despite the independence of the tree from the splitting, we keep $(A,\Sigma)$ in the notation because later we will be gluing copies of this tree together that correspond to different splittings.

\begin{figure}
  \centering
  \includegraphics[width=\textwidth]{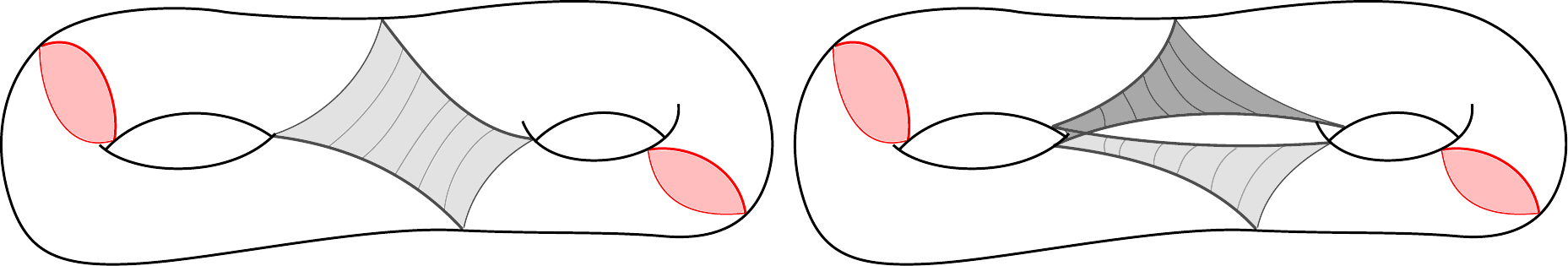}
  \caption{Rational braiding around the two red discs: the $ 1/1 $ disc is transformed into the $1/2 $ disc by a Dehn half-twist.\label{fig:cabling}}
\end{figure}

We now observe that $\tilde{\mc{T}}(A,\Sigma)$ is locally (around each vertex) a union of countably many copies of $ \mc{F} $.
\begin{cons}[Rational braiding]\label{cons:cabling}
  Let $ \Delta=\{ D,E,F \} $ be a simplex in $ \mc{D}(A,\Sigma) $. Then a simplex $ \Delta'=\{ D,E',F \} $ is said to be obtained by $p/q$-\df{rational braiding} $\Delta$
  around the edge $ \{D,F\} $ if there is an identification of the dual of the curve complex of the four-holed sphere $ \Sigma \setminus (\partial D \union \partial F) $ with the Farey tree $ \mc{F} $
  so that $ E $ is identified with $ 0/1 $ if it is non-separating and with $0/1$ if it is separating, and $ E'$ is identified with $ p/q $. See \zcref{fig:cabling} for an example.
\end{cons}
\begin{rem}
  Cho and McCullough use the word `cabling' to refer to a sequence of rational braidings, which unfortunately clashes with another older construction from knot theory~\cite[Example~15.21]{burde3e}.
\end{rem}

Since all simple closed curves on the four-holed sphere bound discs on the interior, every vertex of $\mc{F}$ corresponds to a simplex of $ \mc{D}(A,\Sigma) $ obtained
by such a rational braiding construction. This gives a way of attatching a copy of the Farey tree to every $2$-simplex of $\mc{D}(A,\Sigma)$. When this is done for all $2$-simplices, the result
is the tree $ \tilde{\mc{T}}(A,\Sigma) $.

Fix an edge $ \mc{E} $ of the curve complex, and consider the subcomplex consisting of all $2$-simplices containing $ \mc{E} $.
If $ \Delta $ and $ \Delta' $ lie in this subcomplex, then the vertices of $ \Delta' $ are obtained from the vertices of $\Delta $ by fixing two
curves and varying the third. The possible choices of third curves are indexed by the curve complex of the four-holed sphere. Since
every simple closed curve on the four-holed sphere bounds a disc, the possible motions through the subcomplex are just rational braidings, \zcref{cons:cabling},
and each produces a Heegaard diagram for a different link in the same ambient manifold. In particular, any subpath of the path $ (\Delta_i)_{i=0}^n$ consisting
of adjacent simplices intersecting in edges lies in a single disc complex $ \mc{D}(A,\Sigma) $ and so corresponds to a path in the corresponding tree
$ \tilde{\mc{T}}(A,\Sigma) $: hold two discs constant and
cable the third, then choose a different disc out of the three to cable while holding the others fixed, and so forth.

\begin{rem}
  The Cho--McCullough structure on the quotient $ \tilde{\mc{T}}(\IS^3,\Sigma)/\Goer(\IS^3,\Sigma) $ is obtained by observing that these local Farey trees are preserved by
  the Goeritz group action. In general the action of the Goeritz group of a Heegaard splitting seems to act on the tree $\tilde{\mc{T}}(A,\Sigma)$ in a very complicated way.
  Understanding the quotient $\tilde{\mc{T}}(A,\Sigma)/\Goer(A,\Sigma)$ would give a combinatorial structure on the unknotting tunnels of tunnel number $1$ links up to symmetry;
  since we only study the covering space, we will see the same link and tunnel multiple times. Because our interest is only in understanding the structure of an individual link,
  this is not a big problem (at worst, we obtain several different structures on the same link). The most general results on the space of equivalence classes of Heegaard
  splittings seem to be due to Johnson and McCullough~\cite{johnson13}.
\end{rem}

\subsection{Cutting sequences and word recursion under rational braiding}
We will now consider the following question: given a path $ \Delta_1,\ldots,\Delta_n $ of simplices indexed by a ns-geodesic in $\tilde{\mc{T}}(A,\Sigma)$,
can we write down a nested sequence of subgroups $ \Gamma_1 < \cdots < \Gamma_n $ where each $ \Gamma_i $ is a surface group which includes the loops
indexed by the vertices of $ \Delta_1 \cup \cdots \cup \Delta_i $?

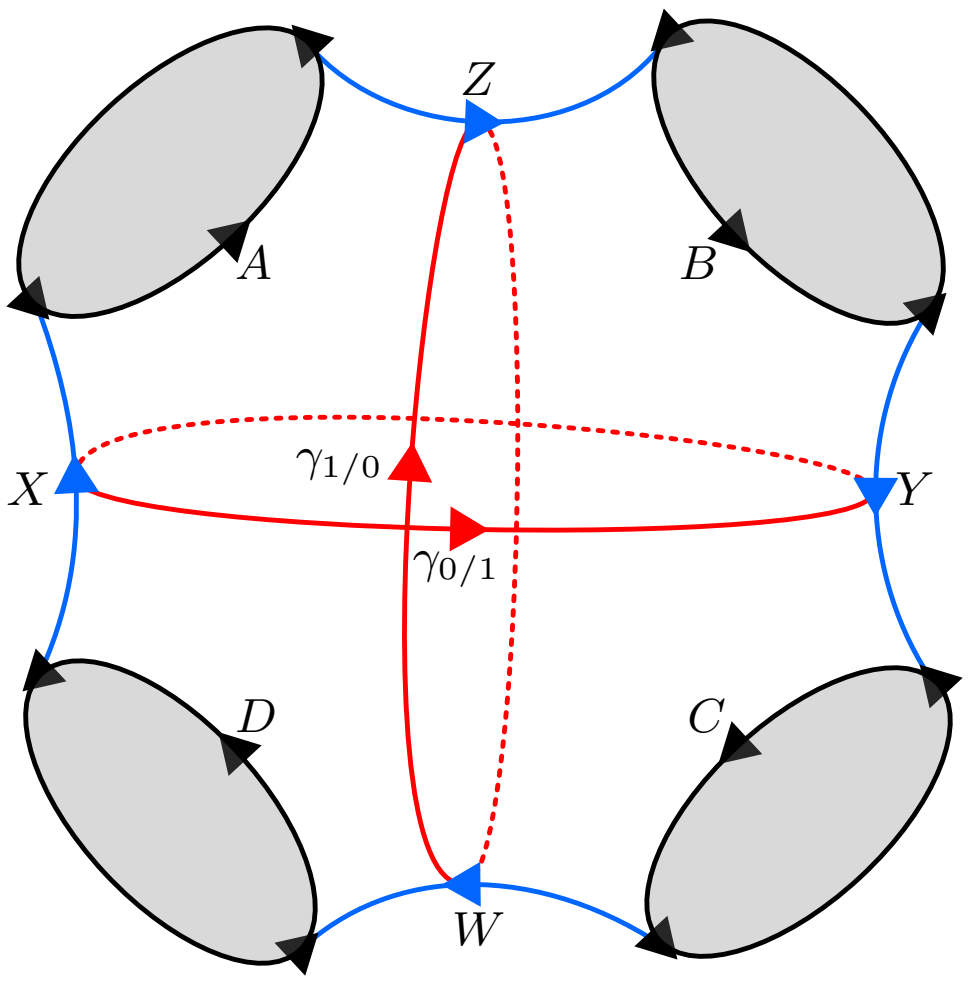
\begin{figure}
  \labellist
  \small\hair 2pt
  \pinlabel {$A$} [tl] at 60 200
  \pinlabel {$Z$} [b] at 128 236
  \pinlabel {$B$} [tr] at 195 200
  \pinlabel {$Y$} [l] at 238 132
  \pinlabel {$C$} [br] at 197 64
  \pinlabel {$W$} [t] at 128 20
  \pinlabel {$D$} [bl] at 60 64
  \pinlabel {$X$} [r] at 15 132
  \pinlabel {$\gamma_{1/0}$} [r] at 105 138
  \pinlabel {$\gamma_{0/1}$} [t] at 122 118
  \endlabellist
  \centering
  \includegraphics[width=0.5\textwidth]{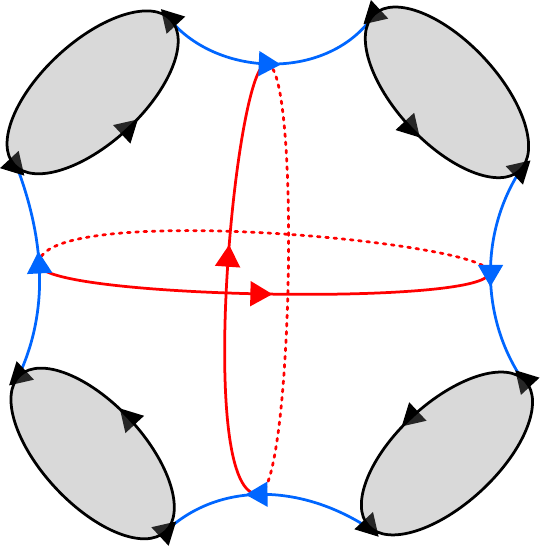}
  \caption{Loops and arcs on the four-holed sphere. The actual loops on the sphere can be written in terms of the oriented arcs $ W $, $ X $, $ Y $, and $ Z $ as follows: $ \gamma_{0/1} = XY^{-1} $; $ \gamma_{1/0} = WZ^{-1} $; $ A = XZ^{-1} $; $B = ZY^{-1} $; $C = YW^{-1} $; and $ D= WX^{-1} $.\label{fig:four_holed_sphere}}
\end{figure}

Consider a $4$-holed sphere with the loops around the punctures represented by words $ A,B,C,D \in \pi_1(S_{0,4}) $, so that $ ABCD = 1 $. We cannot represent loops on the sphere
using cutting sequences in the alphabet $ \mathfrak{P} = \{A,B,C,D\} $. Instead we define oriented arcs $ X, Y, W, Z $ on the sphere so that
the complement of the arcs is a pair of octagons. We choose the labelling to be as in \zcref{fig:four_holed_sphere}. We then choose a
$\IZ$-homology basis $ \{\gamma_{0/1}, \gamma_{1/0}\} $ so that $ \gamma_{1/0} = WZ^{-1} $ and $ \gamma_{0/1} = XY^{-1} $. This gives us `factorised' cutting
sequences for the alphabet $ \mathfrak{P} $ where each letter is written as a product of letters in $ \mathfrak{F} = \{X,Y,Z,W\} $.

\begin{defn}
  For $ (p,q) = 1 $ define $ \gamma_{p/q} $ to be the unique oriented simple closed curve (up to homotopy on $ S_{0,4}$) with
  homology $ q \gamma_{0/1} + p \gamma_{1/0} $ on $ S_{0,4} $. Define $ \Farey_{p/q} $ to be the $ \mathfrak{F}$-cutting sequence
  of the curve $ \gamma_{p/q} $.
\end{defn}
Here, `SW' stands for `\textbf{S}phere \textbf{W}ord'.

\begin{figure}
  \centering
    \begin{tikzpicture}[scale=2]
    \draw (0,0)--(0,1)--(0,2)--(1,2)--(1,1)--(1,0)--cycle;
    \node[label={[label distance=-3pt]left:$y$}] at (1,1.5) {};
    \node[label={[label distance=-3pt]right:$x$}] at (0,1.5) {};
    \node[label={[label distance=-3pt]left:$Y$}] at (1,0.5) {};
    \node[label={[label distance=-3pt]right:$X$}] at (0,0.5) {};

    \draw (1,0)--(1,1)--(1,2)--(2,2)--(2,1)--(2,0)--cycle;
    \node[label={[label distance=-3pt]left:$X$}] at (2,1.5) {};
    \node[label={[label distance=-3pt]right:$Y$}] at (1,1.5) {};
    \node[label={[label distance=-3pt]left:$x$}] at (2,0.5) {};
    \node[label={[color=red,label distance=-3pt]right:$y$}] at (1,0.5) {};

    \draw (2,0)--(2,1)--(2,2)--(3,2)--(3,1)--(3,0)--cycle;
    \node[label={[label distance=-3pt]left:$y$}] at (3,1.5) {};
    \node[label={[color=red, label distance=-3pt]right:$x$}] at (2,1.5) {};
    \node[label={[label distance=-3pt]left:$Y$}] at (3,0.5) {};
    \node[label={[label distance=-3pt]right:$X$}] at (2,0.5) {};

    \draw (3,0)--(3,1)--(3,2)--(4,2)--(4,1)--(4,0)--cycle;
    \node[label={[label distance=-3pt]left:$X$}] at (4,1.5) {};
    \node[label={[color=red, label distance=-3pt]right:$Y$}] at (3,1.5) {};
    \node[label={[label distance=-3pt]left:$x$}] at (4,0.5) {};
    \node[label={[label distance=-3pt]right:$y$}] at (3,0.5) {};

    \draw (0,2)--(0,3)--(0,4)--(1,4)--(1,3)--(1,2)--cycle;
    \node[label={[label distance=-3pt]left:$y$}] at (1,3.5) {};
    \node[label={[label distance=-3pt]right:$x$}] at (0,3.5) {};
    \node[label={[label distance=-3pt]left:$Y$}] at (1,2.5) {};
    \node[label={[label distance=-3pt]right:$X$}] at (0,2.5) {};

    \draw (1,2)--(1,3)--(1,4)--(2,4)--(2,3)--(2,2)--cycle;
    \node[label={[label distance=-3pt]left:$X$}] at (2,3.5) {};
    \node[label={[label distance=-3pt]right:$Y$}] at (1,3.5) {};
    \node[label={[label distance=-3pt]left:$x$}] at (2,2.5) {};
    \node[label={[label distance=-3pt]right:$y$}] at (1,2.5) {};

    \draw (2,2)--(2,3)--(2,4)--(3,4)--(3,3)--(3,2)--cycle;
    \node[label={[label distance=-3pt]left:$y$}] at (3,3.5) {};
    \node[label={[label distance=-3pt]right:$x$}] at (2,3.5) {};
    \node[label={[label distance=-3pt]left:$Y$}] at (3,2.5) {};
    \node[label={[label distance=-3pt]right:$X$}] at (2,2.5) {};

    \draw (3,2)--(3,3)--(3,4)--(4,4)--(4,3)--(4,2)--cycle;
    \node[label={[label distance=-3pt]left:$X$}] at (4,3.5) {};
    \node[label={[label distance=-3pt]right:$Y$}] at (3,3.5) {};
    \node[label={[label distance=-3pt]left:$x$}] at (4,2.5) {};
    \node[label={[label distance=-3pt]right:$y$}] at (3,2.5) {};

    \node[label={[label distance=-3pt]right:$X$}] at (4,0.5) {};
    \node[label={[label distance=-3pt]right:$x$}] at (4,1.5) {};
    \node[label={[color=red,label distance=-3pt]right:$X$}] at (4,2.5) {};
    \node[label={[label distance=-3pt]right:$x$}] at (4,3.5) {};


    \node[label={[label distance=-3pt]below:$w$}] at (1.5,1) {};
    \node[label={[label distance=-3pt]below:$W$}] at (0.5,1) {};
    \node[label={[label distance=-3pt]above:$z$}] at (1.5,0) {};
    \node[label={[label distance=-3pt]above:$Z$}] at (0.5,0) {};

    \node[label={[label distance=-3pt]below:$Z$}] at (1.5,2) {};
    \node[label={[label distance=-3pt]below:$z$}] at (0.5,2) {};
    \node[label={[color=red,label distance=-3pt]above:$W$}] at (1.5,1) {};
    \node[label={[label distance=-3pt]above:$w$}] at (0.5,1) {};

    \node[label={[label distance=-3pt]below:$w$}] at (1.5,3) {};
    \node[label={[label distance=-3pt]below:$W$}] at (0.5,3) {};
    \node[label={[label distance=-3pt]above:$z$}] at (1.5,2) {};
    \node[label={[label distance=-3pt]above:$Z$}] at (0.5,2) {};

    \node[label={[label distance=-3pt]below:$Z$}] at (1.5,4) {};
    \node[label={[label distance=-3pt]below:$z$}] at (0.5,4) {};
    \node[label={[label distance=-3pt]above:$W$}] at (1.5,3) {};
    \node[label={[label distance=-3pt]above:$w$}] at (0.5,3) {};

    \node[label={[label distance=-3pt]below:$w$}] at (3.5,1) {};
    \node[label={[label distance=-3pt]below:$W$}] at (2.5,1) {};
    \node[label={[label distance=-3pt]above:$z$}] at (3.5,0) {};
    \node[label={[label distance=-3pt]above:$Z$}] at (2.5,0) {};

    \node[label={[label distance=-3pt]below:$Z$}] at (3.5,2) {};
    \node[label={[label distance=-3pt]below:$z$}] at (2.5,2) {};
    \node[label={[label distance=-3pt]above:$W$}] at (3.5,1) {};
    \node[label={[label distance=-3pt]above:$w$}] at (2.5,1) {};

    \node[label={[label distance=-3pt]below:$w$}] at (3.5,3) {};
    \node[label={[label distance=-3pt]below:$W$}] at (2.5,3) {};
    \node[label={[color=red,label distance=-3pt]above:$z$}] at (3.5,2) {};
    \node[label={[label distance=-3pt]above:$Z$}] at (2.5,2) {};

    \node[label={[label distance=-3pt]below:$Z$}] at (3.5,4) {};
    \node[label={[label distance=-3pt]below:$z$}] at (2.5,4) {};
    \node[label={[label distance=-3pt]above:$W$}] at (3.5,3) {};
    \node[label={[label distance=-3pt]above:$w$}] at (2.5,3) {};

    \draw[dotted] (0,0)--(-0.2,0);
    \draw[dotted] (0,2)--(-0.2,2);
    \draw[dotted] (0,4)--(-0.2,4);
    \draw[dotted] (4,0)--(4.2,0);
    \draw[dotted] (4,2)--(4.2,2);
    \draw[dotted] (4,4)--(4.2,4);

    \draw[dotted] (0,0)--(0,-0.2);
    \draw[dotted] (1,0)--(1,-0.2);
    \draw[dotted] (2,0)--(2,-0.2);
    \draw[dotted] (3,0)--(3,-0.2);
    \draw[dotted] (4,0)--(4,-0.2);
    \draw[dotted] (0,4)--(0,4.2);
    \draw[dotted] (1,4)--(1,4.2);
    \draw[dotted] (2,4)--(2,4.2);
    \draw[dotted] (3,4)--(3,4.2);
    \draw[dotted] (4,4)--(4,4.2);

    \draw[dotted] (-0.2,1)--(4.2,1);
    \draw[dotted] (-0.2,3)--(4.2,3);

    \draw[red,very thick] (0.1,0.1)--(4.1,2.1);
    \draw[red, thick] (1,0)--(1,1);
    \draw[red, thick] (2,1)--(2,2);
    \draw[red, thick] (3,1)--(3,2);
    \draw[red, thick] (4,2)--(4,3);
    \draw[red, thick] (1,0)--(1,1);
    \draw[red, thick] (1,1)--(2,1);
    \draw[red, thick] (3,2)--(4,2);
  \end{tikzpicture}%
  \caption{The cutting sequence for slope $ 1/2 $: $ y W x Y z X $. (To make the image cleaner, we have written $ x = X^{-1} $, $ y = Y^{-1} $, and so on.)\label{fig:sphere_cylinder}}
\end{figure}
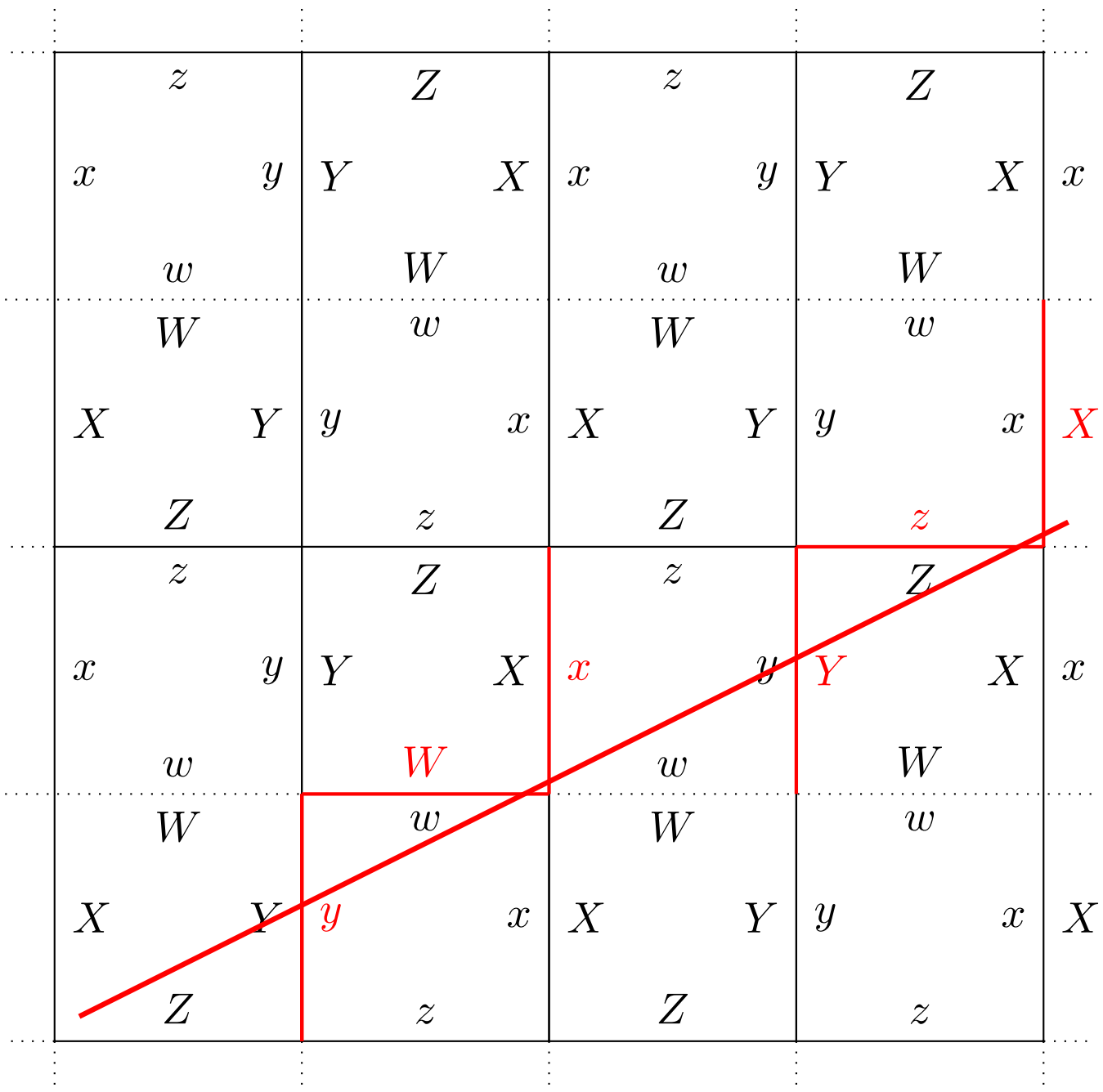

\begin{lem}\label{lem:word_image}
  The map $ p/q \mapsto \Farey_{p/q} $ is injective and has image equal to the set of all finite words in $ \mathfrak{F}^{\pm 1} $
  such that (i) no letter is ever followed by its inverse, and (ii) the signs of exponents in the word alternate.
\end{lem}
\begin{proof}
  These observations are clear from the consideration of the definition of cutting sequences in terms of a tiling of the plane, see
  for instance \zcref{fig:sphere_cylinder}.
\end{proof}

\begin{lem}\label{lem:obvious_symmetry}
  The word $ \Farey_{p/q} $ is of length $ 2p+2q $ and can be decomposed as a product $ R QS R^{-1} TU $ where $ Q,S,T,U \in \mathfrak{F}^{\pm 1} $,
  $ R $ is the first $ p+q-1 $ elements of the word, and $ QS $ and $ TU $ are conjugate to elements of $ \mathfrak{P}^{\pm 1} $.
\end{lem}
\begin{rem}
  This result is a generalisation of the computation of the relator of a two-bridge link group, see Lemma~2.4 of \cite{ems22b}.
\end{rem}
\begin{proof}
  The decomposition as a product of an $ \mathfrak{F}$ is obvious after considering the action of the rotation of angle $\pi $ with centre $ (p, q) $ on $\Farey_{p/q}$.
  The final observation, that $ QS $ and $ TU $ represent loops around the holes, follows from \zcref{lem:word_image} and the computation of the $ \mathfrak{F}$-cutting sequences of $A,B,C,D $.
\end{proof}

We will need the following technical lemma.
\begin{lem}[{\cite[\S 3.3]{hardywright}}]\label{lem:farey_closeness}
  If $ p/q $ and $ r/s $ are Farey neighbours with $ p/q < r/s $, then $ p/q < (p+r)/(q+s) < r/s $, and $ (p+r)/(q+s) $ is the unique fraction of minimal
  denominator between $ p/q $ and $ r/s $. More precisely, let $ u/v $ be any fraction in the open interval $ (p/q,r/s) $; then there exist two positive integers $ \lambda, \mu $
  such that
  \begin{displaymath}
    u = \lambda p + \mu r \text{ and } v = \lambda r + \mu s
  \end{displaymath}
  (so the minimal denominator is obtained when $ \lambda = \mu = 1 $). \qed
\end{lem}

\begin{figure}
  \centering
      \begin{tikzpicture}[scale=.6]
      \draw[step=1.0,black] (0,0) grid (16,6);
      \draw[thick,red] (0,0)--(6,2);
      \draw[thick,green] (6,2)--(16,6);
      \draw[thick,blue] (0,0)--(16,6);
      \draw[very thick] (3,0)--(3,6);
      \draw[very thick] (6,0)--(6,6);
      \draw[very thick] (8,0)--(8,6);
      \draw[very thick] (11,0)--(11,6);
      \draw[very thick] (16,0)--(16,6);
      \draw[very thick] (0,1)--(16,1);
      \draw[very thick] (0,2)--(16,2);
      \draw[very thick] (0,3)--(16,3);
      \draw[very thick] (0,4)--(16,4);
      \draw[very thick] (0,6)--(16,6);
      \node[label={below:$q$}]   at (3,0) {};
      \node[label={below:$2q$}]  at (6,0) {};
      \node[label={below:$q+s$}] at (8,0) {};
      \node[label={below:$2q + s$}]   at (11,0) {};
      \node[label={below:$2q + 2s$}]  at (16,0) {};
      \node[label={left:$p$}]   at (0,1) {};
      \node[label={left:$2p$}]  at (0,2) {};
      \node[label={left:$p+r$}] at (0,3) {};
      \node[label={left:$2p + r$}]   at (0,4) {};
      \node[label={left:$2p + 2r$}]  at (0,6) {};
    \end{tikzpicture}%
  \caption{Farey addition versus word multiplication for $ \Farey_{p/q} $ (red) and $ \Farey_{r/s} $ (green).\label{fig:concatenation}}
\end{figure}
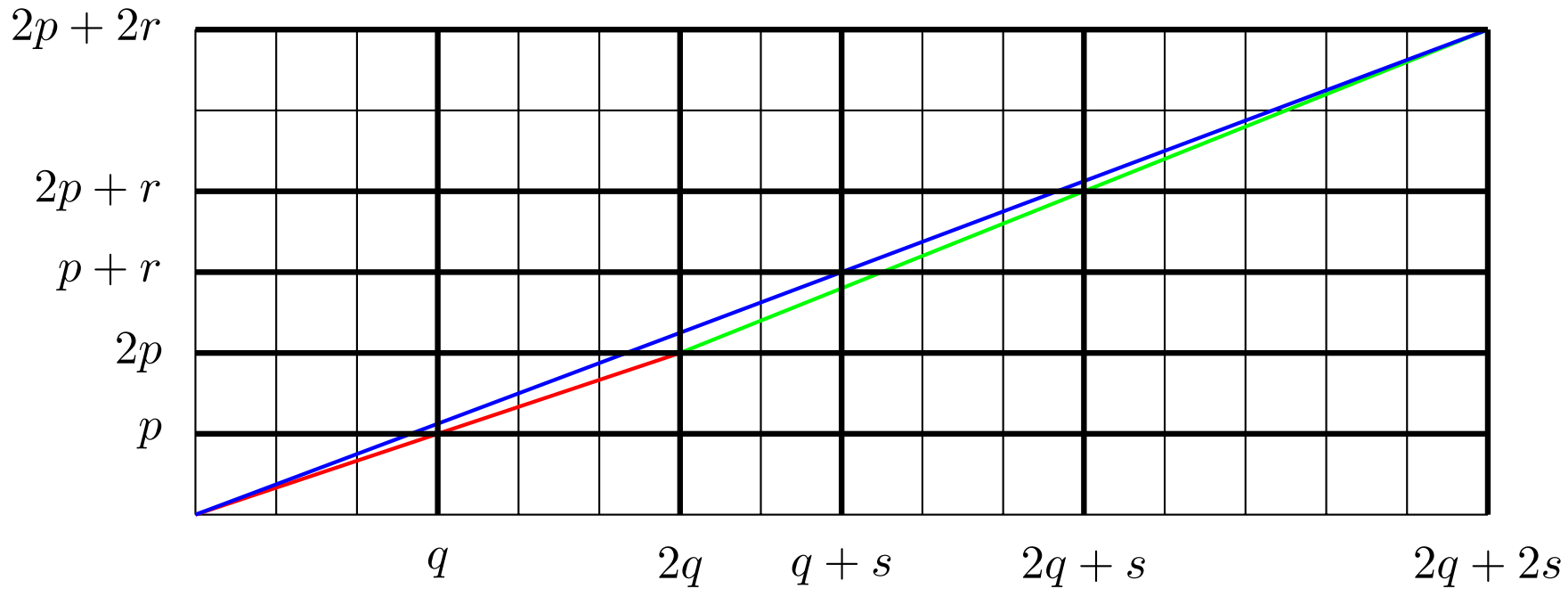

If $ p/q $ and $ r/s $ are neighbouring points in the Farey graph $ \mc{F} $ such that $ \gcd(p,q) = \gcd(r,s) = 1 $, then we define
their \df{Farey sum} to be the mediant $ (p/q)\oplus(r/s) \coloneq (p+r)/(q+s) $~\cite[212]{indras_pearls}.

\begin{rem}\label{rem:more_general_recursion}
  The following recursion is more general than is needed, and in fact is sufficiently general to work for giving a recursion on words in a four-holed sphere
  subset of any surface. That is, it works to give a recursion in $ \KK(S) $ along a fan of $n$-simplices which all share $n-1$ vertices. In our setting
  the curve $A$ (resp.\ $B$) will represented by the inverse of the element in $ \pi_1(\Sigma) $ representing $D$ (resp.\ $C$). The more general result can be
  used to give word recursions in the curve complex of arbitrary surfaces: a top simplex of the curve complex represents a pants decomposition, and one can
  apply the lemma to pairs of pants which share a common bounding curve.
\end{rem}

In the following lemma, the \df{pivots} of the concatenation $ \Farey_{p/q} \Farey_{r/s} $
are the $ (q+s)$th letter in the subsequence of letters from $ \{ X,Y \} $ and the $ (p+r)$th letter in the subsequence of letters from $ \{W,Z\} $.

\begin{lem}\label{lem:product_word}
  Let $ p/q $ and $ r/s $ are Farey neighbours with $ p/q < r/s $. Then the
  word $ \Farey_{p/q \oplus r/s} $ is obtained from $ \Farey_{p/q} \Farey_{r/s} $ by exchanging the two pivots and replacing them with their inverses.
\end{lem}
\begin{proof}
  The result is a generalisation of Lemma~4.1 of \cite{ems22b}, and essentially follows from a detailed consideration of the lines of slopes $ p/q $, $ r/s $,
  and $ (p+r)/(q+s) $ on the lattice, \zcref{fig:concatenation}. We consider separately the subword obtained by considering only `horizontal' cuts (i.e.\ the subword
  consisting only of letters from $ \{W,Z\} $) and the subword obtained by considering only `vertical' cuts.

  Consider first the substitution $ W = Z = 1 $; then the theorem statement is exactly
  \begin{equation}\label{eq:interleaved_cuts}
    \text{$\Farey_{p/q \oplus r/s} $ is  obtained by flipping the $ q+s$th exponent of $ \Farey_{p/q} \Farey_{r/s} $}
  \end{equation}
  (this is exactly Lemma~4.1 of \cite{ems22b}). The equivalent statement for the substution $ X = Y = 1 $ is obtained by applying \eqref{eq:interleaved_cuts} to the
  transpose of the cutting tiling (i.e.\ the reflection across the diagonal of \zcref{fig:sphere_cylinder}). The full statement of the result is now obtained from the additional
  observation that the order of the two interleaved cutting sequences is flipped at the point $ (p+r,q+s) $ when comparing the line of slope $ (p+r)/(q+s) $ to the
  two lines of slopes $ p/q $ and $ r/s $. It remains to recall the proof of \eqref{eq:interleaved_cuts} from \cite{ems22b}.

  Write $ h(i) $ for the height $ (p/q \oplus r/s) i $. Observe that $ h(i) $ is integral only at $ i = 0 $ and $ i = 2q+2s $ (at both
  positions trivially the letters in $ \Farey_{p/q \oplus r/s} $ and $ \Farey_{p/q} \Farey_{r/s} $ are identical) and at $ i = q+s $. The lemma will follow once we check that that at the positions $ i \not\in \{0,2q+2s\} $,
  \begin{displaymath}
    \text{if } 0 < i \leq 2q  \text{ then } (p/q) i < h(i) < \ceil (p/q) i
  \end{displaymath}
  and
  \begin{displaymath}
    \text{if } 0 < i < 2s \text{ then } (r/s)i + 2q < h(i + 2p) < \ceil [(r/s)i + 2p]:
  \end{displaymath}
  indeed, these inequalities show that at every integral horizontal distance the height of the line corresponding to $ \Farey_{p/q \oplus r/s} $ is meeting the same vertical line segment
  as the line corresponding to $ \Farey_{p/q} $ or $ \Farey_{r/s} $, and so the letter chosen is the same except at $ i = q+s $ since at this position the height of the line of
  slope $ (p/q)\oplus(r/s) $, being integral, is rounded up to $ h(i)+1 $ while the height of the line of slope $ r/s $ is non-integral so is rounded up to the integer $ h(i) $.

  Observe now that the two displayed inequalities are equivalent to the following: there is no integer between $ (p/q)i $ and $ (p/q \oplus r/s)i $ (exclusive)
  if $ 0 < i \leq 2q $, and there is no integer between $ (r/s)i + 2q $ and $ h(i + 2p) $ if $ 0 < i < 2s $. But these follow from \zcref{lem:farey_closeness}. Indeed,
  the lemma shows that no integer lies between $ p/q $ and $ (p/q) \oplus (r/s) $; suppose $ a/b $ is a rational between $ (p/q)i $ and $ h(i) $,
  then $ a = i(\lambda p + \mu (p+r)) $ and $ b = (\lambda q + \mu (q+s)) $ for some positive $ \lambda, \mu $; suppose $ a/b \in \IZ $, so $ \lambda q + \mu (q+s) $ divides
  $ i(\lambda p + \mu (p+r)) $. By the case $ i = 1 $, $ \lambda p + \mu (p+r) $ and $ \lambda q + \mu (q+s) $ are coprime, so $ \lambda q + \mu (q+s) $ divides
  $ i $; but $ \lambda q + \mu (q+s) \geq 2q $. The case of $ (r/s)i + 2q $ and $ h(i + 2p) $ is proved in a similar way.
\end{proof}

We would also like to consider the action on dual simplices under Farey moves; i.e.\ given a pair of simplices $ \Delta_1, \Delta_2 $ which are incident
in the Farey graph, and a dual simplex $ \Delta_1^\perp $ of $ \Delta_1 $, we wish to find a dual simplex $ \Delta_2^\perp $ of $ \Delta_2 $ which
is incident to $ \Delta_1^\perp $ in the Farey graph.

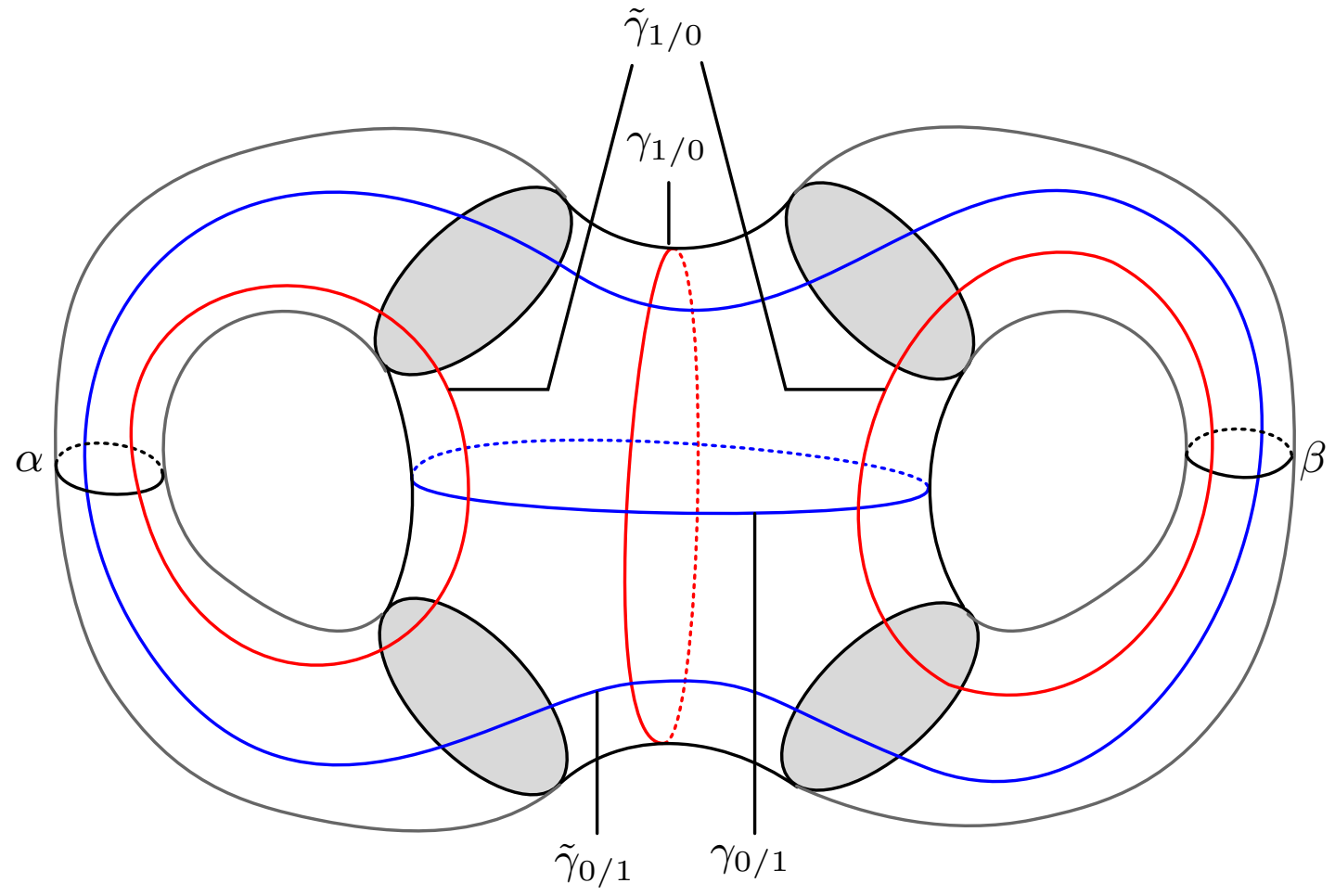
\begin{figure}\vspace{1.5ex}%
  \labellist
  \small\hair 2pt
  \pinlabel {$\gamma_{1/0}$} [b] at 256 275
  \pinlabel {$\gamma_{0/1}$} [t] at 291 0
  \pinlabel {$\tilde{\gamma}_{0/1}$} [t] at 226 0
  \pinlabel {$\tilde{\gamma}_{1/0}$} [b] at 255 323
  \pinlabel {$\alpha$} [r] at 0 156
  \pinlabel {$\beta$} [l] at 516 156
  \endlabellist
  \centering
  \includegraphics[width=0.7\textwidth]{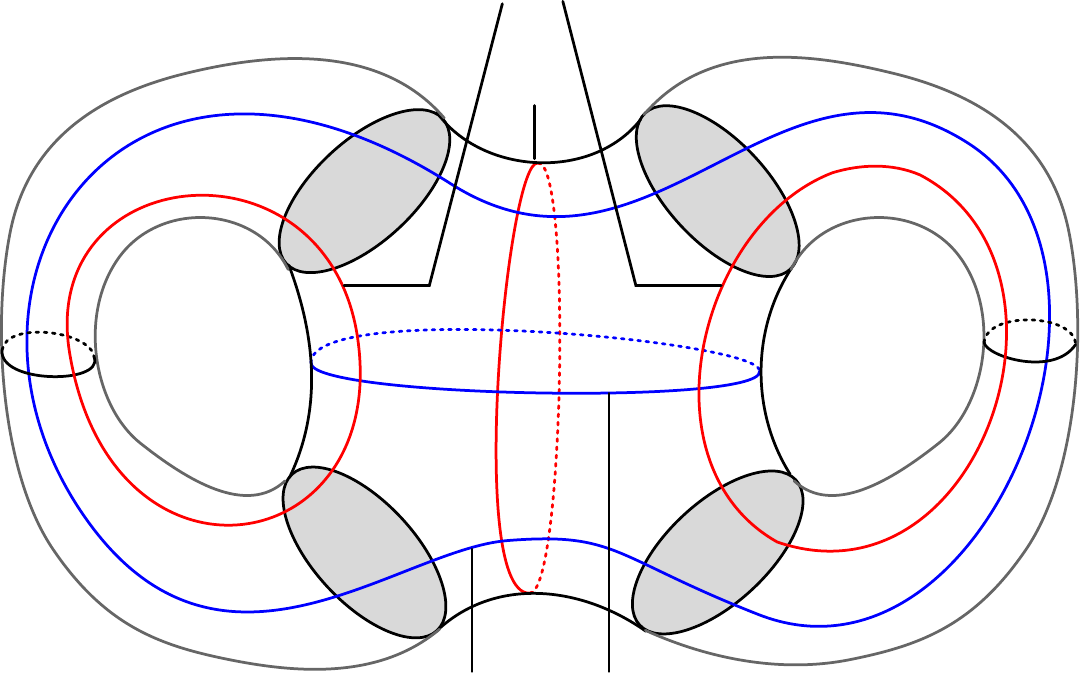}\vspace{1ex}%
  \caption{Passing from an ideal triangulation of the four-holed sphere to a dual, and conversely.\label{fig:ideal_triangulation_gives_dual}}
\end{figure}

The following lemma is analogous to the discussion in \S 1.2 of Akiyoshi, Sakuma, Wada, and Yamashita~\cite{akiyoshi}. The reason we need this result
is that it is a method for (semi-canonically) producing a homology basis for a surface group generated by the representatives in $ \pi_1(S_{2,0}) $
of the vertices of a $2$-simplex of $ \KK(S_{2,0}) $. In the special case that the vertices of the simplex give meridians of the link components,
the dual $2$-simplex will give a canonical choice of longitudes for those components. Once we have a homology basis we can use the curve coding theory
to give word recursions in the surface subgroup, since (by induction) we know explicitly how to expand the alphabet given by the homology basis in terms of the
generators of the higher level groups.
\begin{lem}\label{lem:product_word_dual}
  For each $ p/q \in \IQ $, let $ \hat{\gamma}_{p/q} $ be the lift of $ \gamma_{p/q} $ to $ \IR^2 $ as in \zcref{fig:concatenation}. Let $ \hat{\tilde{\gamma}}_{p/q} $
  be the union of all vertical shifts of $\hat{\gamma}_{p/q}$ which pass through a lattice point of $ \IZ^2 \subset \IR^2 $, and let $ \tilde{\gamma}_{p/q} $
  be the projection of $\hat{\tilde{\gamma}}_{p/q}$ back to the four-punctured sphere. If $ S_{0,4} $ is embedded into a compact genus two surface in the natural
  way by deletion of two non-separating annuli from $ S_{2,0} $ with core loops $ \alpha $ and $ \beta $, then $\tilde{\gamma}_{p/q}$ can be completed to a multicurve
  on $ S_{2,0} $ with one or two connected components that forms a $0$-{} or $1$-simplex of $ \KK(S_{2,0}) $. This construction has the following properties:
  \begin{enumerate}
    \item If $ \gamma_{p/q} $ is any simple closed curve on $ S_{0,4} $, then a dual or pseudodual (depending on whether $ \gamma_{p/q} $ is separating or not)
          $2$-simplex to $ \{ \alpha, \gamma_{p/q}, \beta \} $ is given by $ \tilde{\gamma}_{p/q} \union \tilde{\gamma}_{r/s} $ where $ r/s $ is any Farey neighbour of $p/q$.
    \item If $ r/s $ is a Farey neighbour of $p/q$, then the dual or pseudodual simplices constructed in (1) to $\{ \alpha, \gamma_{p/q}, \beta \}$ and $ \{ \alpha, \gamma_{r/s}, \beta \} $
          are incident simplices in $ \KK(S_{2,0})$.
    \item If the embedding of $ S_{0,4} $ into $ S_{2,0} $ is done such that the handles are glued between $B$ and $ C $, and $A$ and $D$ respectively (so that the arcs $X$ and $Y $
          are completed to the curves $ \tilde{\gamma}_{1/0} $), then the words for the components of $ \tilde{\gamma}_{p/q} \union \tilde{\gamma}_{r/s} $ are obtained from the words
          for $ \gamma_{p/q} $ and $ \gamma_{r/s} $ by a combination of the operations of (i) concatenation and (ii) appending the letters $F$ and $G$ (or their inverses) which respectively
          correspond to crossing the loops $ \alpha $ or $ \beta $.
  \end{enumerate}
\end{lem}
\begin{proof}
  By a standard argument involving transitivity properties of the mapping class group, it suffices to prove (1) for the configuration of curves given in \zcref{fig:ideal_triangulation_gives_dual}
  for which it is obvious. Part (2) may be proved by observing that passing from $\Delta_1 = \{ \alpha, \gamma_{p/q}, \beta \}$ to $\Delta_2 = \{ \alpha, \gamma_{r/s}, \beta \}$ is done by a single
  Dehn half-twist, and this half-twist acts also on the dual or pseudodual simplex to $ \Delta_1$ produced by (1) to give a dual or pseudodual simplex to $\Delta_2$. This leaves part (3), which is proved
  by appealing to the same argument as \zcref{lem:product_word}; the point is that if $ \tilde{\gamma}_{p/q} $ has two topological components then both have the same cutting sequence as $ \gamma_{p/q} $
  except for an addition of $\alpha$ to one and $ \beta $ to the other; if $ \tilde{\gamma}_{p/q}$ has a single component then it corresponds to two copies of $\gamma_{p/q}$ with an $ F^{\pm} $ between them
  and a $ G^{\pm 1} $ appended to the result.
\end{proof}

As a consequence, we immediately obtain the following group filtration:
\begin{thm}\label{thm:filtration_in_component}
  Let $ \Delta_0, \ldots, \Delta_k $ be a ns-geodesic of $2$-simplices in the curve complex of the genus $2$ surface $\Sigma = \partial \mc{H} $, such that
  \begin{enumerate}
    \item $ \Delta_i \inter \Delta_{i+1} $ is a $1$-simplex for each $ 0 \leq i < k $, and
    \item If $ \abs{i-j}>1 $ then $ \dim (\Delta_i \inter \Delta_j) < 1 $.
  \end{enumerate}
  For each $i$ let $\Delta_i^\perp$ be the dual simplex (or possibly, at the final step, the pseudodual simplex) given by part (1) of \zcref{lem:product_word_dual}.

  Let $ \rho : \pi_1(\mc{H}) \to \PSL(2,\IC) $ be a representation. Then the sequence of subgroups
  \begin{equation}\label{eq:filtration_in_component_index}
    \Gamma_{k-i} = \langle \Delta_{k-i}, \Delta_{k-i}^{\perp} \rangle \coloneq \langle \rho(\gamma) \in \Delta_{k-i} \cup \Delta_{k-i}^\perp \rangle
  \end{equation}
  is a filtration of $ \Gamma_k $. That is,
  \begin{displaymath}
    \Gamma_0 < \cdots < \Gamma_k
  \end{displaymath}
  and each $ \Gamma_i $ admits a generating set $ \langle U_i, V_i, F_i, G_i \rangle $ such that $ U_{i-1} = \Farey_{p/q}(U_i, V_i) $ (for some $ p/q \in \IQ $) and $ V_{i-1} \in \{ V_i, U_i \} $,
  and where $ F_{i-1} $ and $ G_{i-1} $ are written in terms of $\Farey_{p/q}(U_i, V_i)$ via part (3) of \zcref{lem:product_word_dual} where the letters $ F $ and $ G $ in that statement
  are substituted with $ F_i $ and $ G_i $. In particular the filtration is made up of groups generated by the Farey rewriting systems of \zcref{lem:product_word} and \zcref{lem:product_word_dual}. \qed
\end{thm}
Note the indexing `swap' in \eqref{eq:filtration_in_component_index}, which means that the `simplest' simplex $ \Delta_0 $ is associated with the largest group $ \Gamma_k $,
and (at the other end) the `complicated' simplex $ \Delta_k $ is associated with $ \Gamma_0 $. As we move along the ns-geodesic, we `move' complexity from the simplex
into the group or vice versa.

\begin{rem}
  For computational reasons, one would really like to have a recursion on the traces of the generators for the groups in the filtration of \zcref{thm:filtration_in_component} as well
  as on the words. In earlier work~\cite{ems22b} we were able to deduce such a recursion from a product formula like \zcref{lem:product_word} in the special case of two-bridge link
  groups. The obstruction to these methods working here is our lack of understanding of the dual curves: these curves depend strongly on global data about the group $G$ (as they
  intersect the curves defining the filtration many times and so pass through many thrice-holed sphere subgroups of the surface groups) and so it is very hard to control their traces.
\end{rem}

\section{Tunnel number one links in the sphere}\label{sec:tno_in_sphere}
In this section we deal with the structure theory of holonomy groups of tunnel number $1$ links in the $3$-sphere. In \zcref{sec:two_bridge} we consider
the two-bridge case, where \zcref{thm:filtration_in_component} is the trivial filtration. The general case is considered in \zcref{sec:general_in_sphere}

\subsection{Two-bridge links}\label{sec:two_bridge}
We begin with the special case of a Heegaard diagam $ \Delta \in \KK^{(2)}(S_{2,0}) $ of the form $ \{X,Y,\gamma\} $; recall that $ X $ and $ Y $ are the handle cores of the exterior
handlebody $ \mc{H}^+ $. Observe that if a $2$-handle is attached to either $X$ or $Y$ then the resulting manifold is non-hyperbolisable (the torus boundary component is compressible),
and so we must attach our $2$-handle along $\gamma$ to obtain a non-elementary Kleinian group. We also require $\gamma$ to be sufficiently complicated.

We see first that there is a path from the trivial simplex $ \Delta_0 $ (\zcref{defn:deltanought}) to the simplex $\Delta$ in the Farey graph $ \mc{F} $, corresponding
to fixing the two curves $ X $ and $Y$ and rational braiding the third curve via Farey moves. The manifold $ \Sigma \setminus (X \union Y) $ is a four-holed sphere, and so has homology
$ \IZ^2 $. It is well-known that simple closed curves on the four-holed sphere are in bijection with homology classes $ (p,q) \in \IZ^2 $ with $ p, q $ coprime (e.g.\ by
combining Propositions 1.5 and 2.6 of Farb and Margalit~\cite{farb}); we choose coordinates such that $ (0,1) $ is the unique simple closed curve which bounds a compression
disc in $ \mc{H}^+$. Since $\gamma$ is a simple closed curve on $\Sigma \setminus (X \union Y)$, it determines such a pair $ (p,q) $, and the continued fraction
decomposition of $ q/p $ determines the path in $ \mc{F} $ as described in \cite[\S 4.5]{graham}: if we restrict to vertices in $ \IQ \inter (0,1) $ then the
edges of the graph at every vertex can be given a labelling $L$ or $R$, so that if one starts at $ 0/1 $ then the geodesic in $ \mc{F} $ to $ p/q = [a_1, a_2, \ldots, a_n] $ is
\begin{equation}\label{eq:conway_word}
  \underbrace{L \cdots L}_{a_1\; \text{times}}\, \underbrace{R \cdots R}_{a_2\; \text{times}}\, \underbrace{L \cdots L}_{a_3\; \text{times}} \cdots.
\end{equation}

To see the knot-theoretic meaning of this fraction, we recall:
\begin{defn}
  A \df{two-bridge link} is a link which admits a diagram in $ \IR^2 $ with exactly two local maxima in the $y$-ordinate~\cite[Chapter~12]{burde3e}.
\end{defn}
The two-bridge links were classfied by Schubert~\cite{schubert56} in terms of two coprime parameters. If these parameters are $ \alpha $ and $ \beta $ then we denote
the corresponding link by $ \mf{b}(\alpha,\beta) $; our notation here follows the modern exposition of Schubert's work in Proposition~12.3 and Theorem~12.6 of \cite{burde3e}.
The link between these parameters and the primitive homology classes of $ S_{0,4} $ comes from the following observation:

\begin{figure}
  \labellist
  \small\hair 2pt
  \pinlabel {$L$} [r] at 138 115
  \pinlabel {$R$} [b] at 88 55
  \endlabellist
  \centering
  \includegraphics[width=.33\textwidth]{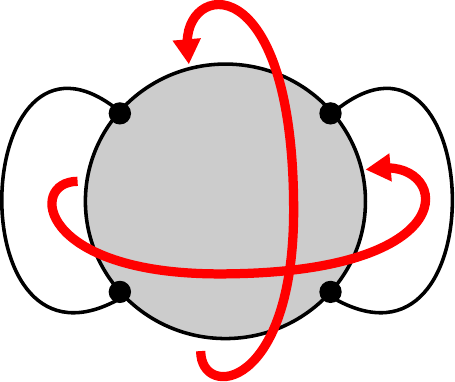}
  \caption{The surface $ S_{0,4} $ has fundamental group virtually generated by the Dehn half-twists $L$ and $R$ indicated. Both of these mapping classes extend by isotopy
  into the shaded handlebody $ \mc{H}^-$.\label{fig:2bridge}}
\end{figure}

\begin{prp}
  Suppose that $ \alpha, \beta \in \IZ $ are coprime, and let $ [a_1,a_2, \ldots, a_n] $ be the continued fraction representation of $ \beta/\alpha $.
  Then the two-bridge link $ \mf{b}(\alpha,\beta) $ is obtained by taking a trivial $2$-tangle (bounded by a four-punctured sphere) and applying
  the word \eqref{eq:conway_word} where $L$ and $R$ are interpreted as the mapping classes of $ S_{0,4} $ in \zcref{fig:2bridge} extended into the
  handlebody $ \mc{H}^- $ by isotopy. Equivalently, $\mf{b}(\alpha,\beta)$ is obtained from $ \mc{H}^+ $ by gluing a $2$-handle along the simple closed curve
  on $ \Sigma \setminus \{ X,Y\} $ with homology class $ (\alpha,\beta) $. \qed
\end{prp}

\begin{figure}
  \begin{subfigure}{.24\textwidth}
    \centering
    \labellist
    \small\hair 2pt
    \pinlabel {$\times 2$} [r] at 0 46
    \endlabellist
    \centering
    \includegraphics[width=30mm]{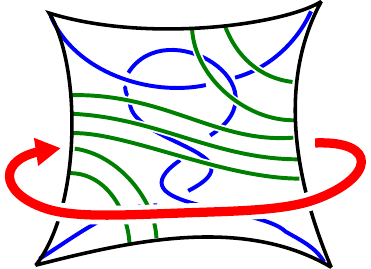}
    \caption{}
  \end{subfigure}\hfill%
  \begin{subfigure}{.24\textwidth}
    \centering
    \includegraphics[width=24.140mm]{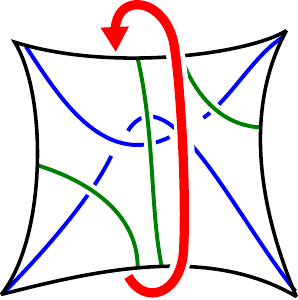}
    \caption{}
  \end{subfigure}\hfill%
  \begin{subfigure}{.24\textwidth}
    \centering
    \includegraphics[width=26.493mm]{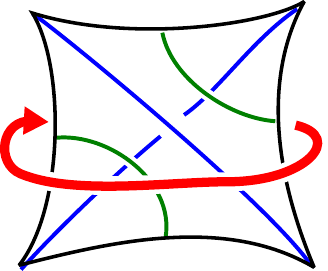}
    \caption{}
  \end{subfigure}\hfill%
  \begin{subfigure}{.24\textwidth}
    \centering
    \includegraphics[width=24.140mm]{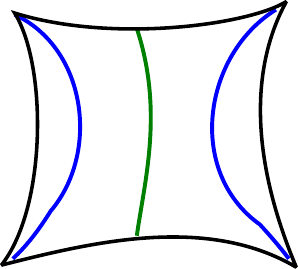}
    \caption{}
  \end{subfigure}
  \caption{Untwisting a figure $8$ knot together with a curve bounding a compression disc.\label{fig:fig8braid}}
\end{figure}

\begin{ex}
  The figure $8$ link has slope $ 2/5 $. This has continued fraction decomposition
  \begin{displaymath}
    2/5 = \cfrac{1}{1+\cfrac{1}{1+\cfrac{1}{2}}} = [1,1,2]
  \end{displaymath}
  and thus the link is produced by alternately performing a sequence of two half twists and one full twist as in \zcref{fig:fig8braid}: here
  we show the rational tangle in blue, and the curve bounding a compression disc in green. This curve is represented by the Farey word $ \Farey_{2/5}(X,Y) $
  in the free group~\cite{ems22b}.
\end{ex}

Since the holonomy group of the gluing has nontrivial kernel, it is a non-free discrete subgroup of $ \PSL(2,\IC) $ generated by two elements. Such groups
have been completely classified by Akiyoshi, Ohshika, Parker, Sakuma, and Yoshida~\cite{akiyoshi2020classification}. They proved the following theorem,
which clarifies computer experiments conducted by Riley in 1979~\cite[\S 3 of preface]{akiyoshi}.

\begin{prp}\label{prp:akiyoshi}
  Suppose that $ \rho : \pi_1(\mc{H}^+) \to \PSL(2,\IC) $ is discrete with finite covolume and solves the system of equations
 \begin{displaymath}
    \tr^2 \rho(X) = \tr^2 \rho(Y) = \tr^2 \rho(\gamma) = 4.
  \end{displaymath}
  Then there exists a hyperbolic two-bridge link $\mf{b}$ and a subgroup $ K \leq \Aut(\mf{b}) \simeq \IZ/2\IZ \oplus \IZ/2\IZ $ such that $ \Gamma = \rho(\pi_1(\mc{H}^+)) $
  is the orbifold holonomy group of the complete hyperbolic metric on $(\IS^3 \setminus \mf{b})/K$. In particular there is a simple closed curve $\gamma$ on the four-holed
  sphere $ \Sigma \setminus (X\cup Y) $ such that $ \IH^3/\Gamma$ is obtained by gluing a $2$-handle along $ \gamma $ and then performing a finite quotient. \qed
\end{prp}
Thus the two-bridge link groups are the only torsion-free groups arising in the solution to \zcref{prob:general} when two of the three curves are handle cores.

A famous construction of Sakuma and Weeks~\cite{sakuma95} gives an exhaustion of $ \IS^3 \setminus \mf{b} $ by embedded $ S_{0,4} \times I $ submanifolds bounded by piecewise geodesic $S_{0,4}$'s.
The homotopy groups of these subsurfaces are computed in Chesebro~\cite[\S 5.3]{chesebro25}. Further, it is known that this triangulation is \emph{canonical} in the sense of Epstein
and Penner~\cite{epstein88}, so it is geometric and dual to the Ford domain of the corresponding group~\cite[appendix by Futer]{gueritaud06}.

\subsection{General tunnel number one links in the sphere}\label{sec:general_in_sphere}
We now recall from \cite{cho09} the relation between the theory of the previous sections and the theory of tunnel number $1$ links. Suppose that $ \mf{k} \subset \IS^3 $
is tunnel number $1$. Then it corresponds to a genus two Heegaard diagram as explained in \zcref{lem:tunnel_number_1_is_heegaard} together with a choice of one of the three
curves to bound the disc dual to the unknotting tunnel, and we can construct a corresponding ns-geodesic. Thus, from \zcref{thm:filtration_in_component} we obtain:
\begin{cor}\label{cor:in_sphere}
  There exists a filtration of the holonomy group of a tunnel number $1$ link complement $ \IS^3 \setminus \mf{k} $ determined by a choice of path in the graph of $2$-simplices
  of $ \KK(\Sigma) $ between the canonical basepoint $\Delta_0$ (\zcref{defn:deltanought}) and the Heegaard diagram $\Delta_n$ associated with $ \mf{k} \union \tau $, where $ \tau $ is a choice
  of tunnel. In particular we can write
  \begin{displaymath}
    \Gamma_0 < \cdots < \Gamma_n = \pi_1(\IS^3 \setminus \mf{k}) < \PSL(2,\IC)
  \end{displaymath}
  where each $ \Gamma_i $ is a $ S_{2,0}$-surface representation. Each $ \Gamma_i $ is the homotopy group of an embedded
  topological genus $2$ surface, and $ \Gamma_0 $ is $ V_4 \coloneq \IZ \oplus \IZ $ (resp.\ $ V_4 * V_4$) if $ \mf{k} $ has one
  component (resp.\ two components). \qed
\end{cor}

\begin{figure}
  \centering
  \includegraphics[width=.8\textwidth]{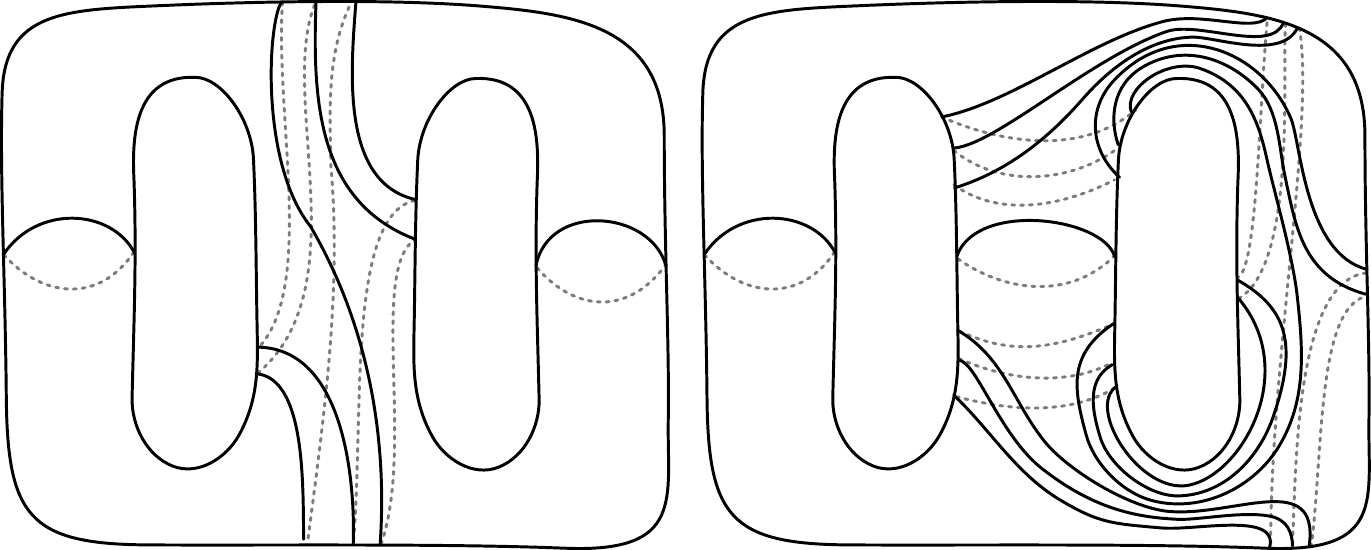}
  \caption{Swapping a handle from the right to the centre increases the complexity of the curves on the surface dramatically.\label{fig:cabling_swap}}
\end{figure}

In terms of the knot theory, the filtration gives a procedure for building the knot $ \mf{k} $ from the unknot via an iterative process. In each step, the knot strands
are separated into three `bundles', one within each `column' of a genus $2$ surface; then rational braiding moves (\zcref{cons:cabling}) are performed,\ i.e.\ the strands in the middle
bundle are cut apart and reglued according to the same procedure as for two-bridge links. This procedure is then repeated by performing cabling moves on a different `column' of the
surface, and can be iterated as many times as desired. When one tries to actually do it in practice, it becomes clear why working directly with word recursion and the filtration is preferrable: the problem is
that swapping the handles by isotopy increases the complexity of the words dramatically (we show a very simple example in \zcref{fig:cabling_swap}), and constantly changing
generators to longer words when swapping the handles around hides this complexity.

\subsection{Presenting tunnel number one links as Dehn fillings}
A fully augmented link is a link $ \mf{f} $ in $ \IS^3 $ whose set of components can be written as the disjoint union of a set $ \{\mf{p}_1, \ldots, \mf{p}_k\} $ of \df{planar components}
and a set $ \{ \mf{a}_1, \ldots, \mf{a}_l \} $ of \df{augmentation components} such that, up to isotopy in $ \IS^3 $:
\begin{itemize}
  \item There exists a $2$-sphere $ P \subset \IS^3  $ such that $ \bigcup_{i=1}^k \mf{p}_i \subset P $.
  \item Each cusp $ \mf{a}_j $ bounds an embedded disc $ D_j $ in $ \IS^3 $ so that
    \begin{itemize}
      \item $ D_j $ intersects $P $ transversely in a line segment $L_j$ and
      \item $ D_j $ intersects $ \mf{f} $ exactly twice, both intersections being transverse and both intersections being with planar components,
            and with the intersection being essential. Here, we say an intersection is inessential if both intersections are with the same component $ \mf{p}_i $ and
            there is an arc $ \alpha \subset \mf{p}_i $ with endpoints $ \mf{p}_i \inter D_j $ such that $ L_j \union \alpha $ bounds an embedded disc in $ \IS^3 \setminus \mf{f} $.
    \end{itemize}
\end{itemize}
Now suppose that $ \mf{k} \subset \IS^3 $ is a hyperbolic link which admits a diagram consisting of even twist regions: that is, there is a diagram $D$ of $ \mf{k} $ consisting of
twist regions with only even crossing numbers, connected together by a family of disjoint arcs embedded in the plane of the diagram. As explained in Purcell~\cite[Chapter~7]{purcell},
one can find a fully augmented link $ \mf{f}(D) $ such that $ \IS^3 \setminus \mf{k} $ is obtained as a Dehn filling
of $ \IS^3 \setminus \mf{f}(D) $---namely, replace each twist region of $D$ with two parallel strands orbited by an augmentation component. Since hyperbolic fully augmented links
have very strongly controlled geometry (e.g.\ they admit a decomposition by regular ideal octahedra) and hyperbolic Dehn filling has well-controlled geometry when filling along
angles $ 1/n $ for $ n \gg 1 $, one can deduce geometric information about $ \mf{k} $ from the corresponding fully augmented link. With some modifications, namely allowing the planar
components to leave the plane and cross each other inside augmentation cusps, one can also deal with odd twist regions (which can be obtained from single crossings by Dehn filling the augmentation cusp)
and hence with all hyperbolic knots in $ \IS^3 $.

\begin{figure}
  \centering
  \includegraphics[width=.8\textwidth]{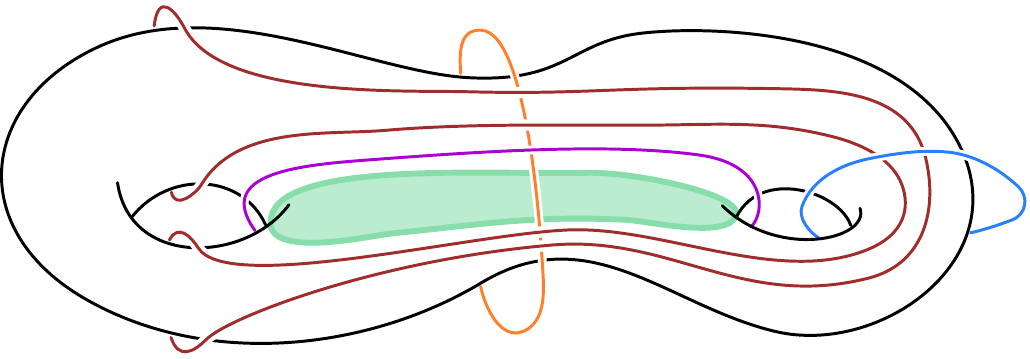}
  \caption{Constructing a tunnel number $1$ link in $ \IS^3 $ by Dehn filling a handlebody (exterior to the genus $2$ surface shown) and gluing a $2$-handle (in green). For simplicity of the
  figure we do not show any stacking of the loops to be Dehn filled, but they should interleave as in \cite[Figure~8]{elzenaar25ej}. We also do not show the loops to perform rational braid
  moves on the left-hand handle of the handlebody, but these loops will be mirror images of those shown for the right-hand handle in light blue and red.\label{fig:twist_regions_s3_1}}
\end{figure}

\begin{figure}
  \centering
  \includegraphics[width=.5\textwidth]{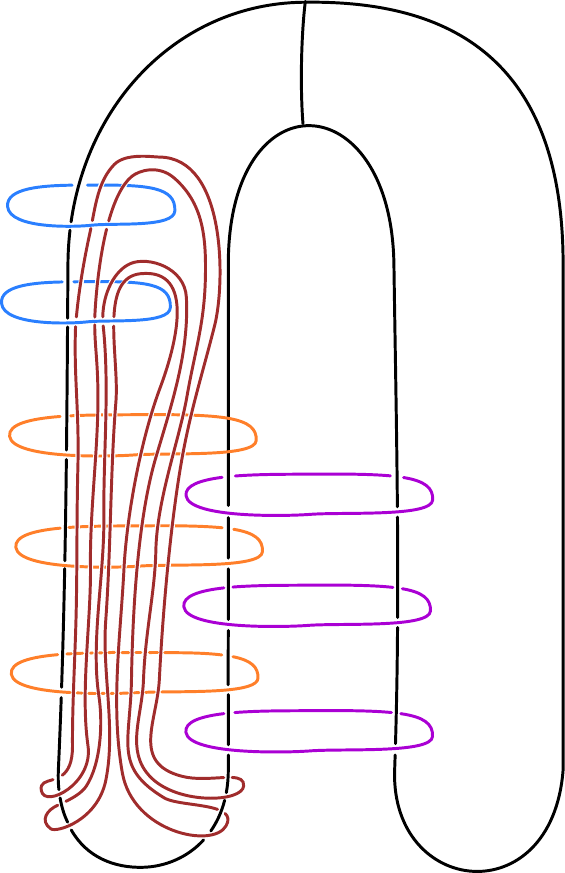}
  \caption{Na\"ive generalisation of a fully augmented link diagram for a $2$-bridge link to arbitrary tunnel number $1$ links in $\IS^3$. Note that the `augmentation cusps' which
  are to be Dehn filled do not bound twice-marked discs in the desired way, so the result is not fully augmented. However, the crossings of the diagram are still heavily
  controlled.\label{fig:twist_regions_s3_2}}
\end{figure}

In the case of $2$-bridge links, there is a very well-known twist diagram and therefore a somewhat canonical choice of Dehn filling presentation, see \cite[Figure~10.2]{purcell}
and also Elzenaar~\cite[\S 3]{elzenaar25ej}. We can use the results of this section to extend this Dehn filling presentation to all tunnel number $1$ links in the sphere. If all
of the Dehn half-twists in the topological intepretation of \zcref{cor:in_sphere} are performed an even number of times, then we have a Dehn filling presentation of a tunnel
number $1$ link $ \mf{k} \subset \IS^3 $ given by Dehn filling loops in $ \mc{H}^+ $ constructed as in \zcref{fig:twist_regions_s3_1} and then gluing in a $2$-handle transverse to the handle.
We can then give a generalisation of the fully augmented link diagrams for $2$-bridge links given by augmenting four planar strands and then capping off the top and bottom. The result
(see \zcref{fig:twist_regions_s3_2}; the diagram can be continued by further stacking loops at the bottom and top of the figure in similar patterns) will not be a fully augmented link
as the discs bounded by the Dehn filling loops no longer have the right intersection numbers (the light blue loops cannot be pulled away from the red loops in the figure without
interacting in some complicated way with the orange loops), but at least the combinatorics of the diagram is highly controlled. An interesting future research direction could be
to try to describe all possible diagrams of this form for all tunnel number $1$ links in $ \IS^3 $ and then to use diagram-theoretic methods like those in Purcell~\cite{purcell}
to triangulate the link exteriors.

\section{The fan of simplices intersecting in a vertex}\label{sec:fan_vertex}
Suppose that two adjacent simplices $ \Delta $ and $ \Delta' $ in $ \KK^{(2)}(S_{2,0})$ intersect in a single vertex. This corresponds to fixing a single
curve on $\Sigma$ and varying the other two, or (equivalently) moving around in the curve complex of the two-holed torus. The combinatorics
of this complex have previously been studied from the point of view of curve coding by Parker, Keen, and Series~\cite{keen99}, but we will take a different point of view.

Let us first describe the sort of construction which we need. We require a combinatorial (graph) structure on the set of all $2$-simplices of $ \KK(S_{2,0}) $
intersecting in a single vertex, or on the $1$-simplices of the curve complex of the twice-punctured torus $\Sigma^* = S_{1,2} $,
where edges do not necessarily come from any incidence data of $\KK(S_{2,0})$. This structure should come with some form of
ordering (like the ordering in the Farey tree induced by `length of continued fraction decomposition') with the property that the vertices of the $2$-simplices that appear later in the ordering have words
that can be written in terms of earlier $2$-simplices and in terms of the vertices of some choice of dual simplices; these dual simplices should also be defined
recursively in terms of words on earlier simplices.

\subsection{The tripants graph}
\begin{figure}
  \centering
  \includegraphics[width=.5\textwidth]{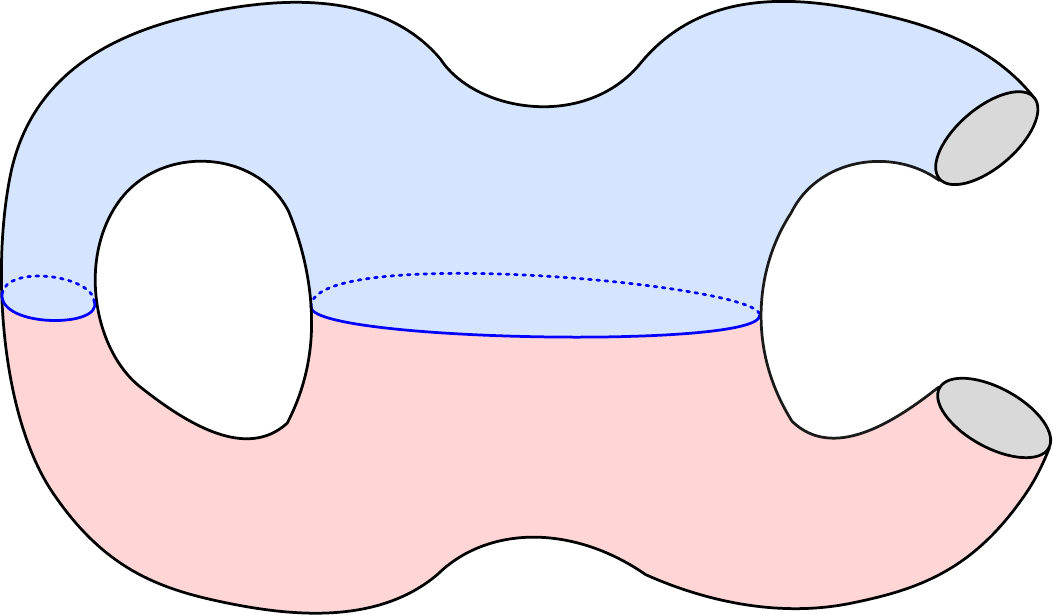}
  \caption{A $1$-simplex of $ \KK(\Sigma^*) $ defines a pants decomposition.\label{fig:twice_punctured_torus_pants}}
\end{figure}

To find a possible construction, note that a $1$-simplex of $ \KK(\Sigma^*) $ is just a choice of pants decomposition of the surface (\zcref{fig:twice_punctured_torus_pants}), and that we can try to mimic
the graph structure on pants decompositions of $ S_{1,1} $ (where two pants decompositions are adjacent iff their boundaries intersect in a single
point). Namely, we place a simplicial complex structure on the set of pants decompositions of $ \Sigma^* $ such that a set $\Delta$ of vertices forms
a simplex if and only if whenever $ \{\alpha, \alpha'\} $ and $ \{\beta,\beta'\} $ are the pairs of simple closed curves defining distinct elements of $ \Delta$,
\begin{displaymath}
  i(\alpha,\beta) = i(\alpha',\beta) = i(\alpha,\beta') = i(\alpha',\beta') = 1
\end{displaymath}
where $ i(\cdot,\cdot) $ denotes intersection number. Unlike the usual curve complex simplicial structure, this structure will turn out to have a
connected $2$-skeleton (i.e.\ graph where vertices are $2$-simplices and edges connect vertices that intersect in $1$-simplices) which is closely connected to $ \mc{F} $
and so has the properties we asked for above. We will now define this dual graph in terms of the combinatorics of the surface, and then prove its properties.

\begin{defn}\label{defn:tripants_graph}
  Let $ \Sigma^*$  be a twice-holed torus. The \df{tripants graph} on $\Sigma^*$ is the graph $ \TP(\Sigma^*) $ given by the following data.
  \begin{itemize}
    \item The vertex set of $\TP(\Sigma^*)$ is the set of unordered triples
      \begin{displaymath}
        \big\{\{\alpha_1,\alpha_1'\}, \{\alpha_2,\alpha_2'\}, \{\alpha_3,\alpha_3'\}\big\}
      \end{displaymath}
      where for each $ i $, $ \alpha_i $ and $ \alpha_i' $ are disjoint nonseparating simple closed curves on $ \Sigma^* $ that define a pants decomposition of the surface
      and where the intersection matrix of the six curves is (missing values are $0$):
      \begin{displaymath}\setlength\arraycolsep{2pt}
        \begin{matrix}
                       & \alpha_1 & \alpha_1' & \alpha_2 & \alpha_2' & \alpha_3 & \alpha_3'\\\addlinespace
             \alpha_1  &   &   & 1 & 1 & 1 & 1\\
             \alpha_1' &   &   & 1 & 1 & 1 & 1\\
             \alpha_2  & 1 & 1 &   &   & 1 & 1\\
             \alpha_2' & 1 & 1 &   &   & 1 & 1\\
             \alpha_3  & 1 & 1 & 1 & 1 &   &  \\
             \alpha_3' & 1 & 1 & 1 & 1 &   &
        \end{matrix}
      \end{displaymath}
    \item Two vertices are joined by an edge if and only if they share two out of the three pants decompositions which they bound.
  \end{itemize}
  See the introduction to \cite{betts21}.
\end{defn}

\begin{figure}
  \centering
  \includegraphics[width=.6\textwidth]{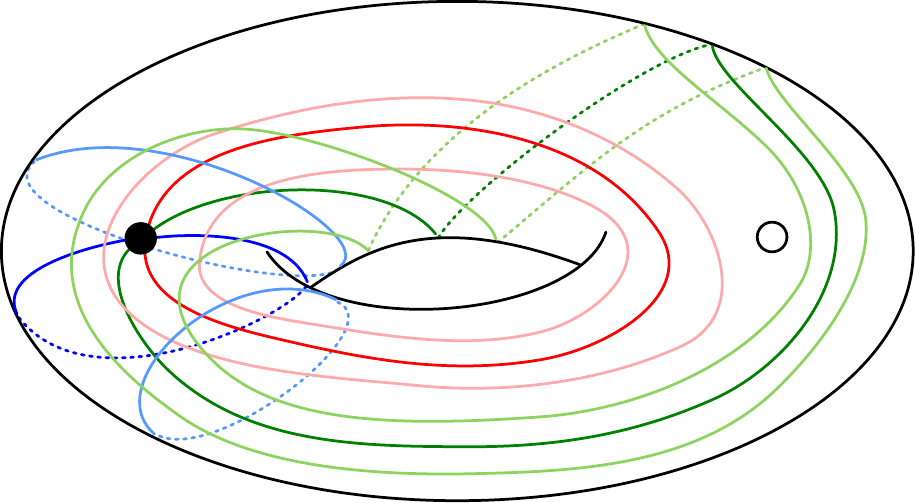}
  \caption{A tripants (bounded by the six lighter curves) and the corresponding triarc through $ \bullet $ (the three darker curves).\label{fig:tripants}}
\end{figure}

\begin{figure}
  \begin{subfigure}{.33\textwidth}\centering
    \begin{tikzpicture}[scale=.6]
      \begin{scope}[very thick,decoration={
          markings,
          mark=at position 0.5 with {\arrow{>}}}
          ]
          \draw[postaction={decorate}, red] (0,0)--(4,0);
          \draw[postaction={decorate}, blue] (4,0)--(4,4);
          \draw[postaction={decorate}, red] (0,4)--(4,4);
          \draw[postaction={decorate}, blue] (0,0)--(0,4);
          \draw[postaction={decorate}, green] (0,0) to[bend left=20] (4,4);

          \node at (0,0)[circle,fill,inner sep=1.5pt]{};
          \node at (0,4)[circle,fill,inner sep=1.5pt]{};
          \node at (4,0)[circle,fill,inner sep=1.5pt]{};
          \node at (4,4)[circle,fill,inner sep=1.5pt]{};
          \node at (2,2)[circle,inner sep=1.5pt]{$\circ$};
      \end{scope}
    \end{tikzpicture}
    \caption{}
  \end{subfigure}\hspace{1cm}%
  \begin{subfigure}{.33\textwidth}\centering
    \begin{tikzpicture}[scale=.6]
      \begin{scope}[very thick,decoration={
          markings,
          mark=at position 0.5 with {\arrow{>}}}
          ]
          \draw[postaction={decorate}, red] (0,0)--(4,0);
          \draw[postaction={decorate}, blue] (4,0)--(4,4);
          \draw[postaction={decorate}, red] (0,4)--(4,4);
          \draw[postaction={decorate}, blue] (0,0)--(0,4);
          \draw[postaction={decorate}, green] (0,0) to[bend right=20] (4,4);

          \node at (0,0)[circle,fill,inner sep=1.5pt]{};
          \node at (0,4)[circle,fill,inner sep=1.5pt]{};
          \node at (4,0)[circle,fill,inner sep=1.5pt]{};
          \node at (4,4)[circle,fill,inner sep=1.5pt]{};
          \node at (2,2)[circle,inner sep=1.5pt]{$\circ$};
      \end{scope}
    \end{tikzpicture}
    \caption{}
  \end{subfigure}
  \caption{The two lifts of an ideal triangulation of $ S_{1,1} $ to a triarc on $ S_{1,2} $ (c.f.\ \zcref{fig:farey_top_upstairs}).\label{fig:farey_to_tripants}}
\end{figure}
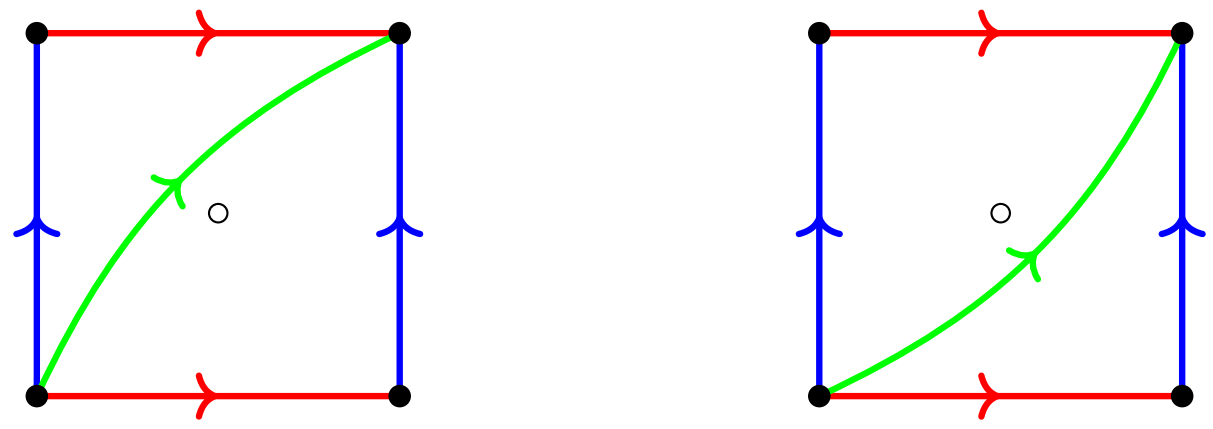

The following is Theorem~A of \cite{betts21}.
\begin{lem}\label{lem:tripants_connected}
  The graph $ \TP(\Sigma^*) $ is connected.
\end{lem}
\begin{proof}
  We sketch the proof, for the details see \cite{betts21}. View $ \Sigma^* $ as a torus $T$ with two marked points $ \circ $ and $ \bullet $.
  A \df{triarc} on $ \Sigma^* $ is a triple of distinct isotopy classes of curves through $ \bullet $ on $ T \setminus \{\circ\} $ which have mutual intersection
  number $1$. Given any triarc, take a small annular neighbourhood of each curve; each annulus forms a pair of pants on $ \Sigma^* $ where one of the
  `holes' is the puncture $ \bullet $ and the boundaries of these annuli give six curves on $ \Sigma^* $ which have the correct intersection number matrix
  to be a tripants (\zcref{fig:tripants}). This map from triarcs to tripants is a bijection, and edges in the tripants graph correspond to pairs of triarcs which share two out of their
  three loops. This gives a graph $ \TA(\Sigma^*) $ of triarcs which is isomorphic to $ \TP(\Sigma^*) $ and it suffices to prove connectedness of $ \TA(\Sigma^*) $.

  Consider next the graph with vertices consisting of triples of isotopy classes of loops on $ T $ through $\bullet$ with mutual intersection number $1$, and where
  two vertices are joined by an edge if the corresponding triples share two out of three loops. This is exactly the Farey graph $ \mc{F} $, c.f.\ \zcref{defn:farey_graph},
  so the map $ T \setminus \{\circ,\bullet\} \to T \setminus \{\bullet\} $ induces a map
  \begin{displaymath}
    \phi : \TA(\Sigma^*) \to \mc{F}.
  \end{displaymath}
  The connectedness of $ \TA(\Sigma^*) $ now follows from the connectedness of $ \mc{F} $, together with the observation that
  if $ \phi(x) = \phi(y) $ (where $x$ and $y$ are triarcs) then either $ x = y $ or $ x $ and $y$ form an edge in $ \TA(\Sigma^*)$.
  This is because a choice of lift of an ideal triangulation $ I $ of the once-punctured torus $ T \setminus \{\bullet\} $ through
  the map $ \phi $ is the same as a choice of which of the two components of $ T \setminus I $ to puncture as shown in \zcref{fig:farey_to_tripants},
  and the two choices form an edge in the triarc graph since they will share two out of the three ideal arcs: in the figure, the vertical and horizontal
  ideal arcs.
\end{proof}

We extract the following from the proof of \zcref{lem:tripants_connected}.
\begin{prp}\label{prp:tripants_summary}
  The tripants graph is the dual graph to a simplicial complex supported on the set of pants decompositions of $ \Sigma^* $. It
  is a cover of $\mc{F}$, with $ \TP(\Sigma^*) \to \mc{F} $ given by contracting certain edges. In particular, the preimage of any vertex of $ \mc{F} $ under
  the covering map consists of at most two vertices. Further, if the preimage consists of two vertices, then they form an edge of the tripants
  graph and become isotopic on the once-punctured sphere obtained by filling in one of the punctures of $ \Sigma^* $. \qed
\end{prp}

\subsection{Word recursion down the tripants graph}
We now relate the tripants graph to the actual problem which we are trying to solve. Recall that our goal is to give a recursion down a ns-geodesic in the set of $2$-simplices of the
complex of curves $ \KK(\Sigma) $, where $ \Sigma $ is a genus two surface, which starts at the trivial simplex $ \Delta_0 $ (\zcref{defn:deltanought}) and which ends at an arbitrary simplex $\Delta$; this path
consists of alternating runs of $2$-simplices which are joined along edges, and runs of $2$-simplices which are joined across a vertex. Further, in the definition of ns-geodesic
we assume that all simplices except possibly for $\Delta$ lie in the \emph{nonseparating} curve complex, and the join between $ \Delta $ and the preceeding simplex is along an edge rather than
a vertex. Thus, in our recursion along a fan of simplices intersecting in a vertex, we may assume that all the simplices have nonseparating curves as vertices.

The recursion should allow us to compute
the words in $ \pi_1(\Sigma) $ corresponding to the vertices of $ \Delta $ in such a way that it is obvious that the group generated by these words is contained in the group generated by the
words associated to the simplex at the end of the previous run (of either edge-intersecting or vertex-intersecting simplices) of $2$-simplices in the path.
In addition, since we need to keep track of an entire surface subgroup, at every step we also need to have in hand a dual simplex to every simplex (\zcref{defn:dual_simplex})
for which the corresponding words are also obtained from previous words in a recursion; it is less important for the path of dual simplices to have any clear relation to each
other, and we only need to make some choice for them which is computable in some way from the data at hand.

We are therefore given a nonseparating $2$-simplex $\Delta \in \KK^{(2)}(\Sigma) $ together with a dual simplex $ \Delta^\perp $ (with the associated words in the surface group $ \pi_1(\Sigma)$)
and a second $2$-simplex $ \Delta' $ such that $ \abs{\Delta \inter \Delta'} = 1 $, and we want to find a natural recursion from $ \Delta $ to $ \Delta' $. We will do this by
producing a point in the tripants graph from $ \Delta $ and $ \Delta^\perp $, walking through the tripants graph (via \zcref{prp:tripants_summary}) to one containing the
pants decomposition encoded by $ \Delta' $ as one of its three pants decompositions, and using the other two pants decompositions to produce a dual simplex for $ \Delta' $. At
each step we will use the identification of the tripants graph with a simple cover of $ \mc{F} $ to compute the words of each simplex.

\begin{figure}\vspace{1.5ex}%
  \begin{subfigure}[t]{.47\textwidth}
    \labellist
    \small\hair 2pt
    \pinlabel {$\sigma$} [b] at 117 300
    \pinlabel {$\beta$} [b] at 220 300
    \pinlabel {$\rho$} [b] at 307 300
    \pinlabel {$\tau$} [b] at 443 300
    \pinlabel {$\alpha$} [r] at 0 153
    \pinlabel {$\gamma$} [l] at 546 158
    \endlabellist
    \centering
    \includegraphics[width=\textwidth]{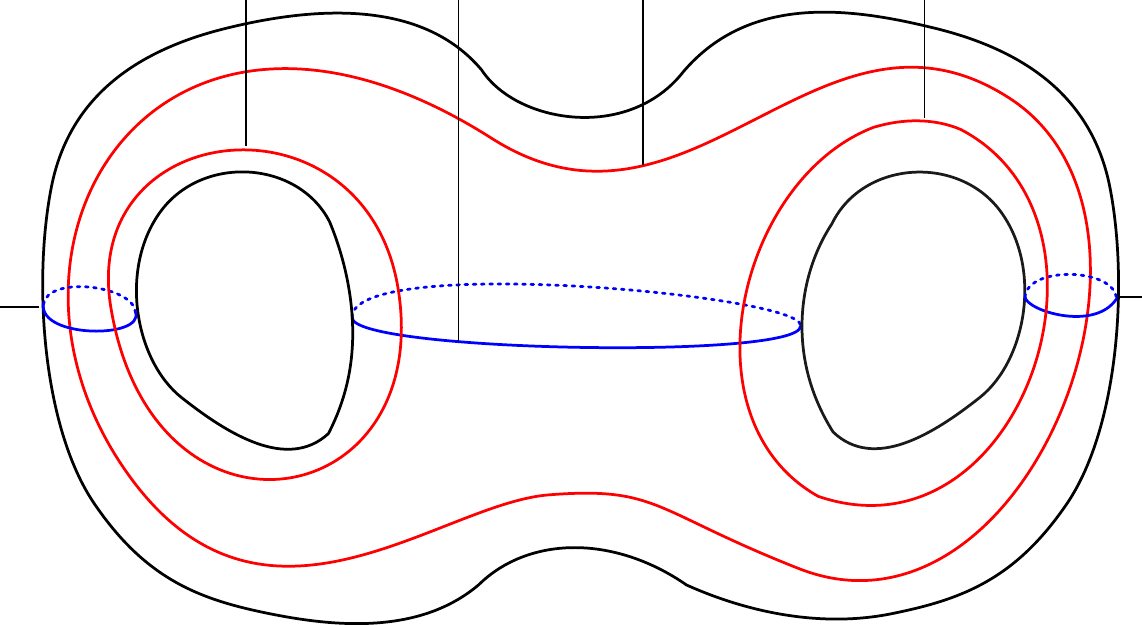}
    \caption{The original pair $ (\Delta,\Delta^\perp)$.}\label{fig:dual_curves_to_tripants_1}
  \end{subfigure}\hfill%
  \begin{subfigure}[t]{.47\textwidth}
    \labellist
    \small\hair 2pt
    \pinlabel {$\sigma$} [b] at 117 300
    \pinlabel {$\beta$} [b] at 220 300
    \pinlabel {$\sigma\gamma$} [b] at 322 300
    \pinlabel {$\alpha$} [r] at 1 153
    \pinlabel {$\gamma$} [l] at 525 158
    \endlabellist
    \centering
    \includegraphics[width=\textwidth]{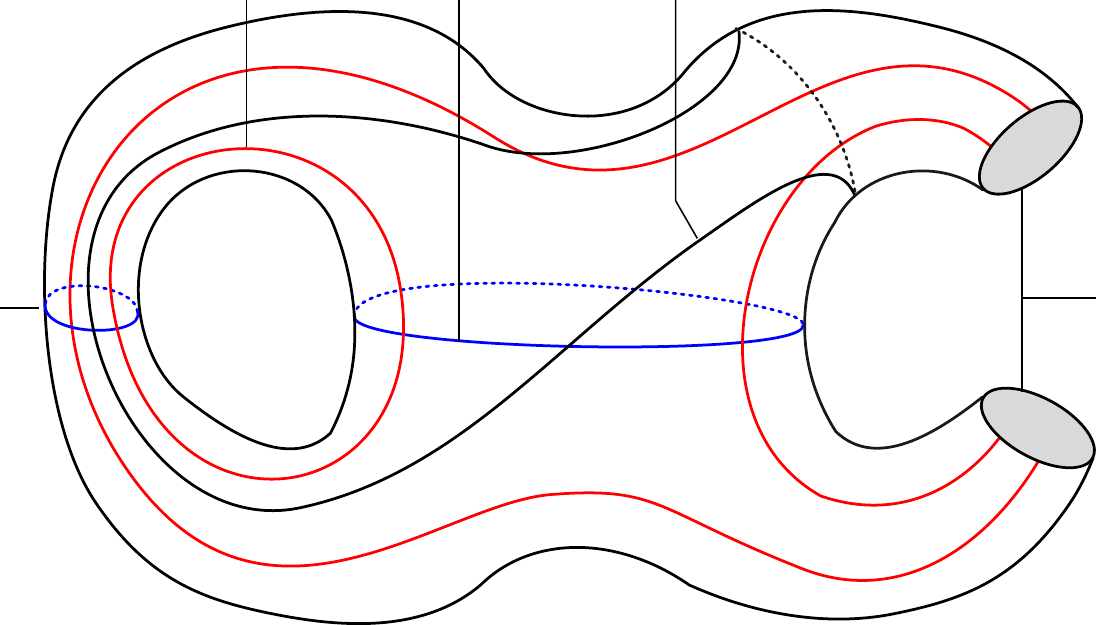}
    \caption{Two pants in the tripants are bounded by $ \{ \alpha, \beta \} $ and $ \{ \sigma, \sigma\gamma \} $.}
  \end{subfigure}
  \caption{Production of two pants in a tripants from a simplex $ \Delta = \{\alpha,\beta,\gamma\} $ with $ \Delta^\perp = \{ \sigma,\tau,\rho \} $
  where $\gamma$ is the curve to be fixed.\label{fig:dual_curves_to_tripants}}
\end{figure}
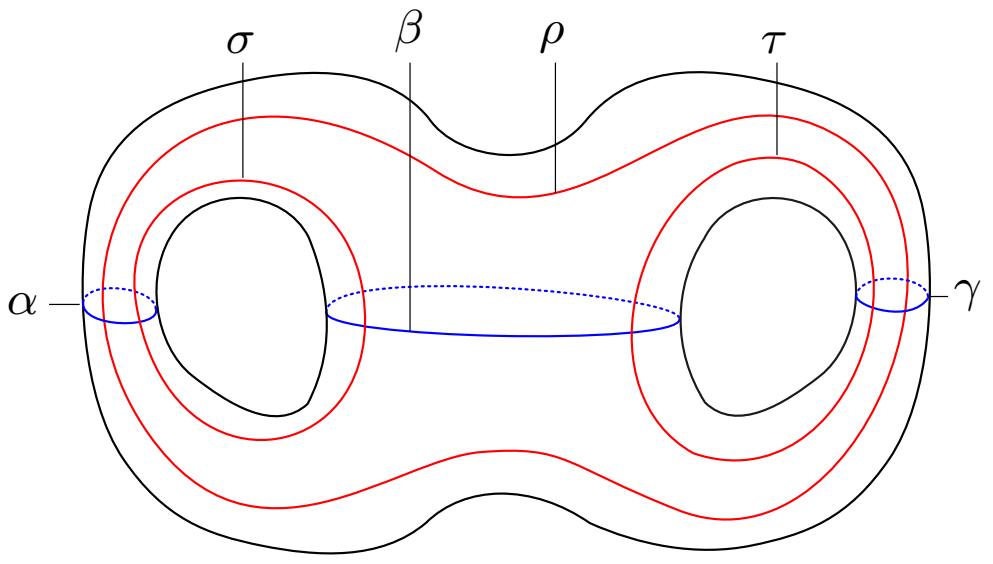
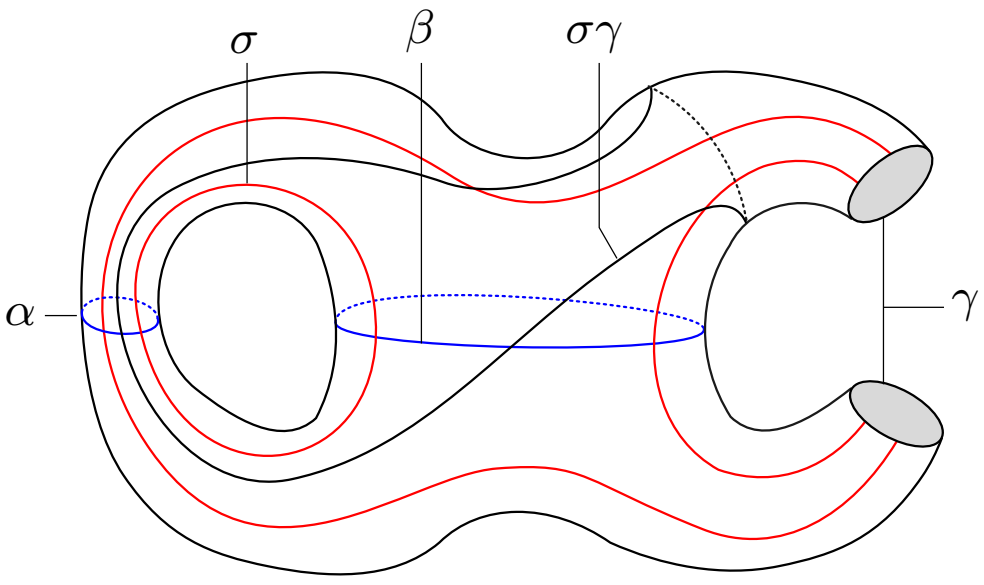

\begin{cons}[Construction of a tripants from a simplex and its dual.]\label{cons:simplex_to_tripants}
  The construction of two pants decompositions with the right intersection data is shown in \zcref{fig:dual_curves_to_tripants}. Suppose that we have a nonseparating simplex
  $ \Delta = \{\alpha,\beta,\gamma\} $ and a dual simplex $ \Delta^\perp = \{ \sigma,\tau,\rho \} $.
  Suppose that $ \gamma $ is the vertex that lies in the intersection of $ \Delta $ with the desired next
  simplex in the path; i.e.\ $\gamma$ is the curve along which we cut $ \Sigma$ to obtain a twice-holed torus $\Sigma^*$.
  Then $ \alpha $ and $ \beta $ together define a pants decomposition of $ \Sigma^* $. In addition, one of the three
  curves in $ \Delta^\perp $ (namely, the unique one which does not intersect $ \gamma$), say $ \sigma $, is still a simple
  closed curve on $ \Sigma^* $ which intersects each of $ \alpha $ and $ \beta $ once. We next find a second curve which
  intersects $ \alpha $ and $ \beta $ once and which is disjoint from, and not isotopic to, $ \sigma $. Direct inspection shows
  that the `twist' of $ \sigma $ given by $ \sigma\gamma $ works.

  A third pants decomposition by nonseparating curves is needed to complete the tripants. Such a decomposition is given by $ \{ \beta\sigma, \beta\sigma\gamma \} $. By a standard
  argument involving transitivity properties of $ \Mod(\Sigma^*) $, it is essentially enough to check that this works in the configuration shown in \zcref{fig:dual_curves_to_tripants}.
\end{cons}

Now we consider the recursion of words bounding tripants as we walk along the tripants graph to a tripants which contains $ \Delta' $ as one of the three pants decompositions.
For this, we use the map $ \Sigma^* \to S_{1,1} $ as in the proof of \zcref{lem:tripants_connected}.

On the once-punctured torus $S_{1,1}$ with longitude and meridian $ A $ and $ B $, we have a natural map $ \mc{F} \to \pi_1(S_{1,1}) $ written as $ p/q \mapsto \Jorgensen_{p/q} $
so that $ A = \Jorgensen_{0/1} $, $ B = \Jorgensen_{1/0} $, and
\begin{equation}\label{eq:jorgensen_recursion}
  \Jorgensen_{p/q \oplus r/s} = \Jorgensen_{p/q} \Jorgensen_{r/s}.
\end{equation}
The proof is very similar to that of \zcref{lem:product_word} above, and may be found in Chapter~9 of Mumford, Series, and Wright~\cite{indras_pearls}.
Here, `TW' stands for `\textbf{T}orus \textbf{W}ord'.

By the Siefert--van Kampen theorem, the fundamental group $ \pi_1(\Sigma^*) $ is $ \langle A,B,K \rangle $ where $K$ is a primitive loop
around a distinguished puncture $\bullet$. (In the setup of \zcref{cons:simplex_to_tripants}, $ K $ is represented by $ \gamma $.) Choose a spanning
tree $ \TP^\dagger(\Sigma^*)$ of $ \TP(\Sigma^*) $ and colour the edges white if they are lifts of edges from $ \mc{F} $ and black if they are contracted when
projecting to $ \mc{F} $. By \zcref{prp:tripants_summary}, every edge is coloured and no black edges are incident.
\begin{lem}\label{lem:tripants_rewriting}
  There is a recurrence on words representing disjoint simple closed curves on the twice-punctured sphere, given by the following rule:
  \begin{enumerate}
    \item Walking along a white edge of $ \TP^\dagger(\Sigma^*) $: fill in the distinguished puncture $ \bullet $ and apply \eqref{eq:jorgensen_recursion}
          to obtain a recursion on the resulting once-punctured torus.
    \item Walking along a black edge of $ \TP^\dagger(\Sigma^*) $: replace one instance of $B$ with $ KBK^{-1} $, namely the single instance corresponding
          to the local piece of loop which is moved across $ \bullet$ by isotopy. \qed
  \end{enumerate}
\end{lem}

Once we have reached the final vertex in the tripants graph, namely a vertex which contains $ \Delta' = \{ \alpha',\beta',\gamma \} $ as one of the three
pants decompositions, we need to be able to read off a dual simplex to $ \Delta' $ and be able to write down the words for this dual simplex in terms
of the original words in $ \Delta $.

\begin{cons}[Construction of a simplex and its dual from a tripants.]\label{cons:tripants_to_simplex}
  Suppose we are given a tripants containing three pants decompositions of $ \Sigma^* $, including our objective simplex $ \Delta' = \{ \alpha',\beta',\gamma \} $.
  We need to construct a dual simplex to $ \Delta' $. Pick any curve (arbitrarily) from one of the other three pants decompositions, say $ \sigma' $. Then $ \sigma' $
  intersects both $ \alpha' $ and $ \beta' $ exactly once. We now follow the same idea as \zcref{lem:product_word_dual}. First, expand the `radius' of $ \sigma' $ by isotopy
  until it hits both punctures of $ \Sigma^* $---i.e.\ we find a pair of disjoint  ideal arcs on $ \Sigma^* $, say $ \nu $ and $ \mu $, with the property that the closure of $ \nu \union \mu $
  is isotopic to $ \sigma' $ on the two-point compactification of $ \Sigma^*$ (similarly to the once-punctured torus setting shown in \zcref{fig:farey_top_downstairs}
  and \zcref{fig:farey_top_upstairs} above). By the fact that $ \sigma' $ is obtained by a word recursion from the original simplices $ \Delta $
  and $ \Delta^\perp $, just as in part (3) of \zcref{lem:product_word_dual} the cutting sequences of $ \nu $ and $ \mu $ come from a factorisation of $ \sigma' $ in terms of the original
  generators $ \{ \alpha,\beta,\sigma,\gamma \} $. Thus the desired dual simplex to $ \Delta' $ can be chosen to be $ \{ \sigma', \nu \tau, \mu \tau \} $ where $ \tau $
  is the dual curve in $ \Delta^\perp $ which represents the handle that was cut when passing from $ \Sigma $ to $ \Sigma^* $ (c.f.\ \zcref{fig:dual_curves_to_tripants_1})---in other
  words, the letter in the alphabet $ \Delta^\perp $ which is appended to words when crossing the curve $ \gamma $ on $ \Sigma $.
\end{cons}

\subsection{Application to the filtration theorem}
We may now extend \zcref{thm:filtration_in_component} to the whole curve complex.
\begin{thm}\label{thm:filtration_global}
  Let $ \Delta_0, \ldots, \Delta_k $ be a ns-geodesic of $2$-simplices in the non-separating curve complex of the genus $2$ surface $\Sigma = \partial \mc{H} $, such that
  \begin{enumerate}
    \item $ 0 \leq \dim (\Delta_i \inter \Delta_{i+1}) \leq 1 $ for each $ 0 \leq i < k $, and
    \item If $ \abs{i-j}>1 $ then $ \dim (\Delta_i \inter \Delta_j) < 1 $.
  \end{enumerate}
  Let $ \Delta_i^\perp$ be defined by \zcref{lem:product_word_dual} if $ \dim (\Delta_i \inter \Delta_{i-1}) = 1 $, and by \zcref{cons:tripants_to_simplex}
  if $ \dim (\Delta_i \inter \Delta_{i-1}) = 0 $.

  Let $ \rho : \pi_1(\mc{H}) \to \PSL(2,\IC) $ be a representation. Then the sequence of subgroups
  \begin{displaymath}
    \Gamma_{k-i} = \langle \Delta_{k-i}, \Delta_{k-i}^\perp \rangle \coloneq \langle \rho(\gamma) \in V(\Delta_{k-i}) \rangle
  \end{displaymath}
  is a filtration of $ \Gamma_k $, i.e.
  \begin{displaymath}
    \Gamma_0 < \cdots < \Gamma_k,
  \end{displaymath}
  and each $ \Gamma_i $ admits a generating set $ S_i = \{ U_i, V_i, F_i, G_i \} $ such that $ S_i $ is obtained from $ S_{i-1} $
  either via \zcref{thm:filtration_in_component} (if $ \dim (\Delta_i \inter \Delta_{i-1}) = 1 $), or via a sequence of rewritings
  as in \zcref{lem:tripants_rewriting} followed by construction of two duals as in \zcref{cons:tripants_to_simplex} (if $ \dim (\Delta_i \inter \Delta_{i-1}) = 0 $). \qed
\end{thm}

\section{Tunnel number one manifolds}\label{sec:tno_general}
We now connect the theory of the previous sections with the theory of tunnel number $1$ manifolds.

\subsection{Links in arbitrary closed manifolds}
Suppose that $ \mf{k} \subset A $ (where $A$ is a suitable closed $3$-manifold) is tunnel number $1$, and so it corresponds to a genus two Heegaard diagram $\Delta$ as explained
in \zcref{lem:tunnel_number_1_is_heegaard} together with a choice of one of the three curves to bound the disc dual to the unknotting tunnel. Suppose we choose a ns-geodesic
in the curve complex from $ \Delta_0 $ to $ \Delta $ as we have described several times. The only thing we need to be careful of is such a path ending up passing through a
fan of simplices joined along a single vertex in such a way as one of the vertices becomes a nonseparating curve, since our arguments involving the tripants graph involved
considering only nonseparating curves. This is excluded by the definition of ns-geodesic, which requires the separating curve to appear only at the final stage via a rational
braiding move.
\begin{cor}\label{cor:in_arbitrary}
  There exists a filtration of the holonomy group of a tunnel number $1$ link complement $ A \setminus \mf{k} $ determined by a choice of ns-geodesic in the graph of $2$-simplices
  of $ \KK(\Sigma) $ between the canonical basepoint $\Delta_0$ (\zcref{defn:deltanought}) and the Heegaard diagram $\Delta_n$ associated with $ \mf{k} \union \tau $, where $ \tau $ is a choice
  of tunnel. In particular we can write
  \begin{displaymath}
    \Gamma_0 < \cdots < \Gamma_n = \pi_1(A \setminus \mf{k}) < \PSL(2,\IC)
  \end{displaymath}
  where each $ \Gamma_i $ is a $ S_{2,0}$-surface representation. Each $ \Gamma_i $ is the homotopy group of an embedded
  topological genus $2$ surface, and $ \Gamma_0 $ is $ V_4 \coloneq \IZ \oplus \IZ $ (resp.\ $ V_4 * V_4$) if $ \mf{k} $ has one
  component (resp.\ two components). \qed
\end{cor}

\begin{figure}\vspace{1.5ex}%
  \labellist
  \small\hair 2pt
  \pinlabel {longitude} [b] at 109 305
  \pinlabel {$\mu_1$} [b] at 21 305
  \pinlabel {$\mu_2$} [b] at 220 305
  \endlabellist
  \centering
  \includegraphics[width=.6\textwidth]{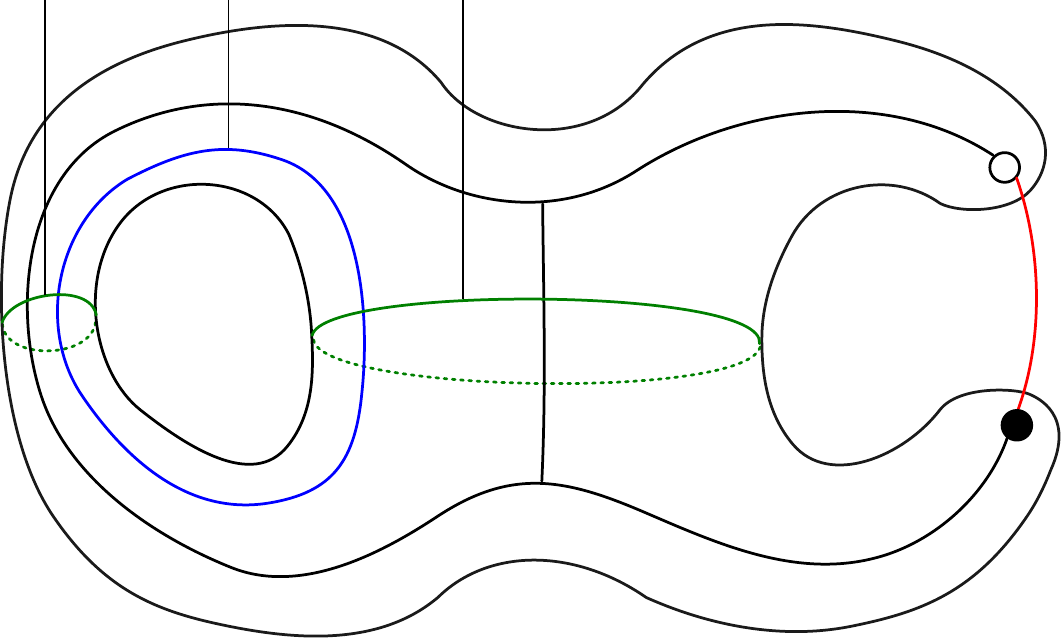}
  \caption{Tripants moves correspond to Dehn fillings along tori which intersect the $\theta$-curve $ \mf{k} \union \tau $ (here, the union of the red and black strands) twice. Instead of a standard Dehn filling where
  one keeps track only of the position that the meridian ends up, here there are two possible choices as $ \mu_1 $ and $ \mu_2 $ are non-isotopic on the twice-punctured torus.
  Thus the filling is parameterised both by a slope and a `colour' as in \zcref{lem:tripants_rewriting}.\label{fig:dehn_surgery}}
\end{figure}

In terms of the knot theory, we have already seen that within a path consting of fans of simplices joined along edges the filtration models a sequence of rational braiding moves
and the topological genus two surfaces separate individual rational braidings. Tripants moves, on the other hand, correspond to Dehn surgeries of the ambient manifold $A$ along
a torus which intersects a single curve of knot or tunnel in two points. The difference between black and white edges in the tripants graph (\zcref{lem:tripants_rewriting})
corresponds to the two choices of position for the meridian after gluing, since now there are two punctures there are two isotopy classes of meridian (the green
loops in \zcref{fig:dehn_surgery}).

Of course our methods depend on choosing a tree structure on the set of $2$-simplices of $ \KK(\Sigma) $. Within each connected component after deleting vertices
we can do this canonically, but as soon as we start requiring such Dehn surgery moves we must make choices. Since $ \KK(\Sigma) $ is a $\delta$-hyperbolic metric space~\cite{masur99} (one
can take $\delta=17$~\cite{hensel15}), any choices that we do make must stay within a uniformly bounded distance of each other. However, we still cannot answer the following.
\begin{question}
  Given a pair of a knot and tunnel in a closed $3$-manifold $A$ represented by some $ \Delta \in \KK^{(2)}(\Sigma) $, what is the minimal number of Dehn surgeries on $A$
  which appear in a path from $ \Delta $ to $ \Delta_0 $?
\end{question}

We wish to note the following implication of this point of view to the topology of the curve complex.
\begin{prp}
  Let $ \KK^{(0)}(\Sigma) $ be the set of vertices. Then $ \KK(\Sigma) \setminus \KK^{(0)}(\Sigma) $ (as a topological space) is disconnected.
\end{prp}
\begin{proof}
  Suppose not; then one could walk from any $2$-simplex to any other by a path that lies entirely in the interior of $2$- and $ 1$-simplices. This
  means that one could walk from the knot represented by $ \Delta_0 $ to any other tunnel number $1$ knot by rational braiding moves. Rational braiding moves stay entirely
  within the same ambient manifold. Thus  every tunnel number $1$ link would be embeddable in $\IS^3$, contradicting (e.g.) the existence of \texttt{m003(0,0)}
  (\zcref{ex:m003}) which is known to be unembeddable in $ \IS^3 $ by the discussion in Example~2.1 of \cite{kegel24}.
\end{proof}

\begin{question}\label{qn:nested}
  Can the nested topological surfaces of \zcref{cor:in_arbitrary} be straightened to nested pleated surfaces, pleated along simple closed curves?
\end{question}
This would have implications for the combinatorial structure of tunnel number $1$ links, as one would like to generalise
the Sakuma--Weeks triangulation for two-bridge links (see \zcref{sec:two_bridge}) to a geodesic triangulation of tunnel number $1$ links by individually
triangulating every interstitial region between these surfaces compatibly with the pleating locus; one would then hope that the result is related to the Ford
domain of the group $\Hol(A \setminus \mf{k}) $ and thus to the canonical triangulation in the sense of Epstein and Penner~\cite{epstein88}.

We can prove only the following weak version of \zcref{qn:nested}:
\begin{thm}\label{thm:nested_weak}
  Let $M=A \setminus \mf{k}$ be a hyperbolic tunnel number $1$ link complement.
  The nested topological surfaces of \zcref{cor:in_arbitrary} be straightened to nested pleated surfaces.
\end{thm}
\begin{proof}
  Since the groups $ \Gamma_i $ are nested, their limit sets are nested. Thus it suffices to prove that every $ \Gamma_i $ is geometrically finite,
  for then we can take the surfaces $ \hconv \Gamma_i/\Gamma_i $ as the desired sequence of pleated surfaces.

  To prove geometric finiteness of the groups in the filtration, we rely on a generalisation of a result of Thurston and Bonahon, namely \zcref{thm:embedded_geom_finite}.
  Since the proof of this result uses entirely different techniques (and is likely to be routine for experts in the theory of pleated surfaces and ends of hyperbolic
  manifolds) we have placed this result in an appendix (\zcref{app:pleated_surfaces}). For us, the main consequence of this result is that, since every surface $ S_i $
  is compressible without trivial handles, the manifold $M$ cannot virtually fibre over the circle with fibres $S_i$ (the compression discs cannot be foliated) and therefore
  the corresponding holonomy groups $ \Gamma_i $ must be geometrically finite.
\end{proof}
Thus the real content of \zcref{qn:nested} is to understand the pleating coordinates of each $ \Gamma_i $. In particular, to understand the combinatorial structure
of the manifold one would like these pleating coordinates to arise from vertices of a geodesic in the curve complex, but it is unclear how tightly the pleated
surfaces we construct in \zcref{thm:nested_weak} hug the geodesics in $M$ corresponding to those vertices. This is a similar problem to that which arises in the theory
of ending laminations. In fact, our entire construction is motivated by taking the analogy between the model manifold theorem for $ S_{0,4} $~\cite{ohshika10}
and the Sakuma--Weeks triangulation~\cite{sakuma95} and transporting it to the more general setting of Schottky space. This is similar to the point of view of Heegaard splittings taken
by Namazi in his doctoral dissertation~\cite{namazi05}. In the ending lamination setting, one takes (relatively) short geodesics and `stretches' pleated surfaces across
them; here there are no `ends' but one can still form these pleated surfaces. It is unclear how these are related to the surfaces we construct in \zcref{thm:nested_weak}.

\subsection{Tunnel number one manifolds}\label{sec:tno_manifold_general}
Recall (\zcref{defn:tno_manifold}) that a tunnel number one manifold $M$ is obtained from Dehn filling between zero and two of the cusps of a tunnel one link $\mf{k} $ in
a $3$-manifold $A$. Each simplex in a family of $2$-simplices in $ \KK(\Sigma) $ which share a vertex $ v $ represents the same
pair $ (\mf{k},\tau) $ of link and tunnel. The two other vertices of each simplex correspond to simple closed curves on the tori boundary component(s) of $ A \setminus \mf{k} $.
If $ \mf{k} $ has two components, then these two closed curves index Dehn fillings of the two components (namely, glue a $2$-handle along each curve and fill the remainder with a $3$-ball).
If $ \mf{k} $ has one component, then the two simple closed curves are non-intersecting on the torus, and hence must be freely homotopic. Thus we again get an indexing of Dehn
fillings by the simplices. We see that the $2$-simplices of the curve complex $ \KK(\Sigma) $ can be viewed as indexing tunnel number $1$ manifolds, and so our filtrations
also give filtrations of all tunnel number $1$ manifold holonomy groups.

\begin{thm}\label{thm:tno_filtration}
  Let $ M $ be any hyperbolic tunnel number $1$ manifold, determined by a $2$-simplex $ \Delta \in \KK^{(2)}(\Sigma) $ together with a choice of between one and three of its vertices.
  There exists a filtration of the holonomy group of $ M $ determined by a choice of path in the graph of $2$-simplices
  of $ \KK(\Sigma) $ between the canonical basepoint $\Delta_0$ (\zcref{defn:deltanought}) and $ \Delta$,
  \begin{displaymath}
     \Gamma_0 < \cdots < \Gamma_n = \pi_1(M) < \PSL(2,\IC),
  \end{displaymath}
  where each $ \Gamma_i $ is a $ S_{2,0}$-surface representation. Each $ \Gamma_i $ is the homotopy group of an embedded
  topological genus $2$ surface and $ \Gamma_0 $ is a quotient of $ V_4 \coloneq \IZ \oplus \IZ $ by killing a generator (resp.\ $ V_4 * V_4$ by killing one or two generators) if
  the tunnel presentation of $ M $ is the Dehn filling of a knot (resp.\ a link of two components). \qed
\end{thm}

As we have already discussed, the choice of path in this theorem is essentially unique in the sense that all choices lie within some bounded distance of each other.
In addition, the proof of \zcref{thm:nested_weak} continues to hold after Dehn filling and so the filtration surfaces can be straightened to nested pleated surfaces.

\appendix
\section{Tameness and geometric finiteness of embedded surfaces}\label{app:pleated_surfaces}
The purpose of this appendix is to give a more general version of a standard result about embedded surfaces in geometrically finite hyperbolic $3$-manifolds (\zcref{thm:embedded_geom_finite}). We
first state the standard result which can be found in, for instance, Purcell~\cite[Theorem~12.15]{purcell}:
\begin{thm}\label{thm:bers_essential}
  Let $S$ be an essential surface in a geometrically finite hyperbolic $3$-manifold $M$. Then one of the following three situations occurs:
  \begin{enumerate}
    \item $\Hol(S)$ is a quasi-Fuchsian group embedded in $\Hol(M)$.
    \item $\Hol(S)$ has an accidental parabolic.
    \item $M$ is virtually $S$-fibred: there is a finite cover of $M$ which is fibred over the circle, with fibres equal to lifts of $S$.
  \end{enumerate}
\end{thm}
\zcref[S]{thm:bers_essential} is obtained by combining two results about realisation of laminations by pleated surfaces in the homotopy class of $S$.
\begin{lem}\label{lem:classical_compactness}
  Let $S$ be an essential surface in a geometrically finite hyperbolic $3$-manifold $M$ and suppose $ \Hol(S) $ has no accidental parabolics.
  Then either $M$ is virtually $S$-fibred, or the entire geodesic lamination space of $S$ is realisable by pleated surfaces in $M$ in the homotopy class of $S$.
\end{lem}
\begin{proof}
  See Thurston~\cite[Corollary~8.10.9]{thurstonN}. A proof of the direction `if $S$ is not a virtual fibre then the entire lamination space is realisable'
  is given more carefully in Canary, Epstein, and Green~\cite[Theorem I.5.3.11]{canary87}: the main step is to check that the space of marked pleated
  surfaces is compact; then one shows that one can realise all finite-leaved maximal laminations and the result follows by taking a sequence of surfaces
  realising finite-leaved approximations to the desired lamination and then passing to a convergent subsequence. The other direction, `if $S$ is a virtual
  fibre then there is a non-realisable lamination', follows from an argument of the following form: consider the pseudo-Anosov mapping class $\phi$ which realises the
  fibration, choose some pleated surface $ S^* $ in the homotopy class of $S$, and look at the iterates $\phi^n(S^*) $; these form a sequence of pleated surfaces
  with no convergent subsequence. The corresponding sequence of \emph{laminations} does have a convergent subsequence, giving an example of a lamination which is
  not realised in the homotopy class of $S$ (see the paragraph subsequent to~\cite[Lemma~I.5.2.16]{canary87}).
\end{proof}

\begin{lem}
  Let $S$ be an essential surface in a geometrically finite hyperbolic $3$-manifold $M$ and suppose $ \Hol(S) $ has no accidental parabolics.
  If the entire geodesic lamination space of $S$ is realisable by pleated surfaces in $M$ in the homotopy class of $S$, then $ \Hol(S) $ is quasi-Fuchsian.
\end{lem}
\begin{proof}
  See Thurston~\cite[Conjecture~8.10.10]{thurstonN}. The result was proved by Bonahon~\cite{bonahon86}, and is a version of the tameness theorem
  for essential surfaces.
\end{proof}

We now state our slightly more general result.
\begin{thm}\label{thm:embedded_geom_finite}
  Let $S$ be an embedded surface in a geometrically finite hyperbolic $3$-manifold $M$ without rank $1$ cusps. Then one of the following two mutually exclusive situations occurs:
  \begin{enumerate}
    \item $\Hol(S)$ is geometrically finite.
    \item $S$ is obtained by adding trivial handles\footnote{A \df{trivial handle} of a surface $ S \subset M $ is an essential (in $S$) subsurface $ E \subset S $ homeomorphic
    to a once-holed torus such that the image of the canonical map $ \pi_1(E) \to \Hol(M) $ is trivial.} to a surface $ S' $  such that $M$ is virtually $S'$-fibred.
  \end{enumerate}
\end{thm}
For the proof, we recall the following definition from Kleineidam and Souto~\cite{kleineidam03} following Masur~\cite{masur86} and Otal~\cite{otal88}.
\begin{defn}
  If $S$ is a compressible surface in a $3$-manifold $M$, then a \df{meridian} is an essential simple closed curve on $S$ which is homotopically trivial
  in $M$ (and hence bounds a compression disc). The set of projective classes of weighted multicurves of
  meridians in $ \PML(S) $ will be denoted $ \mc{M} $ and its closure (in the Hausdorff topology) by $ \mc{M}' $. It is possible to characterise
  the support of laminations in $ \mc{M}' $, see~\cite[Lemma~3.1]{kleineidam03}.

  Restrict to the case that $M$ is a compression body with compression end $S$. The compression body is called \df{small}~\cite[\S 2.2]{kleineidam03} if there is a properly
  embedded essential disc $ (D, \partial D) \subset (M,S) $ such that $ M \setminus D $ is either a thickened surface or a union of two thickened surfaces. Now the
  \df{Masur domain}~\cite[Definition~3.2]{kleineidam03} $ \mc{O}(S) $ of $ S $ is defined in the following way:
  \begin{itemize}
    \item If $M$ is not small, then
          \begin{displaymath}
            \mc{O}(S) = \{ \lambda \in \PML(S) : i(\lambda,\mu) > 0\; \text{for all}\; \mu \in \mc{M}' \}.
          \end{displaymath}
    \item If $M$ is small, then say $ \mu \in \PML(S) $ is \df{$\mc{M}'$-parallel} if there exists $ \nu \in \mc{M}' $ with $ i(\mu,\nu) =0 $. Then
          \begin{displaymath}
            \mc{O}(S) = \{ \lambda \in \PML(S) : i(\lambda,\mu) > 0\; \text{for all}\;\mc{M}'\text{-parallel}\; \mu \in \PML(S) \}.
          \end{displaymath}
  \end{itemize}\
  The \df{parabolic Masur domain}~\cite[Definition~3.11]{kleineidam03} $ \mc{O}_P(S) $ is the set of laminations in $ \mc{O}(S) $ which have positive intersection
  with all simple closed curves on $S$ represented by parabolic elements of $ \Hol(M) $.
\end{defn}

Now the proof of \zcref{thm:embedded_geom_finite} is obtained from a pair of lemmata in the same way as \zcref{thm:bers_essential} above.

\begin{lem}
  Let $S$ be an embedded surface in a geometrically finite hyperbolic $3$-manifold $M$ without rank $1$ cusps such that $ \Hol(M)$ is not free.
  Then either $M$ is virtually $S'$-fibred for some essential surface $S'$ with $ \Hol(S')=\Hol(S) $, or the entire space $ \mc{O}_P(S) $ is realisable by pleated surfaces in $M$ in the homotopy class of $S$.
\end{lem}
\begin{proof}
  We introduce a new property \eqref{eq:approximability}:
  \begin{equation}\label{eq:approximability}\tag{{\tiny\textpmhg{a}\kern-.5ex}}
    \parbox{.75\textwidth}{\centering\slshape There exists a minimal arational lamination $\lambda$ such that no sequence of multicurves on $S$ converging to $\lambda$ can be represented by pleated surfaces in the homotopy class of $S$.}
  \end{equation}

  We now prove three statements to give the result. Namely, the contrapositive of \zcref{subclaim2} together with \zcref{subclaim1} proves one direction, and \zcref{subclaim3} is the other direction.

  \begin{subclaim}\label{subclaim1}
    \eqref{eq:approximability}  $ \implies $ $M$ is virtually a circle bundle with fibre $S'$.
  \end{subclaim}
  This follows the proof of Corollary~8.10.9 of Thurston~\cite{thurstonN}. If there exists such a lamination, then the space of pleated surfaces homotopic to $S$ is not compact. Thus there must be some surface $ S^*$
  which is homotopic to $S$ in more than one way; i.e.\ there is a nontrivial homotopy $ h_t : S \to M $ such that $ h_1 = h_0 \circ \phi $ where $ \phi : S \to S $ is a homotopically nontrivial homeomorphism.
  Up to passing to a finite cover we may assume $ \phi $ is infinite order. The maps $ h_0 $ and $ h_1 $ are therefore conjugate by some $ \beta \in \Hol(M) $. Then $ \langle \beta, \Hol(S) \rangle $ is the
  fundamental group of a circle bundle. \hfill\checkmark

  \begin{subclaim}\label{subclaim2}
    \eqref{eq:approximability} is false $ \implies $ the Masur domain is realised by pleated surfaces homotopic to $S$.
  \end{subclaim}
  This claim is Proposition~4.3 of Kleineidam and Souto~\cite{kleineidam03}.\hfill\checkmark

  \begin{subclaim}\label{subclaim3}
    $M$ is virtually a circle bundle with fibre $S'$ $ \implies $ there is a non-realised lamination in the Masur domain.
  \end{subclaim}
  This claim is proved by the same argument as the final part of the proof of \zcref{lem:classical_compactness}.\hfill\checkmark

  \vspace{1ex}

  The lemma is proved.
\end{proof}

\begin{lem}
  Let $S$ be an embedded compressible surface in a geometrically finite hyperbolic $3$-manifold $M$ without rank $1$ cusps such that $ \Hol(M)$ is not free.
  If the entire space $ \mc{O}_P(S) $ is realisable by pleated surfaces in $M$ in the homotopy class of $S$, then $ \Hol(S) $ is geometrically finite.
\end{lem}
\begin{proof}
  By the tameness theorem~\cite{agol04,calegari06}, $M$ is tame. By Corollary~1.2 of Kleineidam and Souto~\cite{kleineidam03} (they prove one direction and refer to Canary~\cite{canary93}
  for the other), in the presence of tameness, $ \Hol(S) $ is geometrically infinite if and only if there is a lamination in $ \mc{O}_P(S) $ which is not realised.
\end{proof}

\sloppy\printbibliography

\end{document}